\documentclass[times]{article}

\usepackage{amsmath,amssymb}
\usepackage{graphicx}
\usepackage{subfig}
\usepackage{stmaryrd}
\usepackage{bm}
\usepackage{epstopdf}

\usepackage{amsmath}
\usepackage{amssymb}
\usepackage{latexsym}
\usepackage{stmaryrd}
\usepackage{booktabs}
\usepackage{multirow} 
\usepackage{algorithm}
\usepackage{algorithmicx}
\usepackage{algpseudocode}
\usepackage{tabularx}
\usepackage{array}
\usepackage{hyphenat}
\usepackage{float}

\usepackage[authoryear,round]{natbib}
\usepackage[colorlinks,citecolor=blue,linkcolor=blue,urlcolor=blue]{hyperref}
\usepackage{doi}
\DeclareMathOperator{\buffer}{buffer}
\DeclareMathOperator{\opt}{opt}

\DeclareMathOperator{\adaptive}{adaptive}
\DeclareMathOperator{\uniform}{uniform}

\DeclareMathOperator{\refine}{ref}

\DeclareMathOperator{\supp}{supp}
\DeclareMathOperator{\trunc}{trunc}
\DeclareMathOperator{\sign}{sign}

\DeclareMathOperator{\scale}{sc}
\DeclareMathOperator{\shift}{sh}
\DeclareMathOperator{\thres}{th}
\DeclareMathOperator{\SIMP}{}
\DeclareMathOperator{\target}{target}

\DeclareMathOperator{\Lin}{Lin}
\DeclareMathOperator{\Nitsche}{Nitsche}
\DeclareMathOperator{\Ghost}{Ghost}
\DeclareMathOperator{\Spring}{Spring}

\DeclareMathOperator{\cut}{cut}

\DeclareMathOperator{\final}{final}
\DeclareMathOperator{\LoN}{LoN}

\DeclareMathOperator{\xfem}{xfem}
\DeclareMathOperator{\fem}{fem}

\DeclareMathOperator{\ifc}{ifc}
\DeclareMathOperator{\solid}{solid}
\DeclareMathOperator{\void}{void}
\DeclareMathOperator{\init}{init}
\DeclareMathOperator{\current}{current}
\DeclareMathOperator{\bw}{bw}

\DeclareMathOperator{\spring}{S}

\DeclareMathOperator{\up}{up}
\DeclareMathOperator{\low}{low}
\DeclareMathOperator{\vt}{vt}

\DeclareMathOperator{\peak}{peak}

\title{Adaptive level set topology optimization using hierarchical B-splines 
}

\author{L. No\"el, M. Schmidt, C. Messe, J.~A. Evans, K. Maute \\[12pt]
Ann and H. J. Smead Department of Aerospace Engineering Sciences,\\ 
University of Colorado at Boulder,\\
3775 Discovery Dr, Boulder, CO 80303, USA\\
Corresponding author: maute@colorado.edu}

\begin{document}

\maketitle

\begin{abstract}
This paper presents an adaptive discretization strategy for level set topology optimization of structures based on hierarchical B-splines. This work focuses on the influence of the discretization approach and the adaptation strategy on the optimization results and computational cost. The geometry of the design is represented implicitly by the iso-contour of a level set function. An immersed boundary technique, here the extended finite element method, is used to predict the structural response. Both the level set function and the state variable fields are discretized by hierarchical B-splines. While only first-order B-splines are used for the state variable fields, up to third order B-splines are considered for discretizing the level set function. The discretizations of the level set function and the state variable fields are locally refined along the material interfaces and selectively coarsened within the bulk phases. For locally refined/coarsened meshes, truncated B-splines are considered. The properties of the proposed mesh adaptation strategy are studied for level set topology optimization where either the initial design is comprised of a uniform array of inclusions or inclusions are generated during the optimization process. Numerical studies employing static linear elastic material/void problems in 2D and 3D demonstrate the ability of the proposed method to start from a coarse mesh and converge to designs with complex geometries and fine features, reducing the overall computational cost. Comparing optimization results for different B-spline orders suggests that higher order basis functions promote the development of smooth designs and suppress the emergence of small features, without providing an explicit feature size control. A distinct advantage of cubic over quadratic B-splines is not observed.
\end{abstract}

\section{Introduction}\label{intro}

Following the seminal work of \cite{BENDSOE1988}, topology optimization has become, over the last decades, a reliable and efficient design tool in various application fields, see \cite{SIGMUND2013} and \cite{DEATON2014}. In general, an optimization problem is formulated to find the optimal material distribution, within a design domain, that maximizes a given objective while satisfying some given constraints on the geometry and/or the physical response. In specific regions of the design domains, the state variable fields may exhibit large gradients or discontinuities. Additionally, the design variable fields may have to resolve complex shapes, intricate material arrangements, or extremely thin features. In these cases, an increased mesh resolution is necessary to achieve a sufficient accuracy of both the physical responses and the geometry description. Most topology optimization approaches operate on uniformly refined meshes, for which extreme resolution requirements can lead to a drastic increase in computational cost since the number of design and state degrees of freedom (DOFs) is increased. Hence there is a need for efficient adaptive strategies in topology optimization. 

In the 1980s, \cite{BENNET1985} and \cite{KIKUCHI1986} explored mesh adaptation in combination with shape optimization. Their objective was to avoid the distortions of the finite element mesh arising from shape modifications and thus to control the quality of the finite element solution. Later \cite{MAUTE1995, MAUTE1998} and \cite{SCHLEUPEN2000} performed topology optimization and used separate models for the design and the analysis to allow for more flexibility in the optimization process. In these works, the effective design space is adapted at each optimization step with respect to the current material distribution, to both decrease the number of finite elements in the analysis model and the number of design variables used for the geometry description. Ever since, adaptive topology optimization has been an active research topic, as it can provide both precise geometry representations and accurate mechanical responses at a reduced computational cost, substantially speeding up the design process.

Working with density-based approaches and solving compliance minimization problems, several adaptive strategies for topology optimization have been proposed. \cite{COSTA2003} carried out a refinement only approach based on a density criterion. After a given number of optimization steps, the material and boundary elements are refined using \textit{h}-adaptivity, i.e., the elements are split into smaller ones, see \cite{YSERENTANT1986, KRYSL2003}. As no coarsening is introduced in the void regions and as no regularization is used, the optimization leads to mesh dependent designs that differ from the ones obtained with uniformly refined meshes. \cite{STAINKO2006} developed a refinement criterion based on a regularization filter indicator that locates the material interface. In this case, \textit{h}-refinement is only applied around the interface and a reduced number of elements are generated during the optimization process. Extending the previous approaches, \cite{WANG2010} introduced mesh coarsening in the void regions, allowing for further cost reduction and achieving designs that only slightly differs from uniform ones. A similar approach is adopted in \cite{NANA2016} using unstructured meshes and in \cite{NGUYEN2017} using polygonal elements. Also exploiting \textit{h}-refinement, \cite{BRUGGI2011} trigger mesh adaptivity based on error estimators; one related to the geometric error and one to the analysis error.

In the previous approaches, the geometry description and the analysis are strongly coupled, as both the state and the design variables are defined on the same mesh. Following the initial idea by \cite{MAUTE1995}, \cite{GUEST2010} exploited the separation of the geometry and the analysis fields using Heaviside projection methods. After projection, multiple elements are influenced by a single design variable, introducing a redundancy in the design. The latter can be exploited to reduce the number of design variables, and therefore enhance computational performance. \cite{WANG2013, WANG2014} also refined the state and design fields separately by using two distinct meshes adapted by \textit{h}-refinement and resorting to independent geometric and analysis error estimator criteria. 

Later on, adaptive mesh refinement strategies were extended to other optimization techniques, see for example \cite{WALLIN2012} in the context of phase-fields and \cite{PANESAR2017} for Bi-directional Evolutionary Structural Optimization (BESO). Adaptive mesh refinement has also been used to address stress-based optimization problems. \cite{SALAZAR2018} used a density-based approach and adapted the mesh following a stress error estimator, while \cite{ZHANG2018} used moving morphable voids described explicitly by B-splines and refined regions located around the design boundaries.  

So far, most topology optimization approaches have been implemented within the framework of classical finite element, i.e., relying on low order Lagrange interpolation functions. Nonetheless, using B-splines or NURBS as basis functions has become an increasingly popular approach, along with the development of isogeometric analysis (IGA), see \cite{HUGHES2005, COTTRELL2009}. B-splines are used to both describe the geometry of a structure and to solve for its mechanical responses in IGA. From a geometry point of view, using B-splines and NURBS facilitates compatibility with Computer Aided Design (CAD) software, that are based on this type of basis functions. From an analysis point of view, using smooth and higher order bases, such as quadratic and cubic B-splines, leads to more accurate structural responses per DOF than classical $C^0$ finite element approaches, see \cite{HUGHES2008, EVANS2009, HUGHES2014}. For topology optimization, B-splines can offer additional advantages. They tend to promote smoother designs, prevent the development of spurious features and limit the need for filtering techniques. Several works successfully apply topology optimization in combination with IGA, see for example \cite{DEDE2012} for a phase field approach, \cite{QIAN2013} for a density approach, \cite{WANG2016} or \cite{JAHANGIRY2017} for a level set approach, \cite{WANG2017} for lattice structure designs or \cite{LIEU2017} for multi-material designs. A comprehensive review is presented in \cite{WANGY2018}.

Classical B-splines approaches do not allow for local refinement, as their tensor product structure inherently enforces global refinement. However, in solving topology optimization problems, local mesh adaptivity is key to achieve both a precise description of the design boundaries and accurate mechanical response computations while maintaining a reasonable computational cost. Hierarchical B-splines (HB-splines), as introduced by~\cite{FORSEY1988}, naturally support local refinement, but do not form a partition of unity (PU), a beneficial property for numerical purposes. Therefore, the concept was extended to so-called truncated hierarchical B-splines (THB-splines) to recover the PU and other advantageous properties, see \cite{GIANNELLI2012}. Local adaptive mesh refinement relying on HB-splines or THB-splines was implemented by \cite{SCHILLINGER2012} and \cite{GARAU2018}.

Recently research efforts have been dedicated to the exploitation of B-spline refinement in conjunction with level set topology optimization. Solving compliance minimization problems, \cite{BANDARA2016} used B-spline shape functions for both the geometry representation and the analysis with immersed boundary techniques. They proposed a global mesh refinement and coarsening strategy based on Catmull-Clark subdivision surfaces. \cite{WANG2019} focused on cellular structures and used B-spline shape functions to represent the geometry of each representative cell. In their approach, local refinement, i.e., cell subdivision, is achieved by knot insertion. 

This paper proposes an adaptive mesh refinement strategy using HB-splines to perform level set topology optimization. We use a level set function (LSF) to describe the design geometry implicitly. The structural analysis relies on an immersed boundary technique to predict the system responses, here the eXtended Finite Element Method (XFEM) with a generalized Heaviside enrichment strategy. We discretize both the design and the state variable fields, i.e., the level set and the displacement fields, with HB-splines. Contrary to most approaches relying on lower order basis functions, we interpolate the design variables with higher order B-splines. For simplicity, the interpolation of the displacement field is restricted to first order functions. For local refinement we consider truncated B-splines. Although not widely used in the literature so far, THB-splines satisfy the PU and form a convex hull, allowing us to conveniently impose bounds on the design variables. Refinement is triggered according to a user-defined criterion, here the location of the design boundaries, and the mesh is refined along the interfaces and coarsened in the solid and void phases. As the design and state variable fields are defined on the same mesh, they are refined or coarsened simultaneously, which allows for a sufficiently accurate resolution of both the geometry and the physical response. 

Two approaches to handle the design space in level set topology optimization are considered: (i) seeding an initial hole pattern and (ii) nucleating holes during the optimization process using a combined level set/density approach. We solve the optimization problems with mathematical programming techniques and in particular the Globally Convergent Method of Moving Asymptotes (GCMMA), see \cite{SVANBERG2002}. The required sensitivity analysis is carried out with an adjoint formulation. 

The ability of the proposed approach to generate complex geometries at a reduced computational cost starting from initial coarse meshes is assessed with two- and three-dimensional topology optimization problems. By varying the adaptive strategy and the underlying B-spline discretization, we can characterize the influence of the mesh adaptivity on the optimization results and the computational cost. Our numerical study results in interesting findings. First the mesh adaptivity strategy influences the optimization process. Finer initial meshes and more frequent refinement operations lead to the development of more complex geometries exhibiting thin structural members. The properties of B-splines have an impact on the generated designs and an educated choice of the B-splines features can lead to advantageous behaviors for topology optimization. Higher order B-splines promote smoother geometries and tend to eliminate small spurious features from the design. Truncation allows us to conveniently handle enforcement of bounds on the design variables. The numerical examples also reveal that hole nucleation is crucial for the computational efficiency of the adaptive strategy, as it allows for starting the optimization process on coarser meshes.

The remainder of the paper is organized as follows. Section~\ref{GeoDescription} presents the level set description of the geometry using B-splines. Section~\ref{Hole} discusses the strategies adopted to achieve a sufficiently rich design space working with level set topology optimization. Section~\ref{meshAdaptivity} details the mesh refinement and coarsening strategies, while Section~\ref{HBSplineBasis} focuses on the construction of hierarchical B-spline bases. The structural analysis is described in Section~\ref{StructAnalysis}; governing equations are detailed and the basic principles of XFEM are recalled. Section~\ref{optProblem} is devoted to the optimization problem formulation, the used regularization schemes and the sensitivity analysis implementation. Section~\ref{numExample} studies the proposed approach with 2D and 3D solid/void design examples. Finally Section~\ref{Conclusions} summarizes the work and draws some conclusions.

\section{Geometry description}\label{GeoDescription}

Since the development of the level set method (LSM) by \cite{OSHER1988}, it has been extensively used in combination with shape and topology optimization to implicitly describe the geometry of the design, see \cite{VANDIJK2013} and~\cite{SIGMUND2013}. 

\begin{figure}[ht]\center
\includegraphics[scale=0.9]{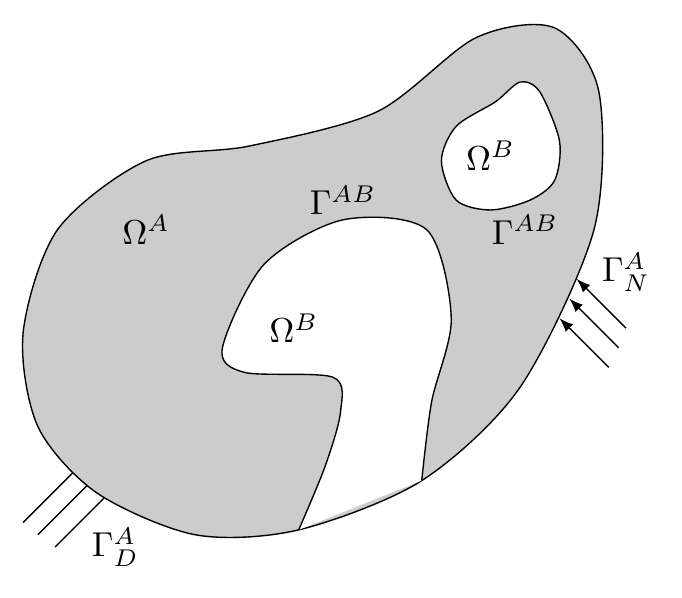}
\caption{Design domain description.}
\label{fig_designDomain}
\end{figure}

A general description of the two-phase problems addressed in this paper is given in Fig.~\ref{fig_designDomain}. Two material phases, $A$ and $B$, are distributed over a $d$-dimensional design domain $\Omega \subset \mathbb{R}^d$. Material subdomains are non-overlapping, i.e., $\Omega = \Omega^A \cup \Omega^B$. The boundaries of material subdomains $A$ and $B$ are denoted $\partial \Omega^A$ and $\partial \Omega^B$. Dirichlet and Neumann boundary conditions are applied to $\Gamma^I_{D} = \partial \Omega^I \cap \partial \Omega_D$ and $\Gamma^I_{N} = \partial \Omega^I \cap \partial \Omega_N$ respectively, with $I = A, B$. The material phases are separated by an interface $\Gamma^{AB} = \partial \Omega^A \cap \partial \Omega^B$.

A LSF is used to describe the material distribution within the design domain $\Omega$. At a given point $\mathbf{x}$ in space, the design geometry is defined as:
\begin{equation}
\renewcommand*{\arraystretch}{1.25}
\begin{array}{llll}
\phi(\mathbf{x}) & < & 0, & \forall \mathbf{x} \in \Omega^{A},\\
\phi(\mathbf{x}) & > & 0, & \forall \mathbf{x} \in \Omega^{B},\\
\phi(\mathbf{x}) & = & 0, & \forall \mathbf{x} \in \Gamma^{AB}.
\end{array}
\end{equation}

The LSF is interpolated using B-spline shape functions $B_i(\mathbf{x})$ as:
\begin{equation}
\phi^h(\mathbf{x}) = \sum_i B_i(\mathbf{x})\ \phi_i,
\end{equation}
where $\phi_i$ are the B-spline coefficients. In this work, the B-spline coefficients and corresponding shape functions are used to evaluated nodal level set values on the integration mesh, that is the mesh over which the weak form of the governing equations are integrated. These nodal values are used to interpolate the level set field linearly within an integration element. The linear approximation within an element simplifies the construction of the intersection of the LSF with the element edges and does not affect the accuracy of the analysis, as linear interpolation functions are used for the state variable field. Further explanations about the level set field interpolation on the analysis mesh are provided in Section~\ref{HBSplineBasis}.

Contrary to classical approaches by \cite{WANG2003} or~\cite{ALLAIRE2004}, where the level set is updated in the optimization process by solving Hamilton-Jacobi-type equations, the coefficients of the LSF are here defined as explicit functions of the design variables. They are updated by mathematical programming techniques driven by shape sensitivity computations.

It should be noted that, in this work, no filter is used to widen the zone of influence of the design variables and thus enhance the convergence of the optimization problem, as proposed in \cite{KREISSL2012}. Proceeding this way, the B-spline basis function support is not altered and its influence on the optimization process and the obtained designs can be emphasized and compared for different B-Spline orders.

\section{Seeding of inclusions}\label{Hole}
The design updates in level set topology optimization is solely driven by shape sensitivities, see \cite{VANDIJK2013}. To allow for a sufficient freedom in the design, level set optimization techniques require either seeding inclusions in the initial design or introducing inclusions during the optimization process. In the context of solid/void problem, the inclusion represents a hole.

Seeding inclusions in the initial design domain is a commonly used strategy and has been successfully implemented in level set topology optimization, see for example \cite{VILLANUEVA2014}. Nonetheless, this strategy leads to several difficulties. First finding a hole pattern that satisfies the design constraints is not a trivial task. Moreover, this strategy requires a sufficiently fine mesh able to resolve both the physics and the design of the initial layout, as the initial hole pattern can easily constitute the most geometrically complex configuration over the entire optimization process. Thus, as will be show in Section~\ref{numExample}, seeding initial holes limits the computational effectiveness of the proposed adaptive discretization strategy. Therefore, allowing for hole seeding during the optimization process proves to be advantageous in the context of mesh adaptivity. 

Topological derivatives constitute a systematic approach to seed new holes in the design domain during the optimization process, see \cite{NOVOTNY2013} for a detailed introduction. The basic concept of evaluating the influence, i.e., the sensitivity, of introducing an infinitesimal hole on the objective and constraint functions was introduced by \cite{ESCHENAUER1994}. A finite size hole is inserted at a location where the topological derivative field associated to an objective or considered cost function is minimal. The shape of the new hole is then optimized along with existing domain boundaries. The reader is referred to \cite{SIGMUND2013} or \cite{MAUTE2017} for further details. Although proven useful, topological derivatives only provide information about where to place an infinitesimal hole and it remains unclear what the size and shape of a finite sized hole should be.

\begin{figure}[h]\center
\includegraphics[scale=1]{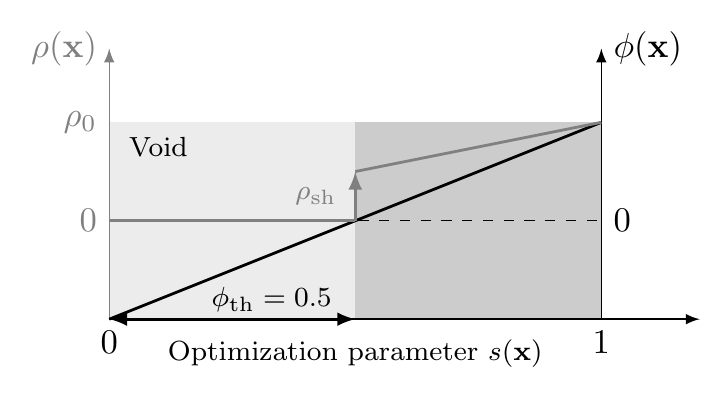}
\caption{Interpolation of the level set $\phi(\mathbf{x})$ and the density $\rho(\mathbf{x})$ for the combined level set/density scheme.}
\label{fig_combinedScheme}
\end{figure}

Alternatively, hole seeding during the optimization process can be achieved by combining level set and density-based techniques, as explored recently in \cite{KANG2013,GEISS2019a,BURMAN2019} or \cite{JANSEN2019}. Following the single-field ap\-proach of \cite{BARRERA2019}, an abstract design variable field, $s(\mathbf{x})$, is introduced with $s \in \mathbb{R}$, $0 \leq s(\mathbf{x}) \leq 1$, to define both the level set field, $\phi$, and spatially variable material properties, such as the density $\rho$ and the Young's modulus $E$, using a Solid Isotropic Material with Penalisation (SIMP) interpolation scheme:
\begin{equation}
\phi( \mathbf{x} ) = \phi_{\scale} \left( \phi_{\thres} - s( \mathbf{x} ) \right),
\end{equation}

\begin{equation}
\rho( \mathbf{x} ) = \left\{
\renewcommand*{\arraystretch}{1.25}
\begin{array}{ll} 
0, & \forall\ s( \mathbf{x} ) < \phi_{\thres},\\
 \rho_{\shift} + \left( \rho_{0} - \rho_{\shift} \right) \frac{\left( s(\mathbf{x}) - \phi_{\thres} \right)}{\left( 1 - \phi_{\thres} \right)}, & \forall\ s( \mathbf{x} ) \geq \phi_{\thres},
\end{array}
\right.
\end{equation}

\begin{equation}
E( \mathbf{x} ) = E_0\  \rho^{\beta_{\SIMP}},
\end{equation}
where $\phi_{\scale}$ is a scaling parameter that accounts for the mesh size $h$ and is set to $\phi_{\scale} \simeq 3h \dots 5h$, $\phi_{\thres}$ is a threshold that defines the void/solid phases in terms of $s(\mathbf{x})$ and $\rho_{\shift}$, a parameter that controls the minimum density in the solid domain. The properties of the bulk material are denoted by $\rho_0$ and $E_0$, and $\beta_{\SIMP}$ is the SIMP exponent. The resulting interpolation scheme is illustrated in Fig.~\ref{fig_combinedScheme}. Throughout the optimization process, the parameter $\phi_{\thres}$ is kept constant and equal to $\phi_{\thres} = 0.5$, while $\rho_{\shift}$ is a continuation parameter and is gradually increased to $\rho_{\shift} = 1$ during the process. To avoid ill-conditioning, its initial value is set to $\rho_{\shift} \simeq 0.1 \dots 0.2$, unless a smaller value is required to satisfy an initial mass or volume constraint. Contrary to classical density-based approaches, the density here is used primarily for hole seeding and for convergence acceleration. A $0-1$ material distribution is achieved through continuation on the parameter $\rho_{\shift}$. A low SIMP exponent, e.g., $\beta_{\SIMP} = 2.0$, is used to reduce the bias of the density method to rapidly separate the material distribution into a solid and void phase. For the proposed mesh adaptation strategy, we have observed that using a low SIMP exponent promotes the formation of fine features, especially when starting from coarse meshes.

In this paper, both approaches, i.e., the seeding of holes in the initial design and the combined level set/density scheme, are considered and compared.

\section{Mesh adaptivity}\label{meshAdaptivity}
This section focuses on the strategies used to perform hierarchical mesh refinement. First, we introduce the notion of a hierarchical mesh. We recall the basic concepts of mesh refinement and detail our implementation of the refinement strategy. Finally, the considered user-defined criteria, used to trigger mesh adaptation, are explained in detail.
 
\subsection{Hierarchical mesh}\label{hierarchicalMesh}
In this work, a regular background tensor grid is used to build hierarchical meshes. The elements of the background grid do not carry any notion of interpolation, i.e., an element of the background grid is not associated with B-spline bases or Lagrange nodes. Starting from a uniform background grid with a refinement level $l^{0}$ and subdividing its elements, elements with a higher refinement level $l > l^{0}$ are created.

\begin{figure}[ht]\center
\includegraphics[width=0.25\columnwidth]{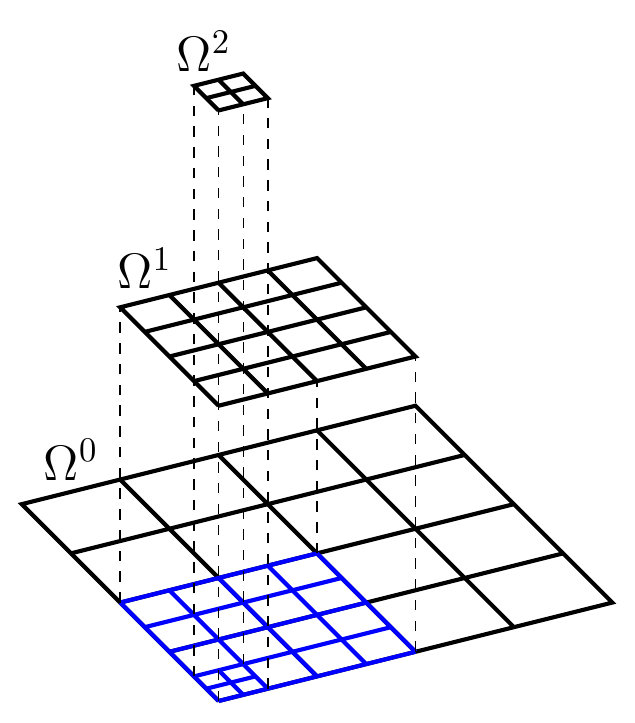}\hspace{0.25cm}
\includegraphics[width=0.175\columnwidth]{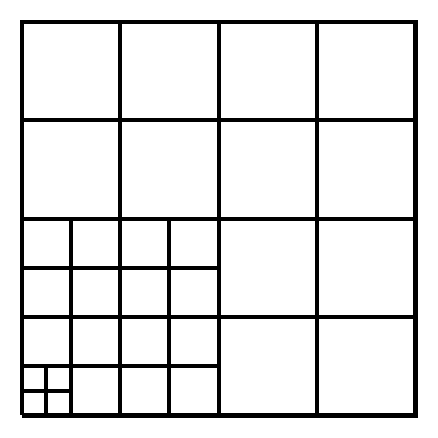}
\caption{Hierarchical refined mesh.}
\label{fig_hierarchicalRefinedMesh}
\end{figure}

Formally to define a hierarchical mesh of depth $n$, a sequence of subdomains $\Omega^{l}$ is introduced:
\begin{equation}
\Omega^{n-1} \subseteq \Omega^{n-2} \subseteq \dots \subseteq \Omega^{0} = \Omega,
\end{equation}
where each subdomain $\Omega^l$ is the subregion of $\Omega$ selected to be refined at the level $l$ and is the union of elements on the level $l-1$. The creation of such a refinement pattern is illustrated in Fig.~\ref{fig_hierarchicalRefinedMesh}.
In Section~\ref{HBSplineBasis}, these subdomains $\Omega^{l}$, refined to a level $l$, are used to create hierarchical B-spline bases to interpolate the design and state variable fields.
 
\subsection{Local refinement strategies}\label{refinementStrategies}
This subsection focuses on the refinement strategies used in this paper and explains their implementation. Before refinement is carried out by element subdivision, elements need to be flagged for adaptation. Elements can be flagged based on one or several user-defined criteria and on additional mesh regularity requirements. 

A complete refinement step proceeds as stated in Algorithm~\ref{alg_refinement}. Initially, elements are flagged for refinement based on user-defined criteria, as further explained in Section~\ref{userDefinedCriteria}. Additionally, elements are flagged within a so-called buffer region, around the previously flagged elements, so that the obtained refined mesh complies with specific regularity requirements, as further detailed in Section~\ref{bufferCriterion}. Finally, a minimum refinement level $l_{\min}$ is enforced and elements with $l< l_{\min}$ are flagged and refined. To increase performance, flagged elements are collected in a refinement queue first and refined afterwards.
\begin{algorithm}[ht]\center
\begin{algorithmic}[1]
\State{Apply refinement strategy in Algorithm~\ref{Alg_refinementStrategy} }
\State{Collect flagged elements into queue for refinement}
\While{change in refinement queue}
              \For{all elements in queue for refinement}
                           \State{Apply refinement buffer in Algorithm~\ref{alg_refinementBuffer}}
                           \State{Collect flagged elements in queue for refinement}
              \EndFor
\EndWhile
\State{Refine elements in refinement queue}
\For{all elements}
\If{element refinement level $l< l_{\min}$}
              \State{Flag element for refinement}
              \State{Refine element in refinement queue}
\EndIf
\EndFor
\end{algorithmic}      
\caption{Refinement algorithm.}
\label{alg_refinement}
\end{algorithm}

In this paper, design and state variables are defined on a common mesh, although the interpolation order may differ. Thus, the discretization of both fields are adapted simultaneously. While the proposed framework allows for separately refining state and design variable fields, this option is not considered here for the sake of simplicity. 

After each adaptation of the mesh, the design variable field used for the geometry description is mapped to the new refined mesh through an $L^{2}$ projection. This projection step requires an additional linear solve for a system which size depends on the number of design DOFs.

\subsubsection{User-defined refinement criteria}\label{userDefinedCriteria}
Refinement can be triggered by one or several user-defined criteria, including geometry criteria or finite element error indicators, although the latter are not considered here. In this paper, refinement is applied after a given number of optimization steps using a purely geometric criterion, i.e., the distance to the solid-void interfaces. 

During each refinement step, an element is flagged for refinement or for keeping its current refinement level by the Algorithm~\ref{Alg_refinementStrategy}. Each element mesh is flagged for refinement depending on its nodal level set values $\phi$. The parameter $\phi_{\bw}$ defines the bandwidth around the interface which is refined. This allows fine-tuning the zone of refinement in addition to the buffer introduced in Section~\ref{bufferCriterion}. Over the course of the optimization process, at the solid-void boundary, the mesh is refined up to a maximum refinement level defined for the interface $l_{\ifc,\max}$. Within the solid phase, the mesh is refined up to a maximum refinement level $l_{\solid,\max}$, that is usually equal or smaller than the interface one, $l_{\solid,\max} \leq l_{\ifc,\max}$. The void phase is never flagged for refinement. However, by applying an initial minimum refinement level, which is lowered to the lowest allowable level of refinement $l_{\min}$ throughout the mesh adaptive process, a coarsening effect is achieved.
\begin{algorithm}[ht]\center
\begin{algorithmic}[1]
\State{Step 0. Set $l_{\ifc}$, $l_{\solid}$, and $l_{\void}$ to an initial mesh refinement level $l_{0}$,
             where $l_{\min} \leq l_{0} \leq \min(l_{\solid,\max}, l_{\ifc,\max})$.}

\State{Step 1. Flag elements}
\For{all elements}
\If{$- \phi_{\bw} \leq \phi \leq + \phi_{\bw}$}
\If{$l_{\current} < l_{\ifc}$}
\State{Refine}
\Else
\State{Keep}
\EndIf
\EndIf

\If{$\phi < -\phi_{\bw}$}
\If{$l_{\current} < l_{\solid}$}
\State{Refine}
\Else
\State{Keep}
\EndIf
\EndIf

\If{$\phi > + \phi_{\bw}$}
\If{$l_{\current} > l_{\void}$}
\State{Do nothing}
\Else
\State{Keep}
\EndIf
\EndIf  
\EndFor

\State{Step 2. Increase/reduce refinement levels}
\State{$l_{\ifc} = \min(l_{\ifc+1}, l_{\ifc,\max})$}
\State{$l_{\solid} = \min(l_{\solid+1}, l_{\solid,\max})$}
\State{$l_{\void} = \max(l_{\void -1}, l_{\min})$}
\end{algorithmic}       
\caption{Refinement strategy.}
\label{Alg_refinementStrategy}
\end{algorithm}
 
In the course of the optimization process, the shape of the solid-void interface changes and the interface may intersect elements previously not refined. In this study, the mesh is adapted to maintain a uniform refinement level for all intersected elements. Although not required by the proposed framework, this strategy is used here as it simplifies the XFEM enrichment, see Section~\ref{XFEM}. The same enrichment algorithm can be used for all intersected elements, as is the case for uniformly refined meshes. To this end, Algorithm~\ref{alg_refinement} is executed, omitting Step 2 of Algorithm~\ref{Alg_refinementStrategy}.

\subsubsection{Mesh regularity requirements}\label{bufferCriterion}
In this paper, we limit the size difference between adjacent elements in the refined mesh, i.e., a refinement of more than a factor 4 in 2D and 8 in 3D is not allowed. Although not mandatory, see studies by \cite{JENSEN2016} and \cite{PANESAR2017}, this rule is adopted in most works dedicated to adaptivity as it promotes accurate analysis results. This requirement can be achieved by enforcing a so-called buffer zone around elements primarily flagged for refinement. The width of the buffer zone $d_{\buffer}$ for a particular flagged element is calculated by multiplying its size with a user-defined buffer parameter $b_{\buffer}$. The width of the buffer zone must be larger than the support size of the considered interpolation functions. Using B-splines, it means that $b_{\buffer} \geq p$, where $p$ is the polynomial degree of the used B-spline basis.

The mesh refinement procedure for building a buffer zone is summarized in Algorithm~\ref{alg_refinementBuffer}. Within a refined mesh, coarse elements are called parents and refined elements, children. The algorithm is applied to each flagged element and starts by determining its parent element. The refinement status of the parent's neighbors, i.e., elements lying within the buffer range of the considered parent, is checked. If these neighbors are neither refined nor flagged for refinement, the distance $d_{\max}$ between the considered parent element and its neighbors is calculated. If the distance $d_{\max}$ is smaller than the buffer size $d_{\buffer}$, the neighbor elements are flagged for refinement. The algorithm is then applied recursively to all newly flagged neighbor elements, until no extra element are flagged. An easy and efficient access to hierarchical mesh information, such as neighborhood relationships, is provided by using a quadtree and an octree data structure in 2D and 3D respectively.
\begin{algorithm}[ht]\center
\begin{algorithmic}[1]
\State{Get parent}
\State{Get parent's neighbors in buffer range $d_{\buffer}$}
\If{neighbor is not refined and not flagged for refinement}
	\State{Calculate distance $d_{\max}$ between the flagged element and the neighbor}
	\If{$d_{\max} < d_{\buffer}$}
		\State{Flag neighbor for refinement}
		\State{Apply refinement buffer Algorithm~\ref{alg_refinementBuffer} for neighbor}
	\EndIf
\EndIf
\end{algorithmic}      
\caption{Refinement buffer algorithm.}
\label{alg_refinementBuffer}
\end{algorithm}
  
\section{Hierarchical B-splines}\label{HBSplineBasis}
This section focuses on hierarchical B-splines for adaptively discretizing design variable and state variable fields. First the basic concepts of B-splines in one and multiple dimensions are recalled. Then the principle of B-spline refinement and the construction of truncated and non-truncated hierarchical B-spline bases are briefly described.

\subsection{Univariate B-spline basis functions}\label{uniBSpline}
Starting from a knot vector $\Xi = \{ \xi_{1},\xi_{2}, \dots ,\xi_{n+p+1} \}$, for which $\xi \in \mathbb{R}$ and $\xi_{1} \leq \xi_{2} \leq \dots \leq \xi_{n+p+1}$, a univariate B-spline basis function $N_{i,p}(\xi)$ of degree $p$ is constructed recursively starting from the piecewise constant basis function:
\begin{equation}
N_{i,0}(\xi) =
\begin{cases}
1, & \text{if}\ \xi_{i} \leq \xi \leq\xi_{i+1},\\
0, & \text{otherwise},
\end{cases}
\end{equation}
and using the Cox de Boor recursion formula (\cite{DEBOOR1972}) for higher degrees, $p > 0$:
\begin{equation}
\begin{split}
N_{i,p}(\xi) = & \frac{\xi - \xi_{i}}{\xi_{i+p} - \xi{i}}\ N_{i,p-1}(\xi)\\
& \hspace{1.25cm} + \frac{\xi_{i+p+1} - \xi}{\xi_{i+p+1} - \xi{i+1}}\ N_{i+1,p-1}(\xi).
\end{split}
\end{equation}
A knot is said to have a multiplicity $k$ if it is repeated in the knot vector. The corresponding B-spline basis exhibits a $C^{p-k}$ continuity at that specific knot, while it is $C^{\infty}$ in between the knots. A knot span is defined as the half open interval $[\xi_i, \xi_{i+1})$. Within this context, an element is defined as a nonempty knot span. 

\subsection{Tensor-product B-spline basis functions}\label{tensorBSplines}
Tensor-product B-spline basis functions $B_{i,p}(\xi)$ are obtained by applying the tensor product to univariate B-spline basis functions. Denoting the parametric space dimension as $d_{p}$, a tensor-product B-spline basis is constructed starting from $d_{p}$ knot vectors $\Xi^{m} = \{\xi_{1}^{m},\xi_{2}^{m},\-\dots,\-\xi_{n_{m}+p_{m}+1}^{m}\}$ with $p_{m}$ donating the polynomial degree and $n_{m}$ the number of basis functions in the parametric direction $m = 1, \dots, d_{p}$. A tensor-product B-spline basis function is generated from $d_{p}$ univariate B-splines $N_{i_{m},p_{m}}^{m}(\xi^{m})$ in each parametric direction $m$ as:
\begin{equation}
B_{\mathbf{i},\mathbf{p}}(\boldsymbol{\xi}) = \prod_{m=1}^{d_{p}} N_{i_{m},p_{m}}^{m}(\xi^{m}),
\label{Eq_multiTensorProd}
\end{equation}
where the position in the tensor product structure is given by the index $\mathbf{i} =\{ i_{1}, \dots, i_{d_{p}} \}$ and $\mathbf{p} = \{ p_{1}, \dots, p_{d_{p}} \}$ is the polynomial degree. Similar to the univariate case, an element is defined as the tensor product of $d_{p}$ nonempty knot spans. Additionally, a B-spline space $\mathcal{V}$ is defined as the span of B-spline basis functions.

\subsection{B-spline refinement}\label{BSplineRefinement}
\begin{figure}[h]\centering
\includegraphics[width=0.45\columnwidth]{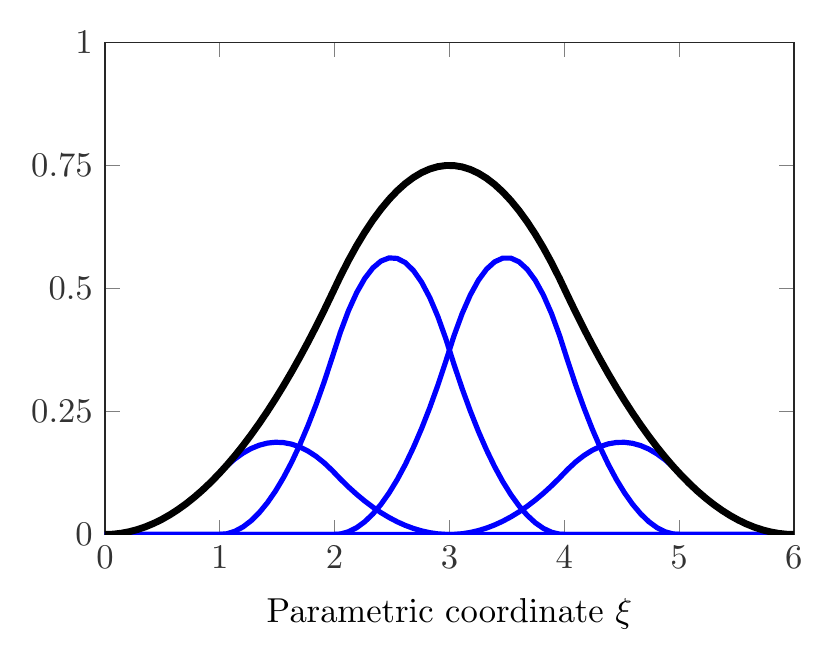}
\caption{Subdivision of a quadratic B-spline basis function (black) into $p+2$ contracted B-spline basis functions of half the knot span width (blue).}
\label{fig_BSplineSubdivision}
\end{figure}

Hierarchical refinement of uniform B-splines can be a\-chieved by subdivision. A univariate B-spline basis function can be expressed as a linear combination of $p+2$ contracted, translated and scaled copies of itself:
\begin{equation}
N_{p}(\xi) = 2^{-p}\sum_{j=0}^{p+1}\binom{p+1}{j}N_{p}(2\xi-j),
\end{equation}
where the binomial coefficient is defined as:
\begin{equation}
\binom{p+1}{j} = \frac{(p+1)!}{j!(p+1-j)!}.
\end{equation}
Figure~\ref{fig_BSplineSubdivision} shows the refinement of a quadratic univariate B-spline basis function obtained by subdivision.

The extension of the subdivision property in Eq.~\eqref{Eq_multiTensorProd} to tensor-product B-spline basis functions $B_{p}$ is straightforward as they exhibit a tensor product structure:
\begin{equation}
B_{\mathbf{p}}(\boldsymbol{\xi}) = \sum_{\mathbf{j}} \left( \prod_{m=1}^{d} 2^{-p_{m}}\ \binom{p_{m}+1}{j_{m}}\ N_{p_{m}}(2\xi^{m}-j_{m}) \right),
\end{equation}
where the indices $\mathbf{j} = \{i_{1}, \dots, i_{d_{p}}\}$ indicate the position in the tensor product structure.
 
\subsection{Hierarchical B-splines}\label{HBSpline}
To build a hierarchical B-spline basis a sequence of tensor product B-spline spaces is introduced:
\begin{equation}
\mathcal{V}^{0} \subset \mathcal{V}^{1} \subset \mathcal{V}^{2} \subset \mathcal{V}^{3} \subset \dots
\end{equation}
Each B-spline space $\mathcal{V}^{l}$ has a corresponding basis $\mathcal{B}^{l}$.

A hierarchical B-spline basis $\mathcal{H}$ can be constructed recursively based on the sequence of B-spline bases $\mathcal{B}^{l}$ that span the domains $\Omega^{l}$. In an initial step, the basis functions defined on the coarsest level, $l=0$, are collected and assigned to $\mathcal{H}^{0}$. The hierarchical B-spline basis $\mathcal{H}^{l+1}$ is constructed by taking the union of all basis functions $B$ in $\mathcal{H}^{l}$ whose support is not fully enclosed in $\Omega^{l+1}$ and all basis functions $B$ in $\mathcal{B}^{l+1}$ whose support lies in $\Omega^{l+1}$. The recursive algorithm reads (\cite{GARAU2018}):
\begin{equation}
\left\{
\begin{array}{lll}
\mathcal{H}^0      & =: & \mathcal{B}^0\\[5pt]
\mathcal{H}^{l+1} & =: & \{ B \in \mathcal{H}^{l}\ |\ \supp(B) \not\subseteq \Omega^{l+1}\}\ \cup \\[5pt]
&&\{ B \in \mathcal{B}^{l+1}\ |\ \supp(B) \subseteq \Omega^{l+1}\},\\[5pt]
&&\mbox{for}\ l = 0, \dots, n-2,\\
\end{array}
\right.
\label{Eq_recursiveAlgorithm}
\end{equation}
where the index $l$ gives the level of refinement. Basis functions collected in $\mathcal{H}$, where $\mathcal{H} = \mathcal{H}^{n-1}$, are called active, while basis functions in $\mathcal{B}^{l}$ not present in $\mathcal{H}$ are said to be inactive.

\begin{figure}\center
\includegraphics[width=0.5\columnwidth]{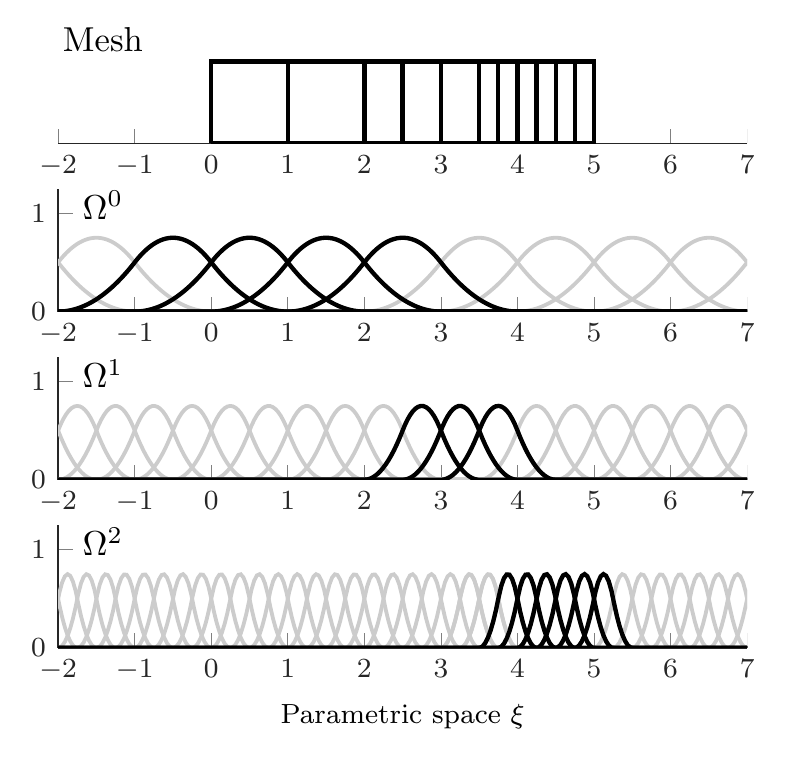}
\caption{Three levels of refinement on a given mesh. The refined mesh is depicted in the top figure. The three underlying figures present the three levels of refinement considered here. Within these three levels, the active B-spline basis $\mathcal{H}$ is displayed in black, while the inactive B-spline basis functions are displayed in gray.}
\label{fig_MeshRefinement1D}
\end{figure}

A hierarchical B-spline basis $\mathcal{H}$ is illustrated for a one dimensional example in Fig.~\ref{fig_MeshRefinement1D}. The top row shows a one dimensional hierarchical refined mesh. Following Algorithm~\ref{Eq_recursiveAlgorithm}, a B-spline basis $\mathcal{H}$ is created through an initialization step with all bases in the subdomain $\Omega^{0}$ refined to a level $l=0$. Recursively all bases in the subdomain $\Omega^{l+1}$ with higher refinement level $l+1$ are added, while existing basis functions of level $l$ fully enclosed in $\Omega^{l+1}$ are discarded.

Contrary to \textit{h}-refinement in classical finite element where extra treatments, such as the introduction of multi-point constraints, are necessary, hanging nodes in hierarchical refined meshes are naturally handled by the B-spline basis.

\subsection{Truncated B-splines}\label{TBSpline}
A major drawback of hierarchical B-spline bases is the loss of the PU property. The truncated hierarchical B-spline basis constitutes an extension of the hierarchical B-spline basis, aiming at recovering the PU principle and at reducing the number of overlapping functions on adjacent hierarchical levels, see \cite{GIANNELLI2012}. Considering a basis function $B^{l}$, part of $\mathcal{B}^{l}$ and defined on the domain $\Omega^{l}$, its representation in terms of the finer basis of level $l+1$ is given as:
\begin{equation}
B^{l} = \sum_{B^{l+1}\ \in\ \mathcal{B}^{l+1}} c_{B^{l+1}}^{l+1} \left( B^{l} \right)\ B^{l+1},
\end{equation}
where $c_{B^{l+1}}^{l+1}$ is the coefficient associated to a basis function $B^{l+1}$.

The truncation of this basis function $B^{l}$, whose support overlaps with the support of finer basis functions $B^{l+1}$, part of $\mathcal{B}^{l+1}$ and defined on $\Omega^{l+1}$, is expressed as (\cite{GIANNELLI2012, GARAU2018}):
\begin{equation}
\begin{array}{lll}
\trunc^{l+1}(B^{l}) & = &\displaystyle \sum_{
\begin{array}{c}
 \scalebox{0.75}{$B^{l+1} \in \mathcal{B}^{l+1},$}\\ 
 \scalebox{0.75}{$\supp(B^{l+1}) \not\subseteq \Omega^{l+1}$}
\end{array}} c_{B^{l+1}}^{l+1} \left( B^{l} \right)\ B^{l+1},\\[40pt]
& = & \displaystyle B^{l} - \sum_{\supp(B^{l+1}) \subseteq \Omega^{l+1}} c_{B^{l+1}}^{l+1} \left( B^{l} \right)\ B^{l+1}.\\
\end{array}
\label{Eq_truncatedBSpline}
\end{equation}

Similar to the creation of a hierarchical B-spline basis $\mathcal{H}$, a truncated hierarchical B-spline basis $\mathcal{T}$ is constructed recursively, but by additionally applying the truncation, see Eq.~\eqref{Eq_truncatedBSpline}, at each iteration (\cite{GIANNELLI2012, GARAU2018}):
\begin{equation}
\left\{
\begin{array}{lll}
\mathcal{T}^0      & =: & \mathcal{B}^0\\[5pt]
\mathcal{T}^{l+1} & =: & \{ \trunc^{l+1}(B) \ |\ B \ \mbox{in} \ \mathcal{T}^{l} \wedge \supp(B) \not\subseteq \Omega^{l+1}\}\\[5pt]
&& \cup\ \{ B \in \mathcal{B}^{l+1} \ |\ \supp(B) \subseteq \Omega^{l+1}\}, \\[5pt]
&& \mbox{for}\ l = 0, \dots, n-2.
\end{array}
\right.
\label{Eq_truncatedRecursiveAlgorithm}
\end{equation}

The truncated basis $\mathcal{T}$ spans the same space as the non-truncated basis $\mathcal{H}$, but its underlying functions exhibit the following properties. They are defined on smaller supports which results in a sparser stiffness matrices in a finite element analysis. The PU property is restored which allows for the imposition of bounds on variables. In addition, truncated bases admit a strong stability property (\cite{GIANNELLI2014}), which enables the construction of optimal multi-level solvers  (\cite{HOFREITHER2016}). However, this properties is not exploited in this paper, as the displacements are interpolated by linear B-splines. 

\begin{figure*}[ht]\center
\includegraphics[width=.45\textwidth]{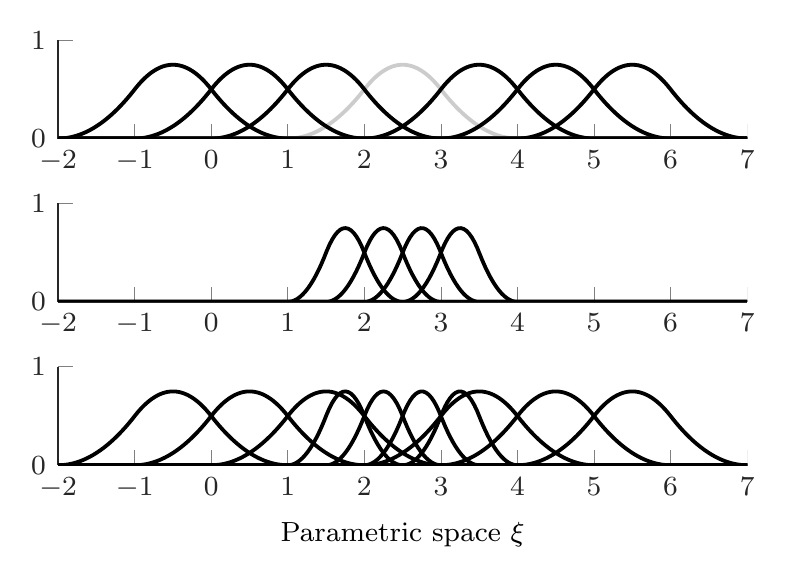}\hspace{0.5cm}
\includegraphics[width=.45\textwidth]{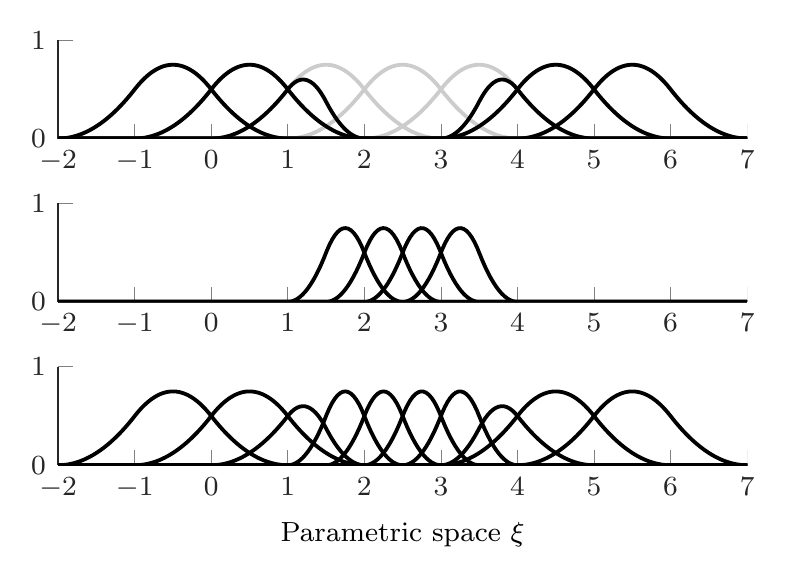}
\caption{Comparison of univariate HB-spline (left) and THB-spline basis functions (right). The first and second levels correspond to $\Omega^0$ and $\Omega^1$ respectively, while the bottom level represents the combination of the functions on these two levels.}
\label{fig_truncation1D}
\end{figure*}

\begin{figure*}[ht]\center
\includegraphics[width=.25\textwidth]{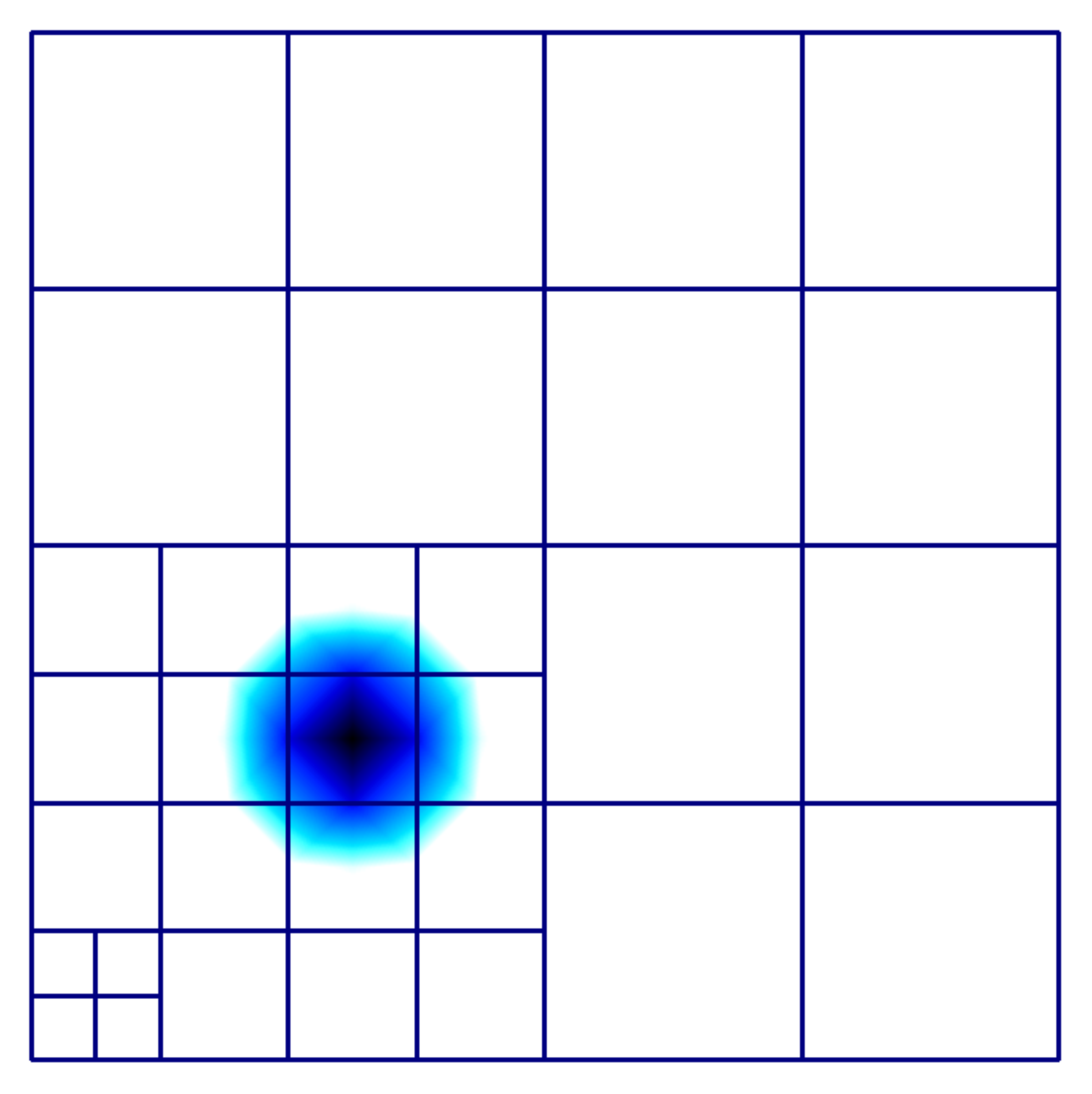} \hspace{0.5cm}
\includegraphics[width=.25\textwidth]{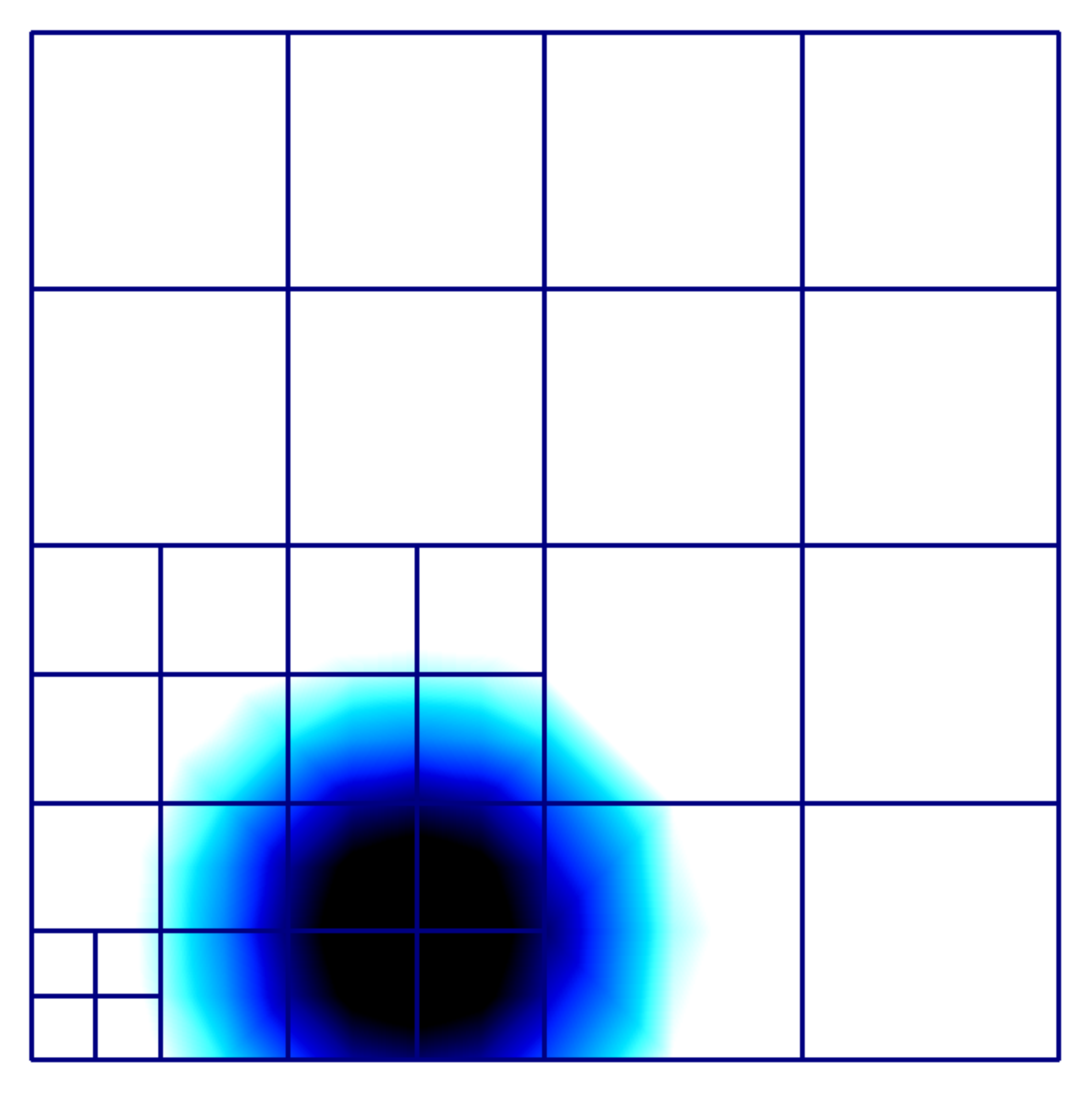} \hspace{0.5cm}
\includegraphics[width=.25\textwidth]{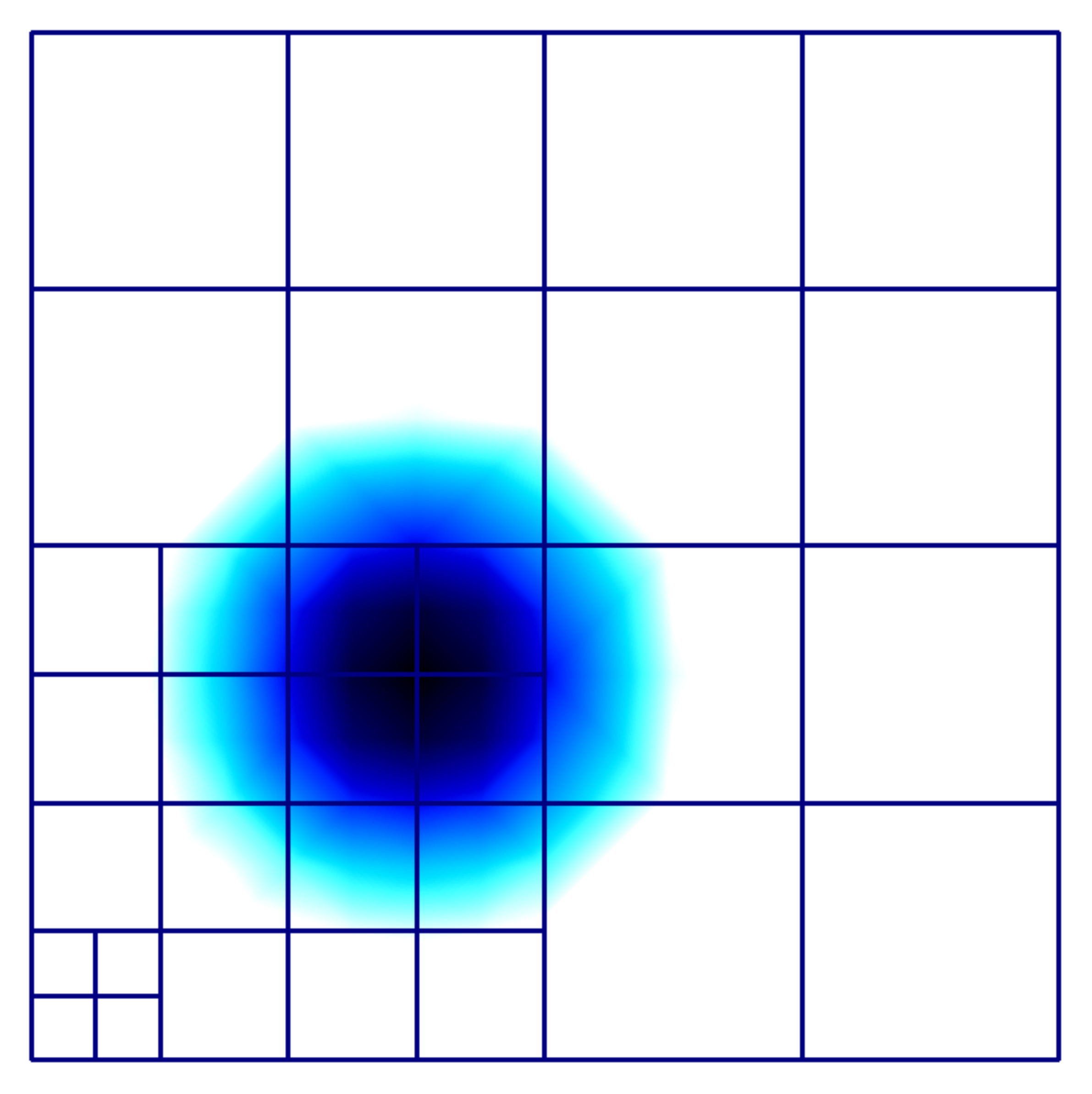}\\ \vspace{0.7cm}
\includegraphics[width=.25\textwidth]{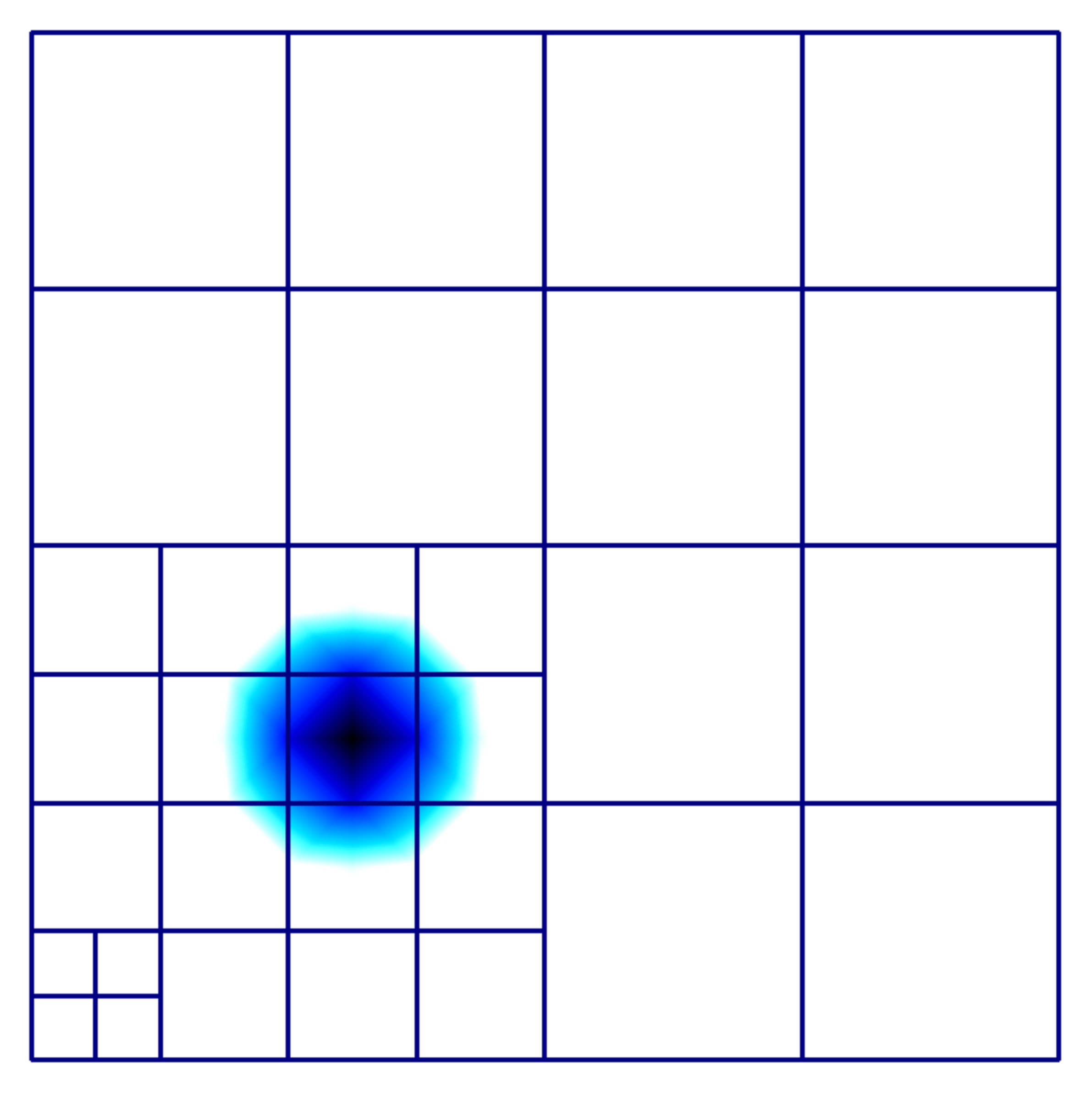} \hspace{0.5cm}
\includegraphics[width=.25\textwidth]{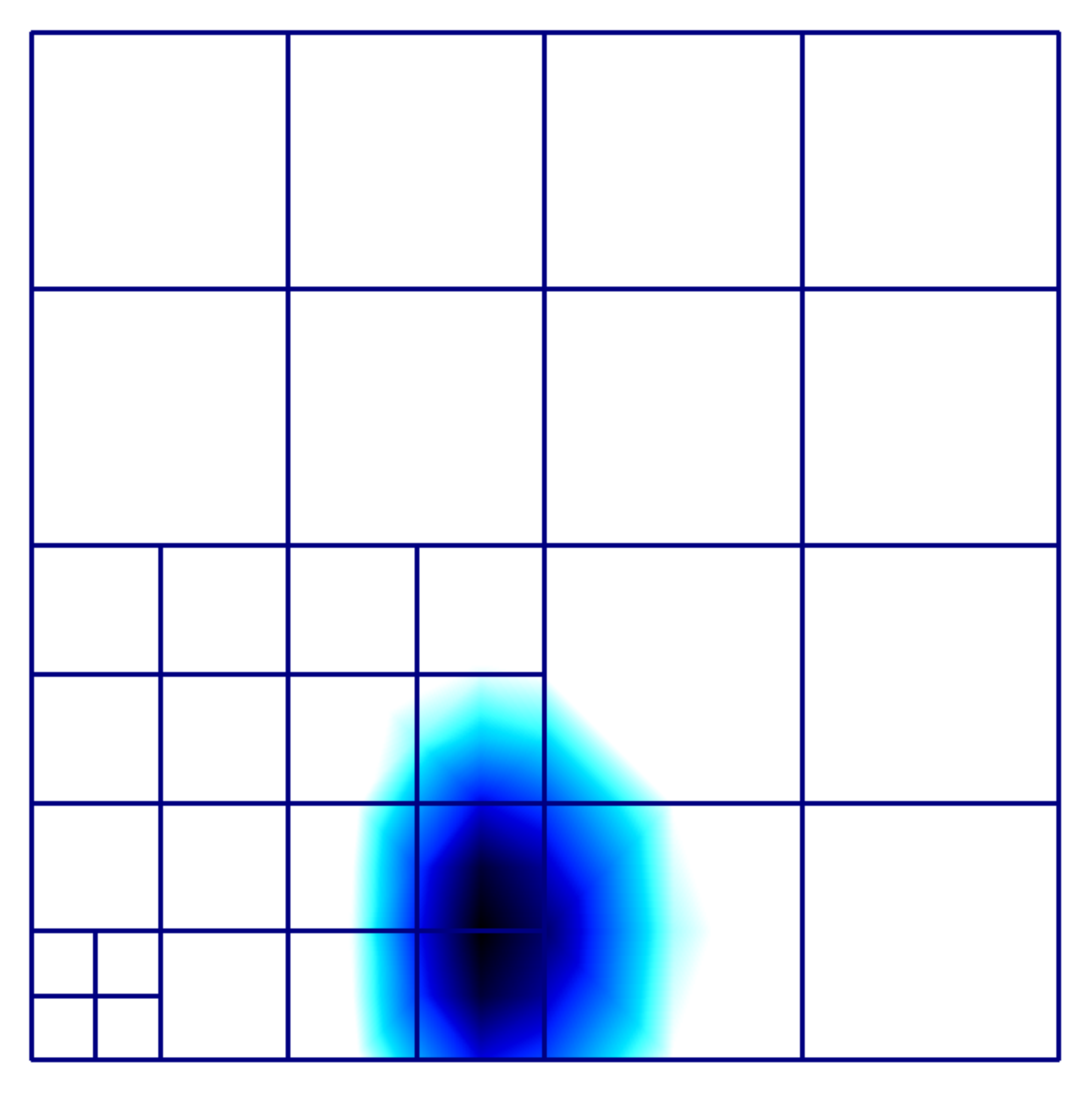} \hspace{0.5cm}
\includegraphics[width=.25\textwidth]{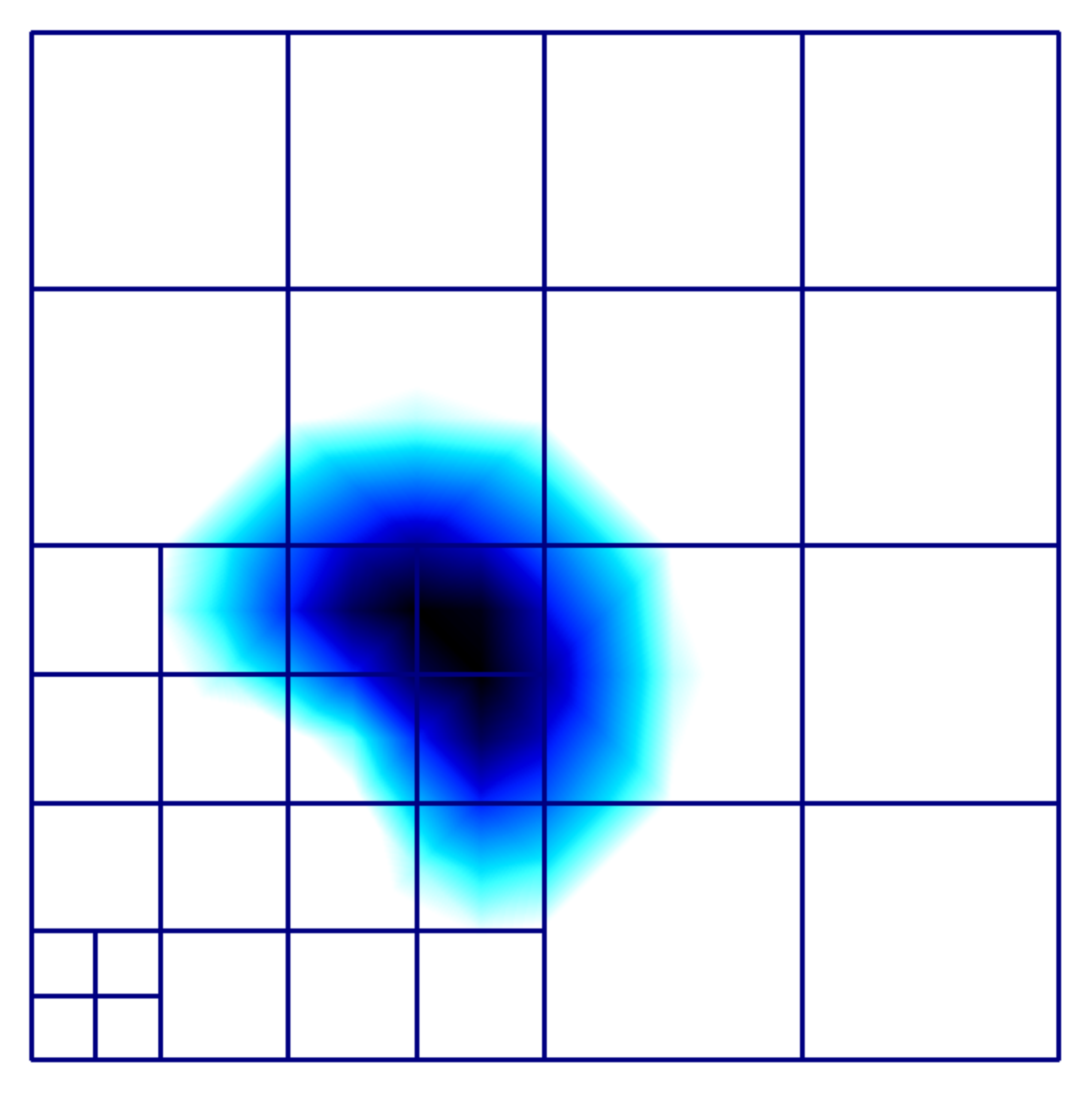}\\
\caption{Comparison of multivariate HB-spline (top) and THB-spline (bottom) supports. While the non-truncated HB-spline basis functions show a uniform support zone, the support zone of THB-spline basis functions exhibit a characteristic kidney shape.}
\label{fig_truncatedKidney}
\end{figure*}

The effect of the truncation is illustrated in Fig.~\ref{fig_truncation1D}. A univariate truncated and non-truncated basis $\mathcal{T}$ and $\mathcal{H}$ are juxtaposed for comparison. The comparison shows the reduced support of the truncated B-splines. The effect of truncation on the support of a multivariate, two dimensional B-spline basis is shown in Fig.~\ref{fig_truncatedKidney}. Truncated multivariate B-splines present a characteristic kidney shaped support zone.

\section{Structural analysis}\label{StructAnalysis}
Whereas the proposed adaptive discretization scheme is not limited to any particular physics problem, in this paper, we consider solid-void problems as illustrated in Fig.~\ref{fig_designDomain}, where the material subdomain $A$ is filled with a linear elastic solid, while the material subdomain $B$ is void. In this section, we present the variational form of the governing equations for the linear elastic analysis model. As the XFEM is used to discretize the state variable fields in the solid domain, the basic principles of the method are briefly recalled.

\subsection{Governing equations}\label{govEq}
We solve for static equilibrium to enforce balance of linear momentum within the solid domain $\Omega^A$. In this work, the total residual $\mathcal{R}$, i.e., the weak form of the governing equations, consists of four terms which are discussed subsequently:
\begin{equation}
\mathcal{R} = \mathcal{R}_{\Lin} + \mathcal{R}_{\Nitsche} + \mathcal{R}_{\Ghost} + \mathcal{R}_{\Spring}.
\label{Eq_RFull}
\end{equation} 

The weak form of the linear elastic governing equation reads:
\begin{equation}
\displaystyle \mathcal{R}_{\Lin} = 
\int_{\Omega^A}\ \delta \boldsymbol{\varepsilon} : \boldsymbol{\sigma}\ d\Omega 
- \int_{\Gamma_N^A} \delta \mathbf{u} \cdot \mathbf{t}_{N}\ d\Gamma,
\label{Eq_RLin}
\end{equation}
where $\mathbf{u}$ and $\delta \mathbf{u}$ are the displacement field and the test function, respectively. Traction forces, $\mathbf{t}_N$, are applied on the Neumann boundary, $\Gamma_N^A$. The Cauchy stress tensor is denoted by $\boldsymbol{\sigma} = \mathbf{D}\ \boldsymbol{\varepsilon}$ and is obtained by multiplication of the infinitesimal strain tensor $\boldsymbol{\varepsilon} = \frac{1}{2} \left( \nabla \mathbf{u} + \nabla \mathbf{u}^T \right)$ with the fourth order constitutive tensor $\mathbf{D}$, here for isotropic linear elasticity. 

To weakly enforce prescribed displacements on Dirichlet boundaries, the static equilibrium in Eq.~\eqref{Eq_RLin} is augmented with Nitsche's method, see \cite{NITSCHE1971}:
\begin{equation}
\renewcommand{\arraystretch}{2.25}
\begin{array}{lll}
\mathcal{R}_{\Nitsche} 
& = & \displaystyle - \int_{\Gamma^{AB}} \delta \mathbf{u} \cdot \left( \boldsymbol{\sigma} \cdot \mathbf{n}_{\Gamma} \right) d\Gamma\\
& & \displaystyle + \int_{\Gamma^{AB}} \delta \left( \boldsymbol{\sigma} \cdot \mathbf{n}_{\Gamma} \right) \cdot \left( \mathbf{u} - \mathbf{u}_D \right)\ d\Gamma\\ 
& & \displaystyle + \gamma_{N} \int_{\Gamma^{AB}} \delta \mathbf{u} \cdot \left( \mathbf{u} - \mathbf{u}_D \right)\ d\Gamma,
\end{array}
\label{Eq_RNitsche}
\end{equation}
where $\mathbf{u}_D$ are the displacement imposed on the Dirichlet boundary $\Gamma^{A}_{D}$. The parameter $\gamma_{N}$ is chosen to achieve a certain accuracy in satisfying the boundary conditions and is a multiple of the ratio $E/h$, where $E$ is the Young's modulus of the considered material and $h$ is the edge length of the intersected elements.

Using immersed boundary techniques, see Sec\-tion \ref{XFEM}, numerical instabilities arise when either the contributions of the DOFs interpolating the displacement field on the residual vanish and/or these contributions become linearly dependent. These issues typically arise when the level set field intersects elements such that small material subdomains emerge. This results in an ill-conditioning of the equation system and inaccurate prediction of displacement gradients along the interface. The face-oriented Ghost stabilization proposed by \cite{BURMAN2014} is used in this work to mitigate this issue. Using a virtual work based formulation, the jump in the displacements gradients is penalized across the faces belonging to intersected elements by augmenting the residual equations with:
\begin{equation}
\mathcal{R}_{\Ghost} = h\ \gamma_{G}\ \sum_{F \in \mathcal{F}_{\cut}} \int_{F} \llbracket \delta \boldsymbol{\varepsilon} \cdot \mathbf{n}_F \rrbracket \llbracket \boldsymbol{\sigma} \cdot \mathbf{n}_{F} \rrbracket\ d\Gamma,
\label{Eq_RGhost}
\end{equation}
where $F$ is a specific face belonging to $\mathcal{F}_{\cut}$, the set of element faces cut by the interface, and $\mathbf{n}_{F}$ is the  outward normal to face $F$. The jump operator $\llbracket \bullet \rrbracket$ computes quantities between adjacent elements $A$ and $B$ as $\llbracket \bullet \rrbracket = \bullet^{A} -\bullet^{B}$. The Ghost penalty parameter is denoted $\gamma_{G}$ and typically takes its value in the range $0.1 \dots 0.001$. In contrast to the displacement based formulation of \cite{BURMAN2014}, Eq.~\eqref{Eq_RGhost} allows for different materials in adjacent elements.
 
During the optimization process, isolated islands of material can emerge and develop in the design domain. This leads to a singular system of equations, as the rigid body modes associated to these free-floating material islands are not suppressed or constrained. To prevent this issue, we adopt the selective structural spring approach proposed and successfully implemented in \cite{VILLANUEVA2017, GEISS2018, GEISS2019a}. An additional stiffness term is added to the weak form of the governing equations: 
\begin{equation}
\mathcal{R}_{\Spring} = \int_{\Omega^A}\ \gamma_{\spring}\ k_{\spring}\ \delta \mathbf{u} \cdot \mathbf{u}\ d\Omega,
\label{Eq_RSpring}
\end{equation}
where the parameter $\gamma_{\spring}$ is set to a value ranging from 0 to 1 depending on the solution of an additional diffusion problem, so that a fictitious spring stiffness $k_{\spring}$ is only applied to solid subdomains disconnected from any mechanical boundary conditions. The spring stiffness $k_{\spring}$ is set to $E/h^2$, where $E$ is the Young's modulus of the considered material and $h$ the element edge width.

\subsection{The extended finite element method (XFEM)}\label{XFEM}
The XFEM was developed by \cite{MOES1999} to model crack propagation without remeshing. The me\-thod allows capturing discontinuous or singular behaviors within a mesh element by adding specific enrichment functions to the classical finite element approximation. Here, following the work by \cite{TERADA2003}, \cite{HANSBO2004}, \cite{MAKHIJA2014}, a generalized Heaviside enrichment strategy with multiple enrichment levels is used. This particular enrichment enhances the classical finite element interpolation with additional shape functions which then results in independent interpolation into disconnected subdomains.

Considering a two-phase problem, the displacement field $\mathbf{u}(\mathbf{x})$ is approximated with the Heaviside enrichment: 
\begin{equation}
\renewcommand*{\arraystretch}{1.5}
\begin{array}{l}
\displaystyle \mathbf{u}^h(\mathbf{x}) = \sum_{m=1}^{M}  \left( H(-\phi(\mathbf{x})) \sum_{i \in I^{\star}}\ B_i(\mathbf{x})\ u_{im}^{A}\delta_{ml}^{Ai} \right.\\ 
\displaystyle \hspace{2.5cm} \left. +\ H(\phi(\mathbf{x})) \sum_{i \in I^{\star}}\ B_i(\mathbf{x})\ u_{im}^{B}\ \delta_{ml}^{Bi} \right),
\end{array}
\end{equation}
where $I^{\star}$ is the set of all the B-spline coefficients in the analysis mesh, $B_i(\mathbf{x})$ is the B-spline basis function associated with the $i^{th}$ coefficient and $u^I_{im}$ is the vector of displacement DOF associated with the $i^{th}$ coefficient for material phase $I = A, B$. The number of active enrichment levels is denoted $M$ and the Kronecker delta $\delta_{ml}^{Ii}$ is used to select the active enrichment level $m$ for the $i^{th}$ coefficient and material phase $I$. This selection ensures the satisfaction of the PU principle (\cite{BABUSKA1997}), as only one set of DOFs is used to interpolate the solution at a given point $\mathbf{x}$. For simplicity, the interpolation of the displacement field is restricted to first order functions. The Heaviside function $H$ is defined as:
\begin{equation}
H(z) = \left\{
\begin{array}{ll}
1, & \quad \mbox{if}\ z > 0,\\
0, & \quad \mbox{if}\ z < 0.
\end{array}
\right.
\end{equation}

Using the Heaviside enrichment, the integration of the weak form of the governing equations has to be performed separately on each material phase. Therefore, intersected elements are decomposed into integration subdomains, i.e., into triangles in two dimensions and into tetrahedra in three dimensions. For further details on the decomposition approach, the reader is referred to \cite{VILLANUEVA2014}.

\section{Optimization problem}\label{optProblem}
In this paper, we focus on minimum compliance designs considering mass either via a constraint or a component of the objective. The optimization problem is formulated as follows: 
\begin{equation}
\renewcommand{\arraystretch}{1.25}
\begin{array}{ccl}
\displaystyle \min_{0 \leq \mathbf{s} \leq 1} & \mathcal{Z}(\mathbf{s}, \mathbf{u}(\mathbf{s}))\ &+\ c_p\ \mathcal{P}_{p}(\mathbf{s})/\mathcal{P}_0\\
&\displaystyle +\  c_{\phi}\ \left(\ \mathcal{P}_{\phi}\left( \phi(\mathbf{s}) \right) \right. 
&+\ \left. \mathcal{P}_{\nabla \phi} \left( \nabla \phi(\mathbf{s}) \right) \right)/\mathcal{P}_{0} \\
\mbox{s.t.}
&\displaystyle \mathcal{M}^{A}(\mathbf{s})/\mathcal{M}(\mathbf{s}) - c_m & \leq 0.\\
\end{array}
\label{Eq_optFormulation}
\end{equation}

The functional $\mathcal{Z}$ is a weighted sum of the compliance and the mass, both evaluated over the solid domain, 
\begin{equation}
\mathcal{Z}(\mathbf{s}, \mathbf{u}(\mathbf{s})) = w_s \frac{\mathcal{S}( \mathbf{s}, \mathbf{u}(\mathbf{s}))}{\mathcal{S}_0} + w_m \frac{\mathcal{M}^{A}(\mathbf{s})}{\mathcal{M}_0}, 
\end{equation}
where $\mathcal{S}$ is the strain energy, $\mathcal{S}_0$ is the initial strain energy value, $\mathcal{M}^{A}$ the mass of material phase $A$, $\mathcal{M}_0$ is the initial mass value and $w_s$ and $w_m$ are the weighting factors. The perimeter associated to the domain boundary is denoted by $\mathcal{P}_p$ and its initial value is $\mathcal{P}_0$. The penalty functions $\mathcal{P}_{\phi}$ and $\mathcal{P}_{\nabla \phi}$ are used to regularize the level set field around and away from the interface and are further discussed in Section~\ref{lsReg}. The parameters $c_p$ and $c_{\phi}$ are the penalties associated with the perimeter and the regularization respectively. The mass constraint ensures that the ratio of material phase $A$, $\mathcal{M}^{A}$, and the total mass, $\mathcal{M}$, is smaller or equal to $c_m$. Note that enforcing the bounds on the design variable $\mathbf{s}$ is conveniently achieved using THB-splines as they form a convex hull, while it requires extra handling when using HB-splines.

On purpose, we do not consider any means to control the feature size. One goal of this paper is to study the influence of the proposed adaptive hierarchical B-spline discretization on the optimization results, including the ability to resolve fine features. To allow the emergence of such features, no feature size control is imposed.  

In this work, we solve the discretized optimization problem in the reduced space of only the design variables by a gradient-based algorithm. For each candidate design generated in the optimization process, the state variables are determined by solving the discretized governing equations; see Section~\ref{StructAnalysis}. 

\subsection{Level set regularization scheme}\label{lsReg}
\begin{figure}[h]\center
\includegraphics[scale=1]{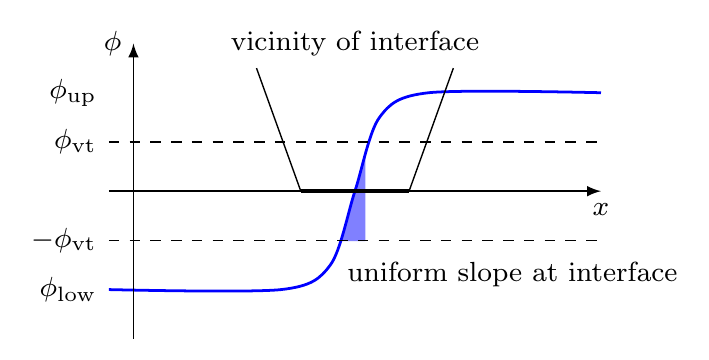}
\caption{Target level set function.}
\label{fig_lsReg}
\end{figure}

To ensure that the spatial gradient near the solid-void interface is approximately uniform and matches a target gradient norm over the span $[-\phi_{\vt},\ \phi_{\vt}]$, a level set regularization scheme is used. The scheme also enforces convergence to a target positive $\phi_{\up}$ or negative $\phi_{\low}$ value away from the interface, as depicted in Fig.~\ref{fig_lsReg}.

\begin{figure*}[ht]\center
\begin{tabular}{ccc}
\includegraphics[scale=1]{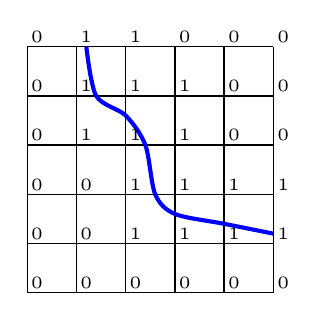} \quad
&\includegraphics[scale=1]{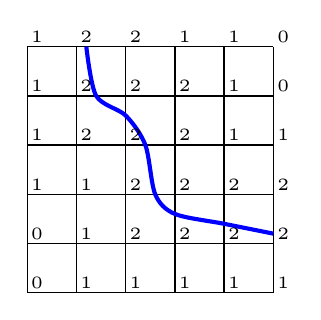} \quad
&\includegraphics[scale=1]{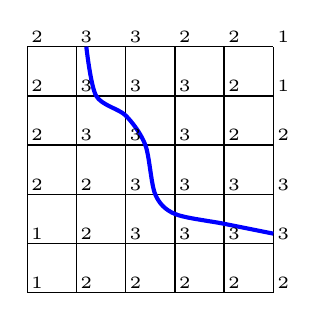} \quad \\
$\mathcal{I}_{\max} = 1$ & $\mathcal{I}_{\max} = 2$ & $\mathcal{I}_{\max} = 3$
\end{tabular}
\caption{Definition of the \textit{level of neighborhood} (LoN).}
\label{fig_lsRegNeighbor}
\end{figure*}

Such a scheme requires the identification of the vicinity of the solid-void interface. The level set values are often used for this purpose. However, unless the LSF is constructed as an approximation or as an exact signed distance field, this approach lacks robustness and may lead to spurious oscillations in the level set field, see \cite{GEISS2019b}. In this study, we do not compute a signed distance field. Instead, we identify the proximity of a point to the solid-void interface by building a so-called \textit{level of neighborhood} (LoN) $\mathcal{I}(\mathbf{x})$ of a point with respect to a point belonging to an intersected element. The LoN is evaluated at the nodes of the considered mesh, providing the LoN value $\mathcal{I}_i$ at each node $i$. The approach is illustrated in Fig.~\ref{fig_lsRegNeighbor}. The process starts with identifying intersected elements and setting the LoN $\mathcal{I}$ of all nodes belonging to these elements to one; otherwise $\mathcal{I}$ is set to zero. In a recursive loop over all elements and nodes, first the elements having a node with a LoN value larger than zero are flagged and then the LoN values of nodes belonging to a flagged element are increased by one. The region with LoN values larger than zero widens with the number of times the loop, $N_{\LoN}$, is executed. Thus, the maximum LoN value $\mathcal{I}_{\max}$ equals $N_{\LoN}$. The use of these nodal LoN values is discussed below.

To promote the convergence of the level set field to a target field, two penalty functions, $\mathcal{P}_{\phi}$ and $\mathcal{P}_{\nabla \phi}$, are added to the optimization problem objective, see Eq.~\eqref{Eq_optFormulation}. The first term penalizes the difference between the level set value $\phi$ and, depending on its sign, the lower or upper target value $\phi_{\target}$ away from the interface:
\begin{equation}
\mathcal{P}_{\phi} = \displaystyle \int_{\Omega} \alpha_1\ (1-w)\ \left( \frac{\phi(\mathbf{s})}{\phi_{\target}} - \sign(\phi) \right)^2\ d\Omega,
\label{Eq_lsReg1}
\end{equation}
where $\phi_{\target}$ is set as a multiple of the element edge size of the initial mesh $h_{\init}$, i.e., the mesh with the initial refinement level.

The second term penalizes both the difference between the norm of the spatial gradient of the level set field $\nabla \phi$ and a target gradient norm $\nabla \phi_{\target}$ near the interface and the norm of the spatial gradient of the level set field $\nabla \phi$ away from the interface:
\begin{equation}
\begin{array}{lll}
\mathcal{P}_{\nabla \phi} & = &\displaystyle \int_{\Omega} \alpha_2\ w\ \left( \frac{\vert \vert \nabla \phi(\mathbf{s}) \vert \vert}{\vert \vert \nabla \phi_{\target} \vert \vert} - 1.0 \right)^2\ d\Omega\\[10pt]
& + &\displaystyle \int_{\Omega} \alpha_3\ (1-w)\ {\vert \vert \nabla \phi \vert \vert\ }^2\ d\Omega.\\
\end{array}
\label{Eq_lsReg2}
\end{equation}
where $\alpha_1$, $\alpha_2$ and $\alpha_3$ are weighting factors. The parameter $w$ is a measure of the distance of a point to the solid-void interface as:
\begin{equation}
\displaystyle w = e^{-\gamma_{I} \left( \mathcal{I}(\mathbf{x})/\mathcal{I}_{\max}-1 \right)^2},
\label{Eq_lsRegWeight}
\end{equation}
where $\mathcal{I}(\mathbf{x})$ is the local LoN value, interpolated by the element shape functions using the nodal $\mathcal{I}_i$ values. The parameter $\gamma_{I}$ determines the region around and away from the interface in which the integrals in~\eqref{Eq_lsReg1} and~\eqref{Eq_lsReg2} are evaluated. As the value of $\gamma_{I}$ is increased, the area considered to be in the vicinity of the interface, i.e., the area where the target slope $\vert \vert \nabla \phi_{\target} \vert \vert$ is promoted, shrinks. The parameters $\phi_{\target}$, $\nabla \phi_{\target}$, and $\gamma_{I}$ are user-defined.

\subsection{Sensitivity analysis}\label{sensitivity}
In this study, we compute the design sensitivities by the adjoint approach. Let us consider an objective or constraint function $\mathcal{F}(\mathbf{u}(\mathbf{s}), \mathbf{s})$ dependent on the design variables. The derivative of this response function with respect to a design variable $s_i$ is computed as follows:
\begin{equation}
\frac{d \mathcal{F}(\mathbf{u}(\mathbf{s}))}{d s_i} = \frac{\partial \mathcal{F}}{\partial s_i} + \frac{\partial \mathcal{F}}{\partial \mathbf{u}} \frac{d \mathbf{u}}{d s_i},
\end{equation}
where the first term of the right-hand side accounts for the explicit dependency on the design variables and the second term for the implicit dependency on the design variables through the state variables. The implicit term is evaluated by the adjoint approach:
\begin{equation}
\frac{\partial \mathcal{F}}{\partial \mathbf{u}} \frac{d \mathbf{u}}{d s_i} = -\boldsymbol{\lambda}^{T} \frac{\partial \mathcal{R}}{\partial s_i},
\label{Eq_adjointApproach}
\end{equation}
where $\mathcal{R}$ is the residual defined in Eq.~\eqref{Eq_RFull} for the forward analysis and $\boldsymbol{\lambda}$ are the adjoint responses, evaluated through: 
\begin{equation}
\boldsymbol{\lambda} = \frac{\partial \mathcal{F}}{\partial \mathbf{u}} \left[ \frac{\partial \mathcal{R}}{\partial \mathbf{u}} \right]^{-1}.
\label{Eq_adjointResponse}
\end{equation}

As truncated B-spline basis functions have smaller support than the non-truncated ones, see Fig.~\ref{fig_truncation1D} and~\ref{fig_truncatedKidney}, the zone of influence associated to a given B-spline coefficient is also smaller. This property can be exploited to accelerate the sensitivity analysis.

When computing the design sensitivities with a LSM, only the intersected elements need to be considered for evaluating Eq.~\eqref{Eq_adjointResponse}, as only the finite element residuals of these elements depend on the level set field. However, when the seeding approach presented in Section~\ref{Hole} is used, the residual contributions of non-intersected elements within the solid domain also depend on the design variables. Thus, in this case, the partial derivative $\frac{\partial \mathcal{R}}{\partial \mathbf{u}}$ needs to be integrated over all elements in the solid domain.

\section{Numerical examples}\label{numExample}
We study the proposed level set topology optimization with hierarchical mesh refinement with 2D and 3D examples, focusing on minimum compliance designs considering mass either via a constraint or a component of the objective, see Eq.~\eqref{Eq_optFormulation}.

The state variable field, i.e., the displacement field, is discretized in space with bilinear Lagrange quadrangular elements for 2D, and with trilinear Lagrange hexahedral elements for 3D design domains. The design variable field, i.e., the level set field, is discretized with HB- and THB-splines. The design variables are the B-splines coefficients and are used to interpolate nodal values on the Lagrange analysis mesh.

The optimization problems are solved by GCMMA (\cite{SVANBERG2002}) and the required sensitivity analysis is performed following the adjoint approach described in Section~\ref{sensitivity}. The parameters for the initial, lower and upper asymptote adaptation in GCMMA are set to 0.05, 0.65 and 1.05 respectively. The optimization problem is considered converged if the absolute change of the objective function relative to the mean of the objective function in the five previous optimization steps drops below $10^{-5}$ and the constraint is satisfied.

The refinement strategy described in Section~\ref{refinementStrategies} is adopted and the maximum refinement level for interface and solid elements is set to $l_{\ifc,\max} = l_{\solid,\max} = l_{\max} = 4$. The minimum refinement level is set to $l_{\min} = 0$, the size of the coarsest mesh and the initial mesh refinement level are given with each example.

The system of discretized governing equations and adjoint sensitivity equations are solved by the direct solver PARDISO for the 2D problems, see \cite{KOUROUNIS2018}, and by a GMRES algorithm for 3D problems, preconditioned by an algebraic multi-grid solver, see \cite{MLGUIDE}. A relative drop of the linear residual of $10^{-12}$ is required. This strict tolerance is imposed to reduce the influence of the iterative solver on the optimization results.

For the first two examples, we provide detailed comparisons between designs generated on uniformly and adaptively refined meshes in terms of performance, a\-chieved geometries and computational cost. Furthermore, we investigate the influence of the B-spline interpolation by varying the interpolation order and applying or omitting truncation. It should be noted that extra treatment is required to impose bounds on the design variables when working with HB-splines, as they do not constitute a PU. Bounds are enforced by clipping the variables values with the upper or lower allowed values. The clipping operation is only applied to the combined level set/density scheme as the density values should remain between 0 and 1.

Providing meaningful computational efficiency measures for the adaptive simulations is not trivial as various factors influence the computational cost. As the simulations could not be carried out on the same hardware with the exact same settings due to limited availability to dedicated computing resources, a simple wall-clock time comparison is not possible. Furthermore, different linear solvers and software implementations lead to different dependencies of the computational time on the number of finite element DOFs. Therefore, we measure the computational cost in terms of the evolution of the number of unconstrained DOFs of the finite element models during the optimization process.

An efficiency factor $E_{\xfem}$ is defined as the ratio between the total number of unconstrained DOFs in the XFEM model for the most refined uniform mesh and the total number of unconstrained DOFs in the XFEM model for an adaptive mesh. This measure is indicative of the computational gains achieved for XFEM solid-void problems and is expressed as follows:
\begin{equation}
\displaystyle E_{\xfem} = \frac{\displaystyle \sum_{k = 1}^{N_{\opt}}\ \left( \#DOFs_{\uniform,\xfem}^{(k)} \right)^{n_s} }{\displaystyle \sum_{k=1}^{N_{\opt}}\ \left( \#DOFs_{\adaptive}^{(k)} \right)^{n_s} },
\label{Eq_compCostXFEM}
\end{equation}
where $\#DOFs_{\uniform,\xfem}^{(k)}$ and $\#DOFs_{\adaptive}^{(k)}$ are the number of unconstrained DOFs in the XFEM models for the uniform and the adaptive meshes at the optimization iteration $k$. $N_{\opt}$ is the number of optimization steps until convergence and is set to $N_{\opt} = \min( N_{\opt, \uniform}, N_{\opt, \adaptive})$. The parameter $n_s$ is an exponent that can be fitted to relate the number of DOFs to the computational effort in terms of floating-point operations or wall clock time, see~\cite{WOZNIAK2014}, for example when a sparse direct linear solver is used and is set to $n_s = 1$ for the 2D and $n_s = 4/3$ for the 3D cases.

Alternatively, a measure of the peak resource requirement $R_{\xfem}$ is defined as the ratio between the maximum number of unconstrained DOFs in the XFEM model for the uniform mesh with highest refinement level and the maximum number of unconstrained DOFs in the XFEM model for the adaptive mesh:
\begin{equation}
\displaystyle R_{\xfem} = \frac{\displaystyle \max \left( \#DOFs_{\uniform,\xfem} \right)^{n_s} }{\displaystyle \max \left( \#DOFs_{\adaptive} \right)^{n_s}}.
\label{Eq_maxCompCostXFEM}
\end{equation}

It should be noted that for an XFEM model, the number of unconstrained DOFs is changing in the course of the optimization process, even on a uniformly refined mesh. For the solid-void problems studied in this paper and because the void domain is omitted in the analysis, the XFEM model reduces the number of DOFs compared to a classical FEM model. To account for these computational gains, we define an efficiency factor $E_{\fem}$, indicative for solid-solid XFEM or solid-void FEM problems, as: 
\begin{equation}
\displaystyle E_{\fem} = \frac{\displaystyle \sum_{k = 1}^{N_{\opt}}\ \left( \#DOFs_{\uniform,\fem}^{(k)} \right)^{n_s} }{\displaystyle \sum_{k=1}^{N_{\opt}}\ \left( \#DOFs_{\adaptive}^{(k)} \right)^{n_s} },
\label{Eq_compCostFEM}
\end{equation}
where $\#DOFs_{\uniform,\fem}^{(k)}$ is an approximation of the number of unconstrained DOFs for the uniform mesh with highest refinement level at the optimization iteration $k$, assuming that all DOFs are active as is the case for a classical FEM model.

A corresponding measure of the peak resource need $R_{\fem}$ is defined as follows: 
\begin{equation}
\displaystyle R_{\fem} = \frac{\displaystyle \max\ \left( \#DOFs_{\uniform,\fem} \right)^{n_s}}{\displaystyle \max\ \left( \#DOFs_{\adaptive} \right)^{n_s}}.
\label{Eq_maxCompCostFEM}
\end{equation}

The measures given by Eq.~\eqref{Eq_compCostXFEM},~\eqref{Eq_compCostFEM},~\eqref{Eq_maxCompCostXFEM},~\eqref{Eq_maxCompCostFEM} for estimating the gains in computational efficiency are solely based on the costs for solving the linear finite element system. It does not account for the costs for building the XFEM model for initializing the linear solver, for the $L^{2}$ projection of the design variable fields, see Section~\ref{refinementStrategies}, for computing the right-hand side of the adjoint system, see Eq.~\eqref{Eq_adjointApproach}, and for the optimization algorithm. The latter computational costs are, however, small typically when compared to building and solving the system of governing equations. 

\subsection{Two dimensional beam}\label{2DBeam}
\begin{figure}[h]\center
\includegraphics[scale=1]{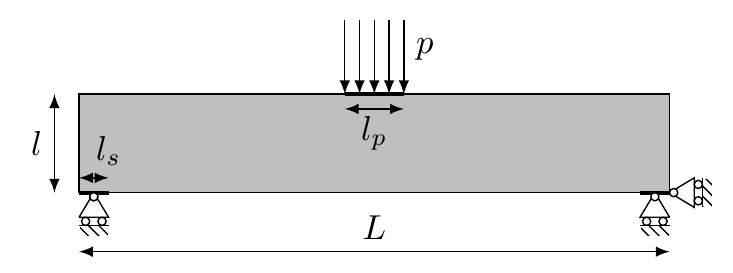}
\caption{Two dimensional beam.}
\label{fig_2D_beam}
\end{figure}

As illustrated in Fig.~\ref{fig_2D_beam}, the first problem setup consists of a solid-void 2D beam with a length $L = 6.0$ and a height $l = 1.0$. The solid is described by a linear elastic model and an isotropic constitutive behavior, with a Young's modulus $E = 1.0$ and a Poisson's ratio $\nu = 0.3$. The beam is supported on its two lower corners over a length $l_s = 0.025$ and a pressure $p = 1.0$ is applied in the middle of its top span over a length $l_p = 2 \times 0.025$. The data are provided in self-consistent units.

The structure is designed for minimum compliance with a mass constraint of $40 \%$ of the total mass of the design domain. Taking advantage of the symmetry of the problem, only one half of the design domain is considered. Table~\ref{table_2DMBB_parameter} summarizes the problem parameters. 
\begin{table}[h]\center
\caption{Parameter list for the 2D beam optimization problem.}
\label{table_2DMBB_parameter}
\renewcommand{\arraystretch}{1.5}
\begin{tabular}{p{0.45\columnwidth}p{0.45\columnwidth}}
\hline
Parameter & Value\\ \toprule
$w_s$ &0.9\\
$w_m$ & 0.0\\
$c_p$ &0.025\\
$c_{\phi}$ & 0.5\\
$c_v$ &0.4\\
$S_0$ &59.12\\
$P_0$ & 6.0\\
$\phi_{\scale}$ & $5h_{\init}$\\
$V_{\max}$ & 6.0\\
$\gamma_N$ & $100E/h$\\
$\gamma_G$ & $0.005$\\
$\phi_{\target}$ &$1.5 h_{\init}$\\
$\nabla \phi_{\target}$ &0.75\\ 
$\mathcal{I}_{\max}$ &1\\ 
$\gamma_{I}$ & 4.61\\
$\alpha_1 = \alpha_2 = \alpha_3$ &0.5 \\ \hline
\end{tabular}
\end{table}

The compliance minimization problem, defined in Eq.~\eqref{Eq_optFormulation}, is solved on uniformly and adaptively refined meshes, considering five levels of refinement $l_{\refine} = 0$, $1$, $2$, $3$ and  $4$, corresponding to $30 \times 10$, $60 \times 20$, $120 \times 40$, $240 \times 80$ and $480 \times 160$ element meshes. The influence of the B-spline interpolation is investigated by solving the same design problem with linear, quadratic and cubic B-spline functions and with truncated and non-truncated B-splines. These influences are studied separately considering the two approaches for introducing inclusions, i.e., the initial hole seeding and the combined level set/density scheme. 

\subsubsection{Design with initial hole seeding}\label{2DMBB_initial_hole}
\begin{figure}[h]\center
\includegraphics[scale=1]{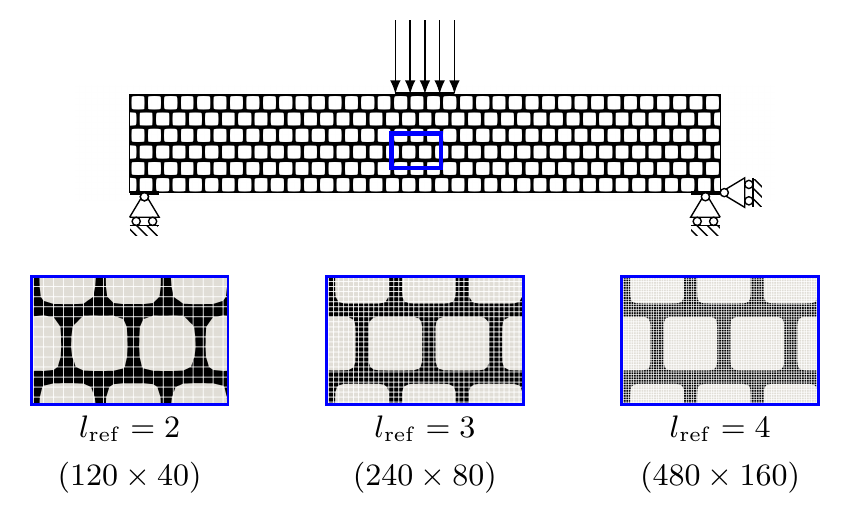}
\caption{Initial hole seeding for the 2D beam. Zoom on the hole pattern representation for the refinement levels $l_{\refine} = 2, 3,$ and $4$.}
\label{fig_2D_MBB_initial_hole}
\end{figure}

In this subsection, the design domain is initially seeded with holes, so that the mass constraint is satisfied at the beginning of the optimization process. The initial configuration is depicted in Fig.~\ref{fig_2D_MBB_initial_hole}. The complexity of such a hole pattern requires a rather fine initial mesh and cannot be represented with the coarsest refinement levels $l_{\refine} = 0, 1$. Thus, the optimization process needs to start from a mesh with $l_{\refine} > 1$.

\begin{figure*}[ht]\center
\renewcommand*{\arraystretch}{1.75}
\begin{tabular}{p{1.5cm}ccc}
& Linear B-splines & Quadratic B-splines & Cubic B-splines\\
$l_{\refine} = 2$
&\includegraphics[width=4.29cm]{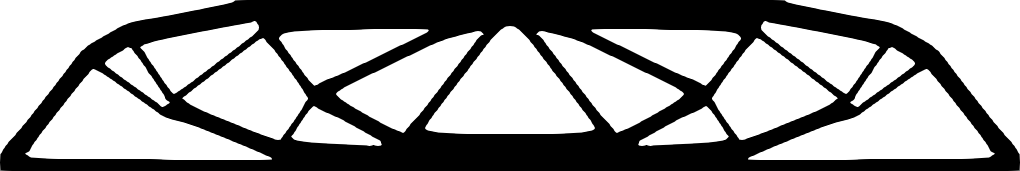}
&\includegraphics[width=4.29cm]{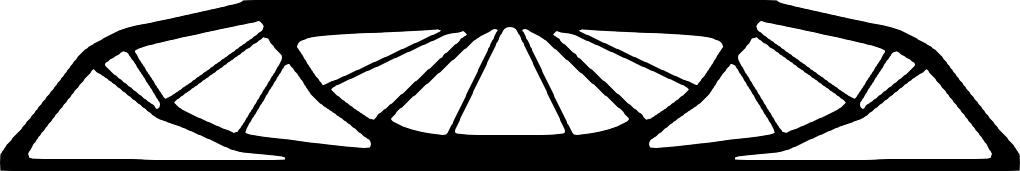}
&\includegraphics[width=4.29cm]{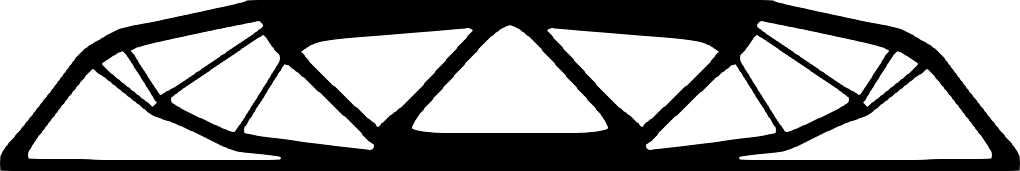}\\[-5pt]
$120 \times 40$ 
&$\mathcal{S} = 101.01$  &$\mathcal{S} = 101.25$  &$\mathcal{S} = 102.53$\\[5pt] 

$l_{\refine} = 3$ 
&\includegraphics[width=4.29cm]{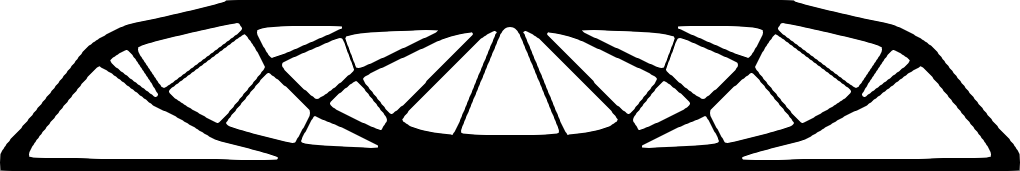}
&\includegraphics[width=4.29cm]{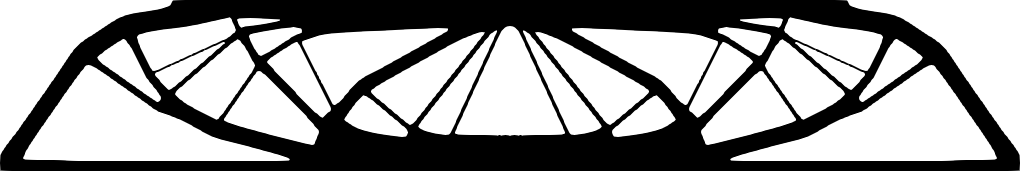}
&\includegraphics[width=4.29cm]{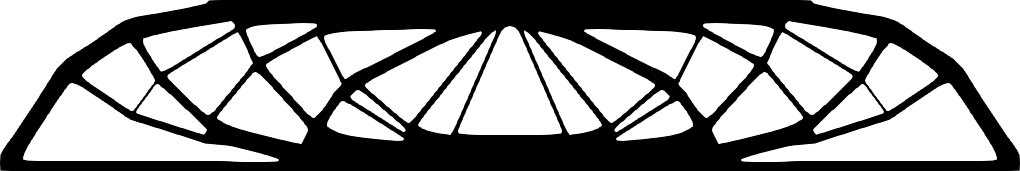}\\[-5pt]
$240 \times 80$
&$\mathcal{S} = 100.63$  &$\mathcal{S} = 101.63$  &$\mathcal{S} = 100.61$\\[5pt] 

$l_{\refine} = 4$ 
&\includegraphics[width=4.29cm]{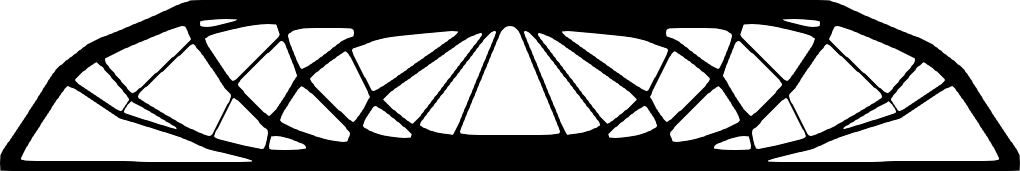}
&\includegraphics[width=4.29cm]{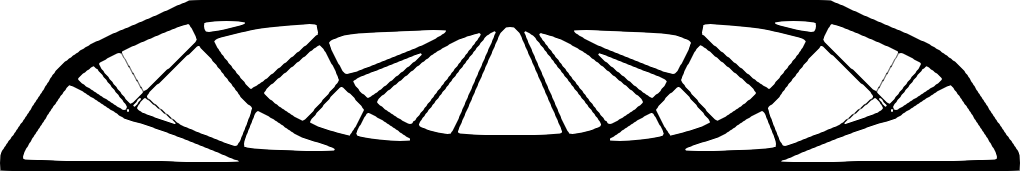}
&\includegraphics[width=4.29cm]{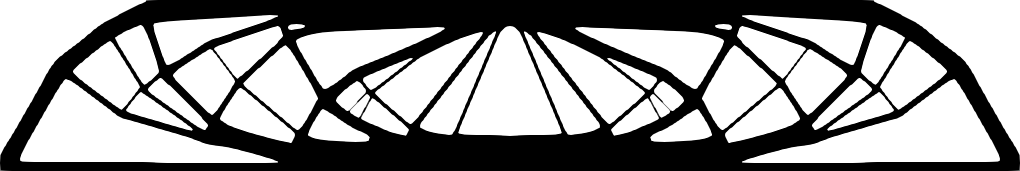}\\[-5pt]
$480 \times 160$
&$\mathcal{S} = 101.02$  &$\mathcal{S} = 100.61$  &$\mathcal{S} = 101.32$ \\
\end{tabular}
\caption{2D beam using uniformly refined meshes with THB-splines and initial hole seeding.}
\label{fig_2D_MBB_uniform}
\end{figure*}

Figure~\ref{fig_2D_MBB_uniform} shows the optimized designs generated on uniformly refined meshes with $l_{\refine} = 2, 3$ and $4$ using THB-splines and for different B-spline orders, along with the associated mechanical performance in terms of strain energy. The solid phase is depicted in black and the meshes are omitted as they are uniform. The designs generated on uniform meshes for different B-spline orders become increasingly similar as the mesh is refined. The performances of the designs are also rather similar. As expected, the results show that refining the mesh allows for a higher geometric complexity, including thinner structural members. 

\begin{figure*}[ht]\center
\renewcommand*{\arraystretch}{1.75}
\begin{tabular}{p{1.5cm}ccc}
& Linear B-splines & Quadratic B-splines & Cubic B-splines\\
$l_{\refine}^{0} = 2$ 
&\includegraphics[width=4.29cm]{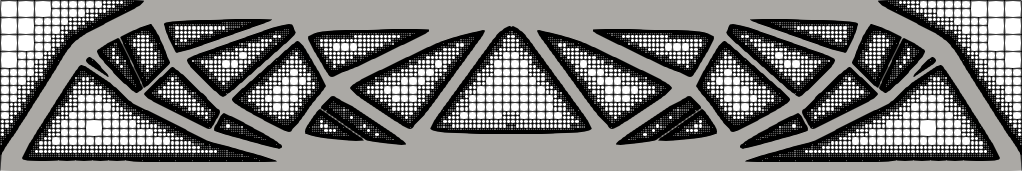}
&\includegraphics[width=4.29cm]{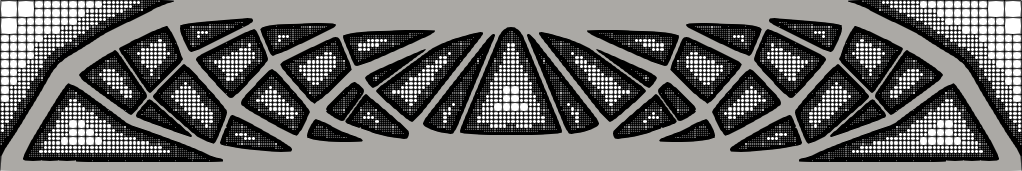}
&\includegraphics[width=4.29cm]{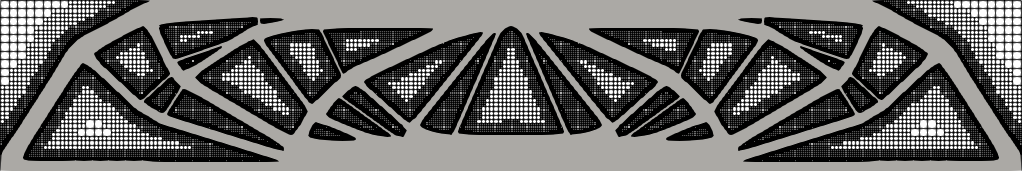}\\[-5pt]
$10$ iter. 
&$\mathcal{S} = 101.23$  &$\mathcal{S} = 99.67$  &$\mathcal{S} = 101.08$ \\[5pt]
 
$l_{\refine}^0 = 2$
&\includegraphics[width=4.29cm]{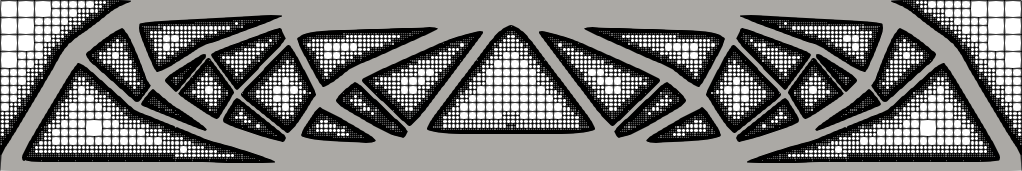}
&\includegraphics[width=4.29cm]{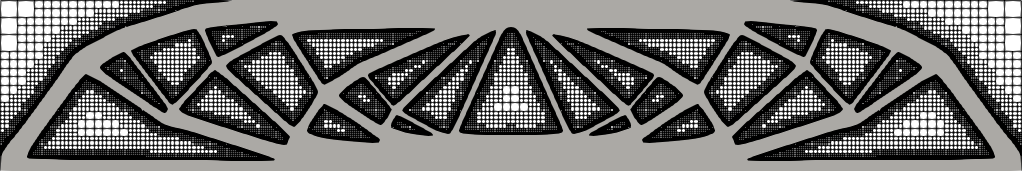}
&\includegraphics[width=4.29cm]{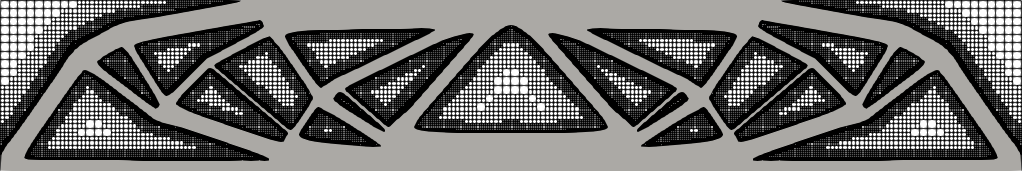}\\[-5pt]
$25$ iter.
&$\mathcal{S} = 102.03$  &$\mathcal{S} = 100.10$  &$\mathcal{S} = 101.75$ \\[5pt]

$l_{\refine}^0 = 2$
&\includegraphics[width=4.29cm]{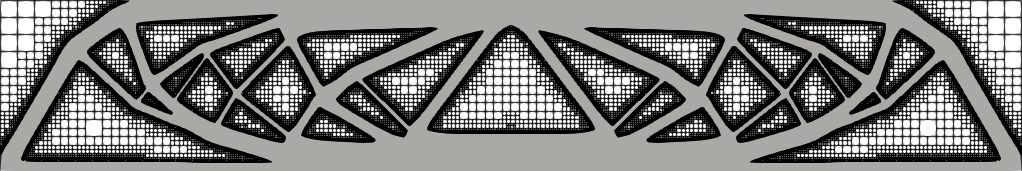}
&\includegraphics[width=4.29cm]{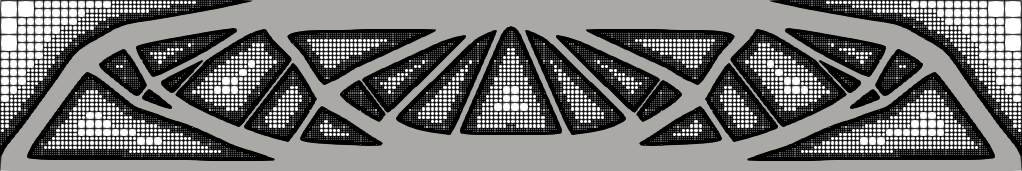}
&\includegraphics[width=4.29cm]{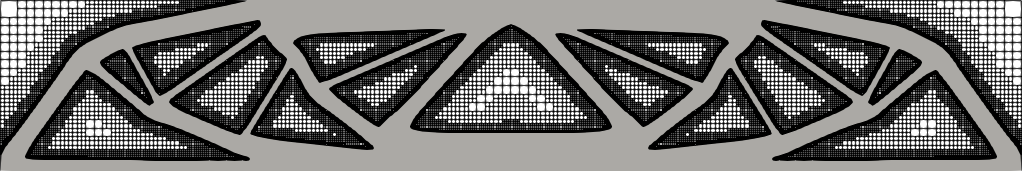}\\[-5pt]
$50$ iter.
&$\mathcal{S} = 102.33$  &$\mathcal{S} = 101.22$  &$\mathcal{S} = 102.78$ \\
\end{tabular}
\caption{2D beam using adaptively refined meshes with THB-splines and initial hole seeding. Initial refinement level of $l_{\refine}^{0} = 2$ and maximum refinement level of $l_{\refine,\max} = 4$.}
\label{fig_2D_MBB_adaptive}
\end{figure*}

Figure~\ref{fig_2D_MBB_adaptive} shows the optimized designs generated on adaptively refined meshes using linear, quadratic and cubic THB-splines. For each layout, the solid phase is depicted in gray and the refined mesh is shown in the void phase. The computational gains achieved with the adaptive strategy are given in Table~\ref{table_2D_MBB_compCost}.

Applying hierarchical mesh refinement, the optimization process is started with an initial uniform mesh with a refinement level $l_{\refine}^{0} = 2$ and is refined up to a refinement level $l_{\refine,\max} = 4$ and coarsened to a refinement level $l_{min} = 0$. The minimum and maximum refinement levels are updated after every 10, 25 or 50 iterations. Additionally, the mesh is adapted to maintain a uniform refinement level in all intersected elements, as mentioned in Section~\ref{userDefinedCriteria}. First, the designs demonstrate the ability of the method to recover similar layouts to the ones obtained on uniform meshes with higher mesh refinement. No significant loss is observed regarding the level of details or the complexity of the geometry achieved. The strain energy values are quite similar regardless of the use of uniform or adaptive meshes. The remaining differences in the strain energy values are not only due to differences in the designs, but also arise from differences in the discretization of the finite element model. 

The results in Fig.~\ref{fig_2D_MBB_adaptive} show, however, that the adaptive designs are sensitive to the number of iterations after which the refinement is applied. More frequent refinement operations, specifically early in the optimization process, allow for the emergence of finer features in the design. In terms of computational efficiency, the cost reduction with the adaptive strategy remains limited as initial configurations require rather fine meshes, i.e., $l_{\refine} > 1$. The efficiency factors in Table~\ref{table_2D_MBB_compCost} show that the achieved gain slightly increases with the B-spline order. The largest factors, in terms of efficiency and peak need, are achieved when the number of iterations before refinement is increased. Performing mesh refinement and updating the minimum and maximum refinement levels every 50 iterations lead to the largest computational gains. Using XFEM alone reduces the computational cost by about 50\%. The adaptation strategy allows for a reduction of the computational cost by more than 3.5 over the uniform FEM model and by more than 1.75 over the uniform XFEM model. The peak resource needs $R_{\xfem}$ and $R_{\fem}$ follow the same trend.

Considering the refined meshes in Fig.~\ref{fig_2D_MBB_adaptive}, it is noticeable that the refined bandwidth around the solid-void interfaces increases in size with the order of the B-spline interpolation. This is a direct consequence of the refined buffer zone described in Subsection~\ref{bufferCriterion}, which size depends on the interpolation order. For both, the uniform and adaptive designs the influence of the B-spline interpolation order on the designs is mainly noticeable when working with lower order B-splines. Using a linear interpolation allows for the emergence of finer members and details in the structure, at the increased risk of converging to a local minimum with inferior performance. Higher order designs tend to be smoother and to eliminate thin members. For this 2D case, it is noteworthy that working with linear B-splines does not lead to any smoothness issues or irregular shapes which might be caused by a spurious interplay between geometry and finite element prediction. 

\begin{table*}[ht]\center
\caption{Performance in terms of computational cost for designs in Fig.~\ref{fig_2D_MBB_uniform}.}
\label{table_2D_MBB_compCost}
\renewcommand*{\arraystretch}{1.5}
\begin{tabular}{l|m{1cm}m{1cm}m{1cm}|m{1cm}m{1cm}m{1cm}|m{1cm}m{1cm}m{1cm}}\hline
&\multicolumn{3}{c|}{Mesh update after 10 iter.}
&\multicolumn{3}{c|}{Mesh update after 25 iter.}
&\multicolumn{3}{c}{Mesh update after 50 iter.}\\
& Lin. & Quad. & Cub. 
& Lin. & Quad. & Cub. 
& Lin. & Quad. & Cub. \\ \hline
$N_{\opt}$ 	&1648 &1158 &840 &1460 &1158 &1254 &1219 &1158 &1137\\ \hline
$E_{\xfem}$ 	&1.85 &1.72 &1.78 &1.93 &1.85 &1.97 &2.08 &2.03 &2.23\\
$R_{\xfem}$ 	&1.35 &1.29 &1.26 &1.65 &1.70 &1.72 &1.81 &1.99 &2.11\\
$E_{\fem}$   	&3.81 &3.56 &3.61 &3.94 &3.81 &4.05 &4.20 &4.19 &4.58\\
$R_{\fem}$    	&2.38 &2.31 &2.22 &2.92 &3.03 &3.05 &3.21 &3.56 &3.74\\ \hline
\end{tabular}
\end{table*}

\begin{figure*}[ht]\center
\renewcommand*{\arraystretch}{1.75}
\begin{tabular}{p{1.75cm}ccc}
& Linear B-splines & Quadratic B-splines & Cubic B-splines\\
HB-splines
&\includegraphics[width=4.29cm]{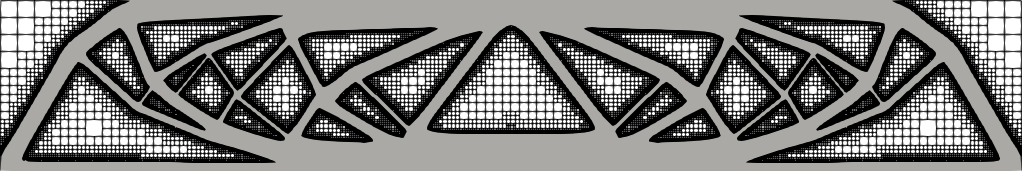}
&\includegraphics[width=4.29cm]{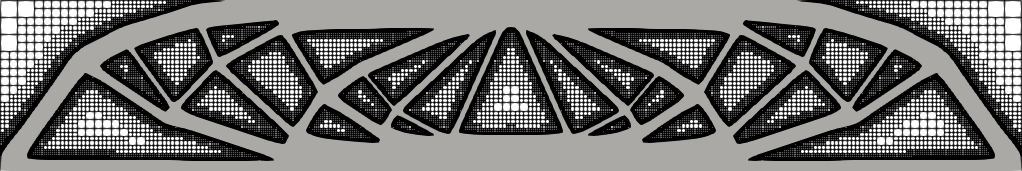}
&\includegraphics[width=4.29cm]{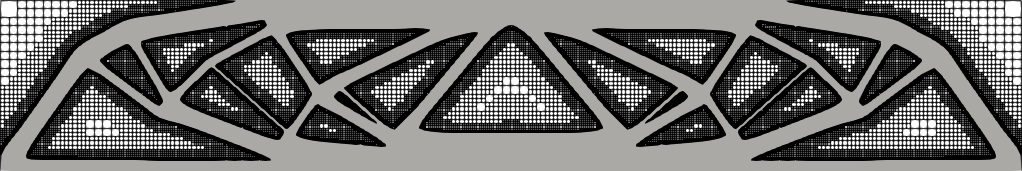}\\[-5pt]
$l_{\refine}^{0} = 2$
&$\mathcal{S} = 101.74$ &$\mathcal{S} = 100.14$ &$\mathcal{S} = 101.51$\\[5pt]
THB-splines
&\includegraphics[width=4.29cm]{figures/2DMBB/levsOnly_Adaptive/Linear_Adap_25_Start_2_mesh}
&\includegraphics[width=4.29cm]{figures/2DMBB/levsOnly_Adaptive/Quadratic_Adap_25_Start_2_mesh}
&\includegraphics[width=4.29cm]{figures/2DMBB/levsOnly_Adaptive/Cubic_Adap_25_Start_2_mesh}\\[-5pt]
$l_{\refine}^{0} = 2$
&$\mathcal{S} = 102.03$  &$\mathcal{S} = 100.10$  &$\mathcal{S} = 101.75$\\
\end{tabular}
\caption{2D beam using adaptively refined meshes with initial hole seeding. Initial refinement level of $l_{\refine}^{0} = 2$, maximum refinement level of $l_{\refine,\max} = 4$. Comparison of HB- and THB-splines.}
\label{fig_2D_MBB_adaptive_truncation}
\end{figure*}

The influence of the truncation operation on the B-splines is investigated by solving the 2D beam problem with HB-splines and THB-splines on an adaptive mesh with an initial refinement level $l_{\refine}^{0} = 2$ and applying the refinement operation every 25 iterations. The obtained designs and corresponding final strain energy values are given in Fig.~\ref{fig_2D_MBB_adaptive_truncation}. The linear, quadratic and cubic designs only differ slightly. The remaining differences between the designs can be explained by the support size of the HB- and THB-splines that differs, as shown in Fig.~\ref{fig_truncatedKidney}. This is further supported by the maximum stencil size recorded for each design variables, i.e., the maximum number of design coefficients affected by the change in a specific design coefficient. The stencil sizes for the HB-splines are 8, 16 and 33 for the linear, quadratic and cubic orders, against 2, 9 and 16 for the THB-splines. The efficiency factors are given in Table~\ref{table_2D_compCost_truncation} and match closely for HB- and THB-splines. 

\begin{table}[ht]\center
\caption{Performance in terms of computational cost for designs in Fig.~\ref{fig_2D_MBB_adaptive_truncation}.}
\label{table_2D_compCost_truncation}
\renewcommand*{\arraystretch}{1.5}
\begin{tabular}{l|ccc|ccc}\hline
&\multicolumn{3}{c|}{HB-splines}
&\multicolumn{3}{c}{THB-splines}\\ 
& Lin. & Quad. & Cub. 
& Lin. & Quad. & Cub. \\ \hline
$N_{\opt}$ 	&1385 &1158 &1254 &1460 &1158 &1254\\ \hline
$E_{\xfem}$  	&1.93 &1.83 &1.97 &1.93 &1.85 &1.97\\
$R_{\xfem}$    	&1.65 &1.70 &1.74 &1.65 &1.70 &1.72\\
$E_{\fem}$     	&3.93 &3.78 &4.06 &3.94 &3.81 &4.05\\
$R_{\fem}$    	&2.92 &3.04 &3.08 &2.92 &3.03 &3.05\\ \hline
\end{tabular}
\end{table}

\begin{figure*}[ht]\center
\renewcommand*{\arraystretch}{1.75}
\begin{tabular}{p{1.5cm}ccc}
& Linear B-splines & Quadratic B-splines & Cubic B-splines\\
Uniform
&\includegraphics[width=4.29cm]{figures/2DMBB/levsOnly_Uniform/Linear_Uniform_4}
&\includegraphics[width=4.29cm]{figures/2DMBB/levsOnly_Uniform/Quadratic_Uniform_4}
&\includegraphics[width=4.29cm]{figures/2DMBB/levsOnly_Uniform/Cubic_Uniform_4}\\[-5pt]
$l_{\refine} = 4$
&$\mathcal{S} = 101.02$  &$\mathcal{S} = 100.61$  &$\mathcal{S} = 101.32$\\[5pt]
Adaptive
&\includegraphics[width=4.29cm]{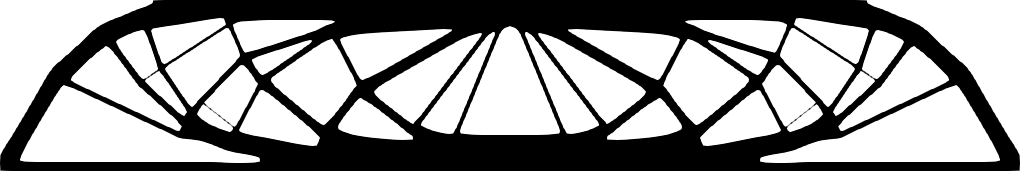}
&\includegraphics[width=4.29cm]{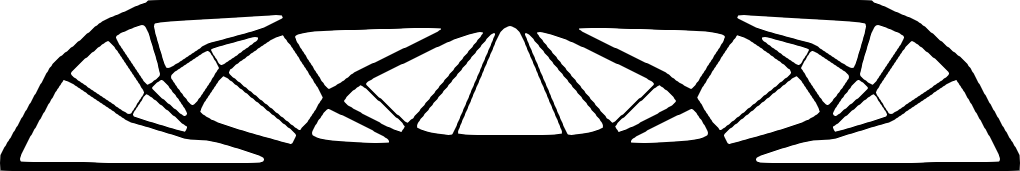}
&\includegraphics[width=4.29cm]{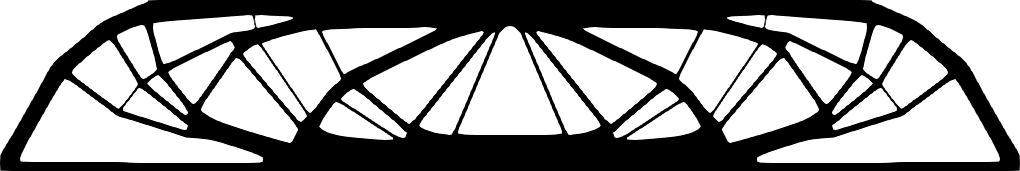}\\[-5pt]
$l_{\refine}^{0} = 4$
&$\mathcal{S} = 101.35$  &$\mathcal{S} = 101.55$  &$\mathcal{S} = 101.95$\\ 
\end{tabular}
\caption{2D beam using adaptively refined meshes with THB-splines and initial hole seeding. Initial refinement level of $l_{\refine}^{0} = 4$ and maximum refinement level of $l_{\refine,\max} = 4$. Only coarsening is applied to adapt the mesh.}
\label{fig_2D_MBB_adaptive_coarsen}
\end{figure*}

To further investigate the influence of mesh adaptivity, the optimization problem is solved starting from a uniform mesh with the highest refinement level, here $l_{\refine,\max}=4$ and using THB-splines. The mesh adaptation operation, triggered here every 25 iterations, involves coarsening of the mesh in void regions. The generated designs are shown in Fig.~\ref{fig_2D_MBB_adaptive_coarsen} along with the ones generated on the finest uniform mesh $l_{\refine}=4$. They slightly differ from the layouts obtained on uniform meshes, but present quite similar performance in terms of strain energy. These differences might result from the restart of the optimization process, as the mesh is adapted. After refinement, the level set field is mapped onto the newly generated mesh, which might result in a slight perturbation of the field values. Furthermore, the GCMMA is restarted with uniform lower and upper asymptotes for all variables which may alter the evolution of the design. Additionally, these results show that the obtained designs are dependent on the initial mesh refinement $l_{\refine}^{0}$. 

The efficiency factors are given in Table~\ref{table_2D_MBB_compCost_coarsening}. They show that the initial refinement level affects the computational cost. In this case, the mesh adaptivity does not provide an actual computational gain with respect to an XFEM model, as $E_{\xfem}$ and $R_{\xfem}$ are about equal to 1. 

\begin{table}[ht]\center
\caption{Performance in terms of computational cost for designs in Fig.~\ref{fig_2D_MBB_adaptive_coarsen}.}
\label{table_2D_MBB_compCost_coarsening}
\renewcommand*{\arraystretch}{1.5}
\begin{tabular}{l|ccc}\hline
&\multicolumn{3}{c}{$l_{\refine}^{0} = 4$}\\ 
& Lin. & Quad. & Cub.\\ \hline
$N_{\opt}$ 	&1911 &1158 &1131\\ \hline
$E_{\xfem}$ 	&0.99 &1.00 &1.00\\
$R_{\xfem}$    	&1.00 &1.00 &1.00\\
$E_{\fem}$    	&2.05 &2.05 &2.04\\
$R_{\fem}$     	&1.77 &1.78 &1.77\\ \hline
\end{tabular}
\end{table}

\subsubsection{Design with simultaneous hole seeding}\label{2DMBB_seeding}
In this section, the combined level set/density scheme, presented in Section~\ref{Hole}, is used to nucleate holes during the optimization process. The combined scheme uses a SIMP exponent $\beta_{\SIMP} = 2.0$. The density shift parameter $\rho_{\shift}$ is initially set to $0.2$ and is increased every 25 optimization steps. In this case, the optimization process can be started with rather coarse meshes. The minimum and maximum refinement levels are updated every 25 iterations. All other parameters are kept the same as in Section~\ref{2DMBB_initial_hole}. Numerical studies have shown that the chosen strategy for updating the density shift parameter and the minimum/maximum refinement levels lead to satisfactory results. Refining this strategy or using alternate approaches may lead to improved performance and should be explored in future studies.

\begin{figure*}[ht]\center
\renewcommand*{\arraystretch}{1.75}
\begin{tabular}{p{1.5cm}ccc}
& Linear B-splines & Quadratic B-splines & Cubic B-splines\\
$l_{\refine} = 1$
&\includegraphics[width=4.29cm]{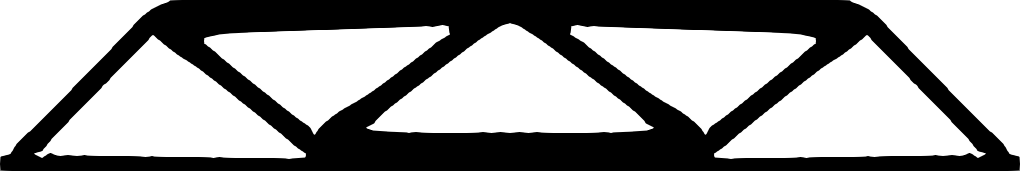}
&\includegraphics[width=4.29cm]{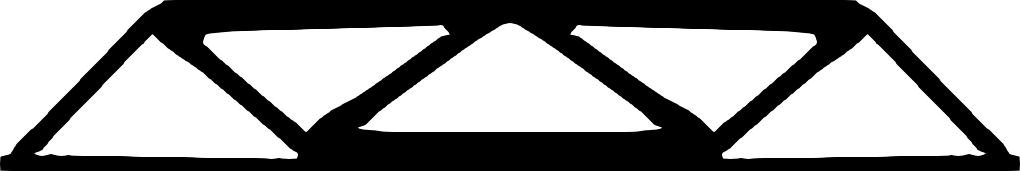}
&\includegraphics[width=4.29cm]{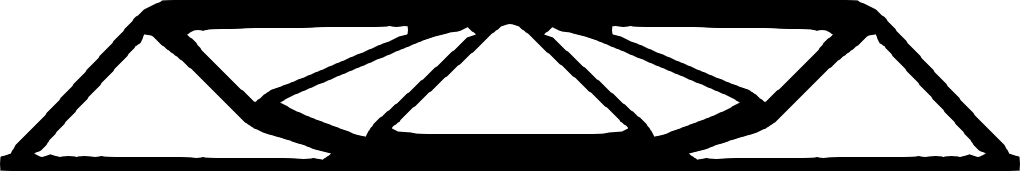}\\[-5pt] 
$60 \times 20$
&$\mathcal{S} = 101.11$ &$\mathcal{S} = 101.22$  &$\mathcal{S} = 99.48$\\[5pt]

$l_{\refine} = 2$
&\includegraphics[width=4.29cm]{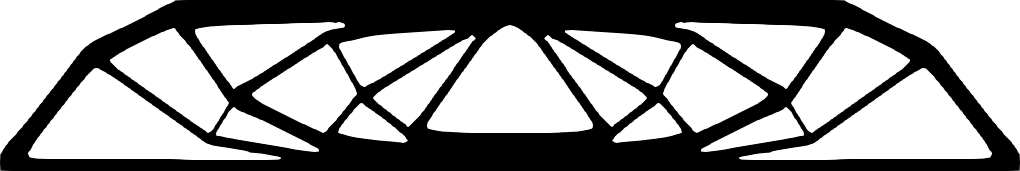}
&\includegraphics[width=4.29cm]{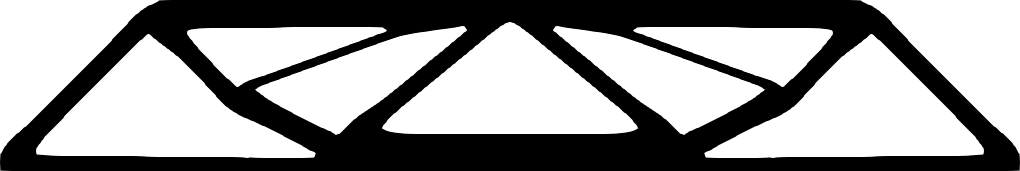}
&\includegraphics[width=4.29cm]{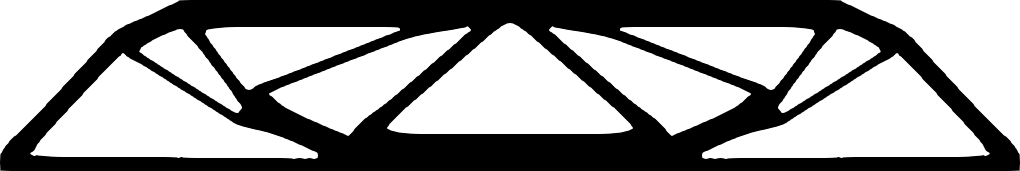}\\[-5pt] 
$120 \times 40$
&$\mathcal{S} = 100.81$ &$\mathcal{S} = 104.94$  &$\mathcal{S} = 103.77$\\[5pt]

$l_{\refine} = 3$
&\includegraphics[width=4.29cm]{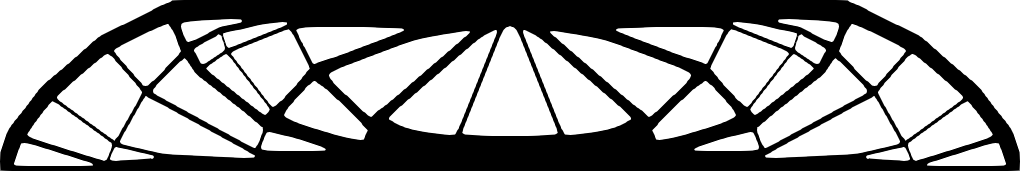}
&\includegraphics[width=4.29cm]{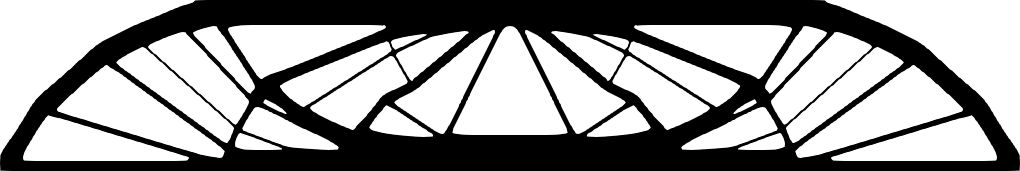}
&\includegraphics[width=4.29cm]{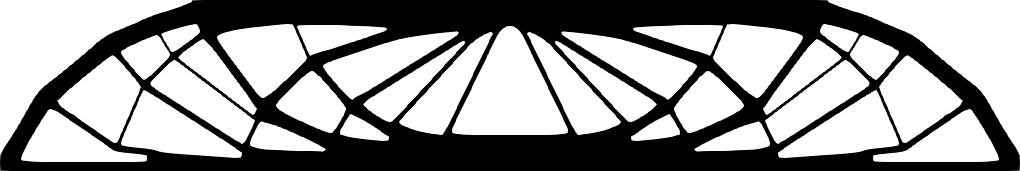}\\[-5pt] 
$240 \times 80$
&$\mathcal{S} = 99.95$   &$\mathcal{S} = 100.94$  &$\mathcal{S} = 100.18$ \\[5pt]

$l_{\refine} = 4$
&\includegraphics[width=4.29cm]{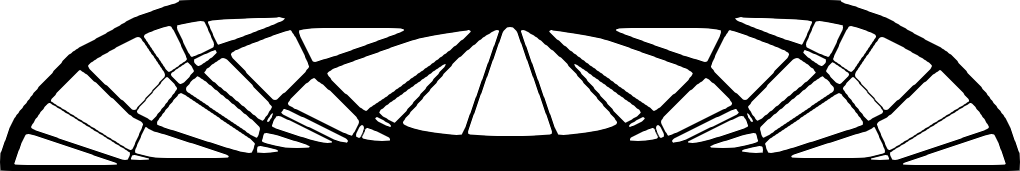}
&\includegraphics[width=4.29cm]{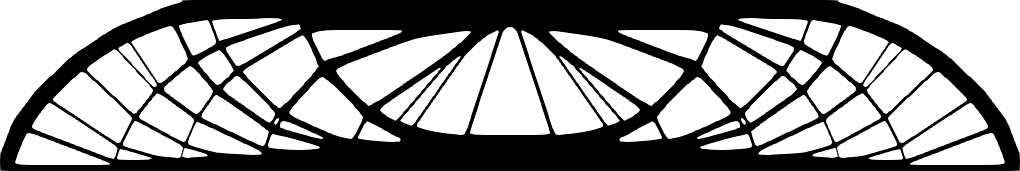}
&\includegraphics[width=4.29cm]{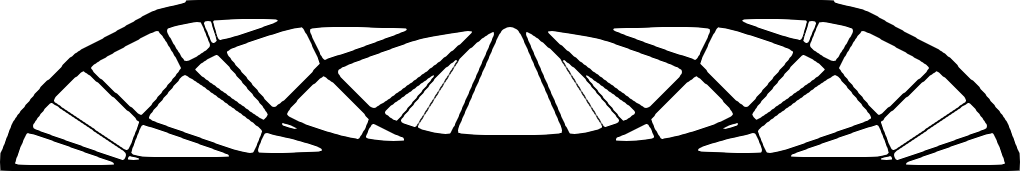}\\[-5pt] 
$480 \times 160$
&$\mathcal{S} = 99.67$   &$\mathcal{S} = 99.49$    &$\mathcal{S} = 99.86$ \\
\end{tabular}
\caption{2D beam using uniformly refined meshes with THB-splines and simultaneous hole seeding.}
\label{fig_2D_MBB_uniform_seeding}
\end{figure*}

Figure~\ref{fig_2D_MBB_uniform_seeding} shows the optimized designs generated on uniformly refined meshes with $l_{\refine} = 1, 2, 3$ and $4$ using THB-splines and considering different B-spline orders, along with the associated performance in terms of strain energy. The solid phase is depicted in black and the meshes are omitted as they are uniform. For coarse meshes, the designs exhibit different performances. The differences in the strain energy measure disappear as the meshes are refined. As expected, the complexity of the designs increases with mesh refinement. Remaining differences in both the geometry and the performance are likely caused by the non convexity of the optimization problem. 

\begin{figure*}[ht]\center
\renewcommand*{\arraystretch}{1.75}
\begin{tabular}{p{1.5cm}ccc}
& Linear B-splines & Quadratic B-splines & Cubic B-splines\\
$l_{\refine}^{0} = 1$
&\includegraphics[width=4.29cm]{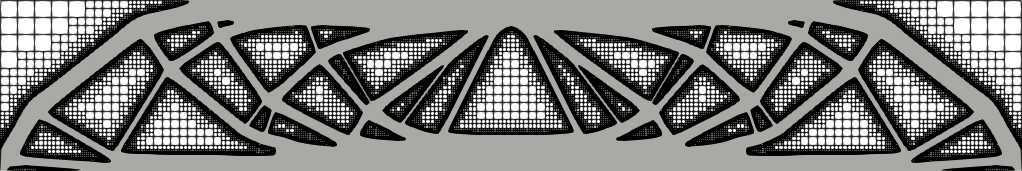}
&\includegraphics[width=4.29cm]{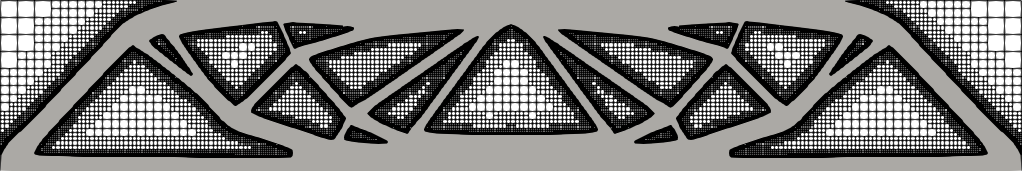}
&\includegraphics[width=4.29cm]{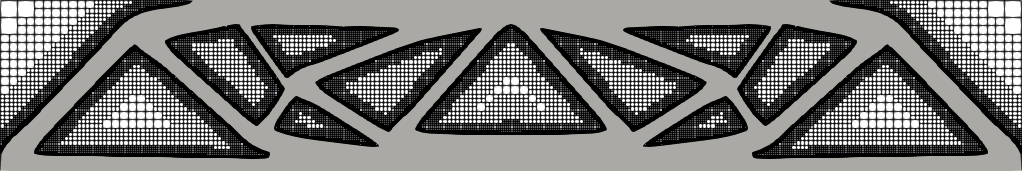}\\[-5pt]
$60 \times 20$
&$\mathcal{S} = 100.52$ &$\mathcal{S} = 102.25$  &$\mathcal{S} = 103.10$ \\[5pt]

$l_{\refine}^{0} = 2$
&\includegraphics[width=4.29cm]{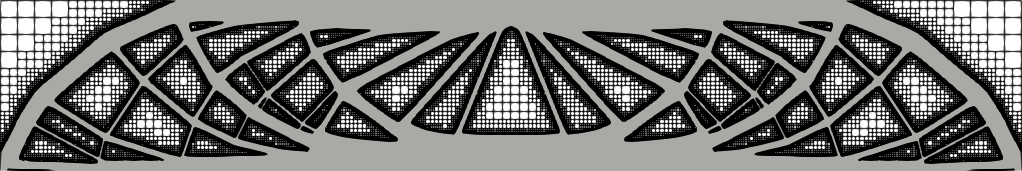}
&\includegraphics[width=4.29cm]{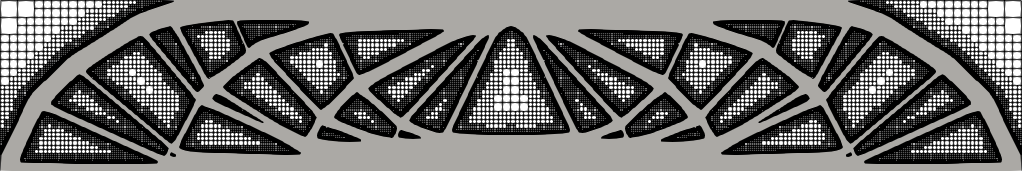}
&\includegraphics[width=4.29cm]{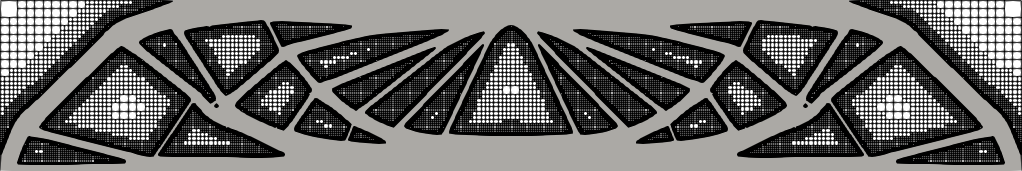}\\[-5pt]
$120 \times 40$
&$\mathcal{S} = 99.61$ &$\mathcal{S} = 99.73$  &$\mathcal{S} = 101.01$ \\[5pt]
  
$l_{\refine}^{0} = 3$
&\includegraphics[width=4.29cm]{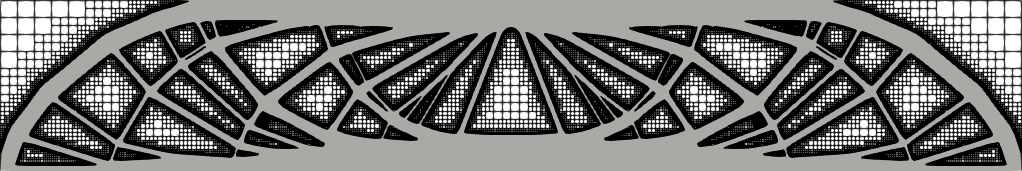}
&\includegraphics[width=4.29cm]{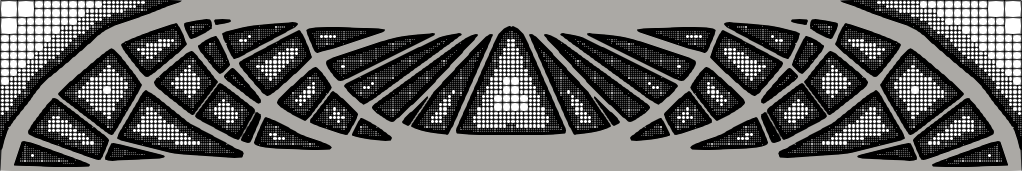}
&\includegraphics[width=4.29cm]{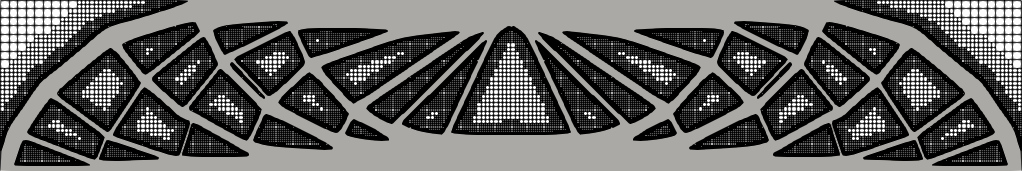}\\[-5pt] 
$240 \times 80$
&$\mathcal{S} = 99.44$   &$\mathcal{S} = 99.61$  &$\mathcal{S} = 99.86$ \\ 
\end{tabular}
\caption{2D beam using adaptively refined meshes  with THB-splines and simultaneous hole seeding. Initial refinement levels of $l_{\refine}^{0} = 1, 2$ and $3$ and maximum refinement level of $l_{\refine,\max} = 4$.}
\label{fig_2D_MBB_adaptive_seeding}
\end{figure*}

Figure~\ref{fig_2D_MBB_adaptive_seeding} shows the optimized designs generated on adaptively refined meshes with initial refinement level $l_{\refine}^{0} = 1, 2$ and $3$ and using linear, quadratic and cubic THB-splines. The performance in terms of strain energy values are provided for each design. For the adaptive designs, the mesh refinement is depicted in the void phase and the solid phase is represented in gray. The gains in terms of computational cost are given in Table~\ref{table_2D_MBB_seeding_compCost}.

Using the adaptive strategy, the results clearly show that, starting from coarser meshes, it is possible to recover designs rather similar, in terms of geometry and performance, to the ones obtained with finer uniform meshes, at least up to one or two refinement level higher. It is also noticeable that the layouts obtained with the refinement strategy depend on the initial mesh refinement. Finer structural members are more likely to develop and be maintained through the optimization process as the initial mesh is refined. 

The efficiency factors given in Table~\ref{table_2D_MBB_seeding_compCost}. They show that more significant gains are obtained when starting from coarser meshes, with lower refinement levels. The adaptation strategy allows for a reduction of the computational cost by more than a factor 3 over the uniform FEM model and by more than 1.75 over the uniform XFEM model. The achieved efficiency factors starting from the same initial refinement with different B-splines order match closely.

As a direct consequence of the buffer zone definition, see Subsection~\ref{bufferCriterion}, the refined meshes in Fig.~\ref{fig_2D_MBB_adaptive_seeding} exhibit a wider refined bandwidth as the interpolation order is increases from linear to cubic. As observed in the previous section and in Fig.~\ref{fig_2D_MBB_uniform_seeding} and ~\ref{fig_2D_MBB_adaptive_seeding}, linear B-splines can resolve finer structural members, while keeping a smooth description of the geometry at least in this 2D setting. Designs obtained with quadratic and cubic B-splines are smoother and similar.

\begin{table*}[ht]\center
\caption{Performance in terms of computational cost for designs in Fig.~\ref{fig_2D_MBB_adaptive_seeding}.}
\label{table_2D_MBB_seeding_compCost}
\renewcommand*{\arraystretch}{1.5}
\begin{tabular}{l|ccc|ccc|ccc}\hline
&\multicolumn{3}{c|}{$l_{\refine}^{0} = 1$}
&\multicolumn{3}{c|}{$l_{\refine}^{0} = 2$}
&\multicolumn{3}{c}{$l_{\refine}^{0} = 3$}\\ 
& Lin. & Quad. & Cub.
& Lin. & Quad. & Cub. 
& Lin. & Quad. & Cub.\\ \hline
$N_{\opt}$ 	&671 &550 &547 &941 &539 &433 &941 &484 &454\\ \hline
$E_{\xfem}$ 	&2.91 &3.58 &3.69 &2.15 &2.54 &2.70 &1.66 &1.79 &1.82\\
$R_{\xfem}$ 	&4.28 &4.98 &5.04 &3.13 &3.16 &3.42 &1.76 &1.81 &1.83\\
$E_{\fem}$  	&5.40 &6.52 &6.80 &4.13 &4.61 &4.80 &3.19 &3.20 &3.27\\
$R_{\fem}$ 	&4.28 &4.98 &5.04 &3.13 &3.16 &3.42 &1.76 &1.81 &1.83\\ \hline
\end{tabular}
\end{table*}

\begin{figure*}[ht]\center
\renewcommand*{\arraystretch}{1.75}
\begin{tabular}{p{1.75cm}ccc}
& Linear B-splines & Quadratic B-splines & Cubic B-splines\\
HB-splines
&\includegraphics[width=4.29cm]{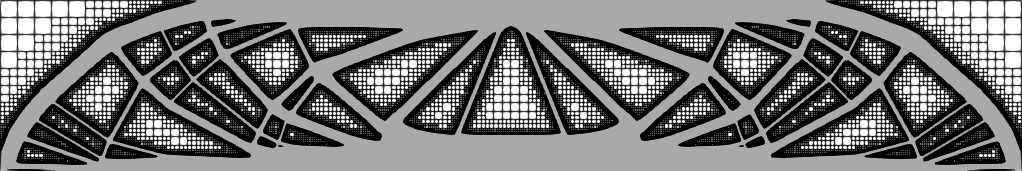}
&\includegraphics[width=4.29cm]{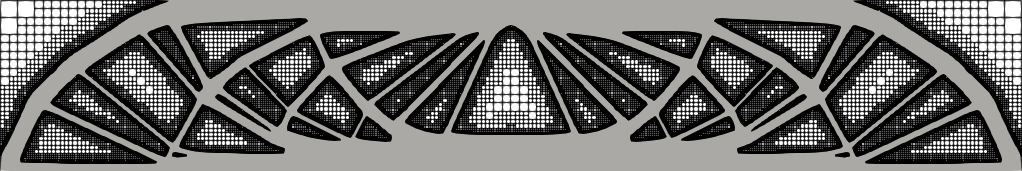}
&\includegraphics[width=4.29cm]{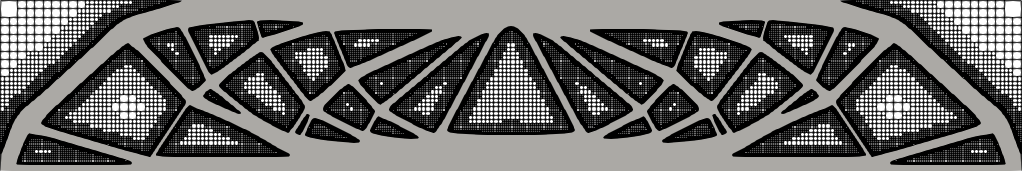}\\[-5pt]
$l_{\refine}^{0} = 2$
&$\mathcal{S} = 99.75$ &$\mathcal{S} = 99.85$   &$\mathcal{S} = 100.20$ \\[5pt] 
THB-splines
&\includegraphics[width=4.29cm]{figures/2DMBB/combined_Adaptive/Linear_Adap_25_Start_2_mesh}
&\includegraphics[width=4.29cm]{figures/2DMBB/combined_Adaptive/Quadratic_Adap_25_Start_2_mesh}
&\includegraphics[width=4.29cm]{figures/2DMBB/combined_Adaptive/Cubic_Adap_25_Start_2_mesh}\\[-5pt] 
$l_{\refine}^{0} = 2$
&$\mathcal{S} = 99.61$ &$\mathcal{S} = 99.73$  &$\mathcal{S} = 101.01$ \\ 
\end{tabular}
\caption{2D beam using adaptively refined meshes with simultaneous hole seeding. Initial refinement level of $l_{\refine}^{0} = 2$ and maximum refinement level of $l_{\refine,\max} = 4$. Comparison of HB- and THB-splines.}
\label{fig_2D_MBB_adaptive_seeding_truncation}
\end{figure*}

The influence of the truncation is explored by solving the problem with an adaptive mesh and both HB- and THB-splines on an initial mesh with initial refinement $l_{\refine}^{0} = 2$ and applying the refinement every 25 iterations. The obtained designs and corresponding final strain energy values are given in Fig.~\ref{fig_2D_MBB_adaptive_seeding_truncation}. The designs differ more significantly for HB- and THB-splines than when working with level set only, see Fig.~\ref{fig_2D_MBB_adaptive_truncation}. These differences can be partly explained by the maximum stencil size difference between HB- and THB-splines, i.e., 7, 16, and 37 for linear, quadratic and cubic HB-splines against 2, 9 and 16 for the THB-splines. On top of the stencil size mismatch, these differences in the final designs can be explained by the clipping operation applied to the density values to keep them between 0 and 1 when working with HB-splines. Clipping yields non differentiability with respect to the clipped values which can influence the optimization process. The computational gain factors are given in Table~\ref{table_2D_MBB_compCost_truncation} and match closely for both HB- and THB-splines.

\begin{table}[ht]\center
\caption{Performance in terms of computational cost for designs in Fig.~\ref{fig_2D_MBB_adaptive_seeding_truncation}.}
\label{table_2D_MBB_compCost_truncation}
\renewcommand*{\arraystretch}{1.5}
\begin{tabular}{l|ccc|ccc}\hline
&\multicolumn{3}{c|}{HB-splines}
&\multicolumn{3}{c}{THB-splines}\\ 
& Lin. & Quad. & Cub. 
& Lin. & Quad. & Cub. \\ \hline
$N_{\opt}$ 	&941 &405 &542 &941 &539 &433\\ \hline
$E_{\xfem}$ 	&2.16 &2.71 &2.52 &2.15 &2.54 &2.70\\
$R_{\xfem}$   	&3.31 &3.15 &3.29 &3.13 &3.16 &3.42\\
$E_{\fem}$ 	&4.14 &4.72 &4.65 &4.13 &4.61 &4.80\\
$R_{\fem}$   	&3.31 &3.15 &3.29 &3.13 &3.16 &3.42\\ \hline
\end{tabular}
\end{table}
A comparison of the results from Subsections~\ref{2DMBB_initial_hole} and~\ref{2DMBB_seeding} suggests that the ability to nucleate holes during the optimization process is crucial to increase the computational cost reduction achieved with mesh adaptation. Creating an initial hole pattern requires the use of rather fine initial meshes, which limits the computational cost savings. Comparing the efficiency ratios in Tables~\ref{table_2D_MBB_compCost} and~\ref{table_2D_MBB_seeding_compCost} shows that starting from a $120 \times 40$ initial mesh and refining every 25 iterations, an efficiency factor $E_{\xfem}$ of about $2.0$ is achieved for the initial seeding approach, while a ratio of about $2.5$ is obtained for the combined scheme. Moreover, although a dependency on the initial mesh refinement is observed in both cases, the combined level set/density scheme seems to mitigate this dependency.

The 2D examples, in Fig.~\ref{fig_2D_MBB_adaptive_truncation} and~\ref{fig_2D_MBB_adaptive_seeding_truncation}, de\-mon\-strated that using THB-splines does not affect the design geometry complexity, the performance or the computational gain achieved with the proposed adaptive strategy. THB-splines offer advantageous properties; they form a PU, a desirable property when solving finite element problems. They also lead to an improved conditioning of the system of equations and offer a convenient way to impose bounds on design variables, i.e., an operation largely applied when solving optimization problems. Therefore, all subsequent optimization problems are solved using THB-splines.

\subsection{Three dimensional beam}\label{3DBeam}
In this section, the previous 2D beam example is extended to 3D. This example allows studying the influence of the mesh adaptivity on 3D designs in terms of geometry, performance and computational cost. In particular, we assess the ability of our approach to generate well-known three dimensional optimal members such as shear webs that usually develop when solving optimization problems on fine uniform meshes. For 3D designs, the dependency on the initial hole pattern is mitigated by the ability to develop holes in the third dimension, which might counterbalance the performance issues encountered when using the initial seeding approach with the refinement strategy.

\begin{figure}[h]\center
\includegraphics[scale=1]{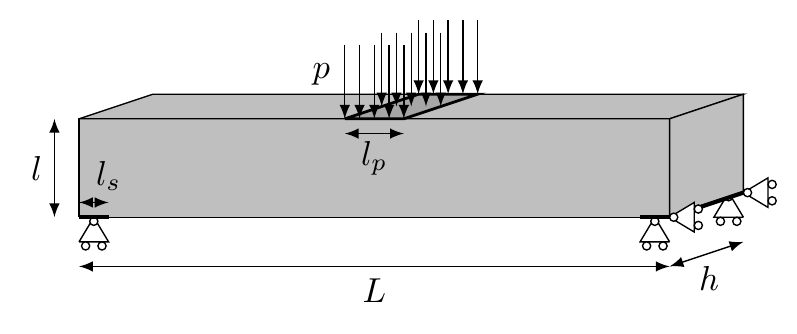}
\caption{Three dimensional beam.}
\label{fig_3D_MBB}
\end{figure}

Figure~\ref{fig_3D_MBB} illustrates the setup for the solid-void 3D beam with a length $L = 6.0$, a height $l = 1.0$ and a thickness $h = 1.0$. The structural response is modeled by an isotropic linear elastic model with a Young's modulus $E = 1.0$ and a Poisson's ratio $\nu = 0.3$. The structure is supported on its two lower sides over a length $l_s = 0.025$ and a pressure $p = 1.0$ is applied in the middle of its top span over an area $l_p \times h = 2 \times 0.025 \times 1.0$. The values above are provided in self-consistent units.

We optimize the beam for minimum compliance under a mass constraint of $10 \%$ of the total mass of the design domain. Taking advantage of the symmetry, we only study a quarter of the design domain. A list of problem parameters is given in Table~\ref{table_3DMBB_parameter}. 
\begin{table}[h]\center
\caption{Parameter list for the 3D beam optimization problem.}
\label{table_3DMBB_parameter}
\renewcommand{\arraystretch}{1.5}
\begin{tabular}{p{0.45\columnwidth}p{0.45\columnwidth}}\hline
Parameter & Value\\ \toprule
$w_s$ & 9.0\\
$c_p$ & 1.0\\
$c_{\phi}$ & 10.0\\
$c_v$ & 0.4\\
$P_0$ & 10.0\\
$S_0$ & 59.12\\
$\phi_{\scale}$ & $3h_{init}$\\
$V_{max}$ & 3\\
$\gamma_N$ & $100E/h$\\
$\gamma_G$ & $0.01$\\
$\phi_{\target}$ &$1.5h_{\init}$\\
$\nabla \phi_{\target}$ &0.75\\ 
$\mathcal{I}_{\max}$ & 1\\
$\gamma_{I}$ & 4.61\\
$\alpha_1 = \alpha_2 = \alpha_3$ & 0.5\\ \hline
\end{tabular}
\end{table}

As for the 2D case, the initial hole seeding and the combined level set/density schemes are studied. The optimization problem, in Eq.~\eqref{Eq_optFormulation}, is solved on uniformly and adaptively refined meshes considering four levels of refinements, $l_{\refine} = 0, 1, 2, 3$, corresponding to $18 \times 6 \times 3$, $36 \times 12 \times 6$, $72 \times 24 \times 12$, and $144 \times 48 \times 24$ element meshes. The influence of the B-spline interpolation order is investigated by solving the same problem with linear, quadratic and cubic B-spline functions. All the simulations are performed with THB-splines.

\subsubsection{Design with initial hole seeding}\label{3DMBB_initial_hole}
\begin{figure*}[ht]\center
\includegraphics[width=0.4\columnwidth]{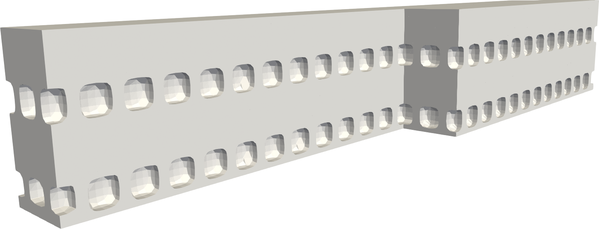}
\includegraphics[width=0.4\columnwidth]{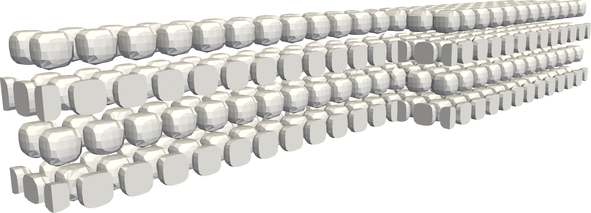}
\caption{Initial hole seeding for the 3D beam. Representation of the solid on the left and the void on the right.}
\label{fig_3D_MBB_initial_hole}
\end{figure*}

The initial seeding of holes in the design domain is presented in Fig.~\ref{fig_3D_MBB_initial_hole}. This configuration satisfies the mass constraint, but requires a rather fine initial mesh to resolve the hole pattern, i.e., $l_{\refine} > 1$.

\begin{figure*}[ht]\center
\begin{tabular}{lccc}
& Linear B-splines & Quadratic B-splines & Cubic B-splines\\
\multirow{ 2}{*}{$\begin{array}{c} l_{\refine} = 2\\[5pt] 72 \times 24 \times 12\\ \end{array}$} 
&\includegraphics[width=4.09cm]{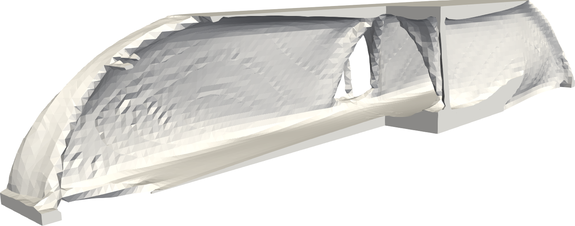}
&\includegraphics[width=4.09cm]{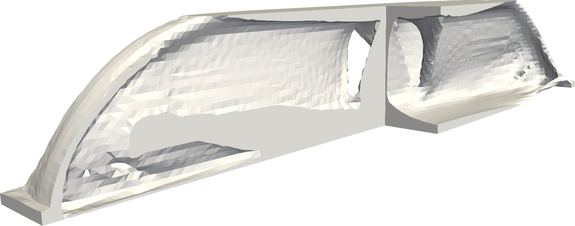}
&\includegraphics[width=4.09cm]{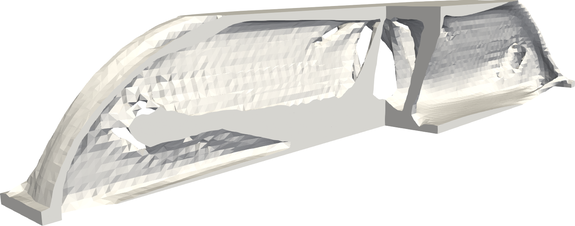}\\
&\includegraphics[width=4.09cm]{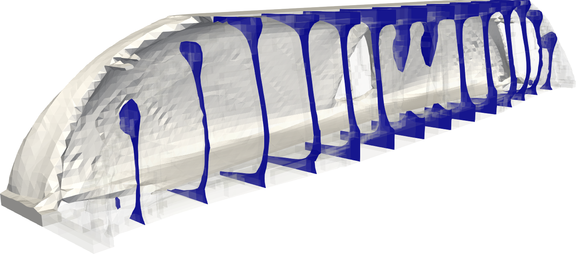}
&\includegraphics[width=4.09cm]{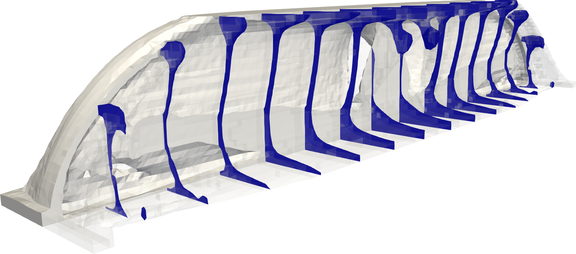}
&\includegraphics[width=4.09cm]{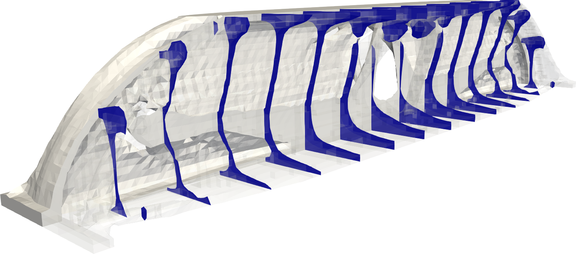}\\
&$\mathcal{S} = 77.20$  &$\mathcal{S} = 78.12$  &$\mathcal{S} = 78.24$\\[10pt] 

\multirow{ 2}{*}{$\begin{array}{c} l_{\refine} = 3\\[5pt] 144 \times 48 \times 24\\ \end{array}$} 
&\includegraphics[width=4.09cm]{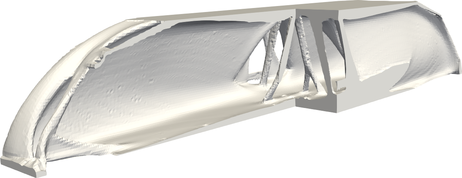}
&\includegraphics[width=4.09cm]{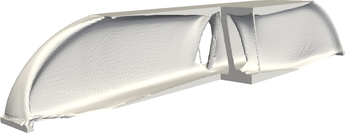}
&\includegraphics[width=4.09cm]{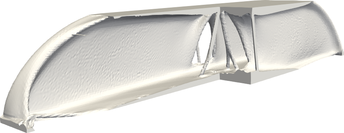}\\
&\includegraphics[width=4.09cm]{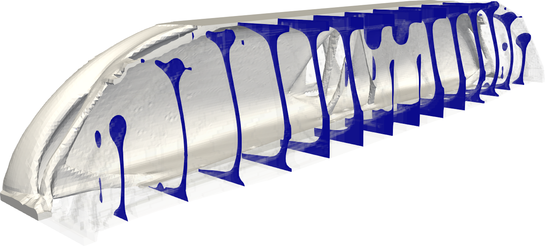}
&\includegraphics[width=4.09cm]{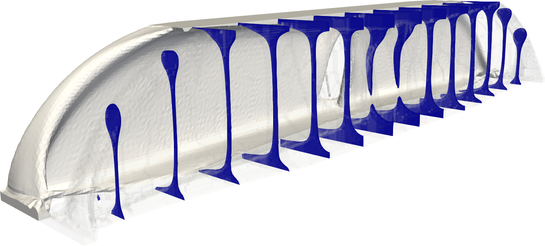}
&\includegraphics[width=4.09cm]{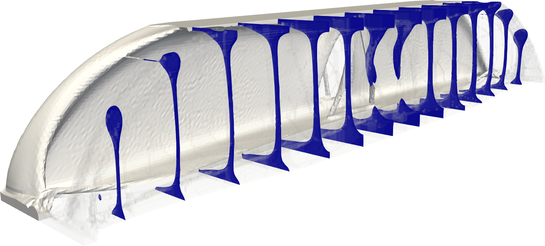}\\
&$\mathcal{S} = 77.50$  &$\mathcal{S} = 76.99$  &$\mathcal{S} = 77.14$\\
\end{tabular}
\caption{3D beam using uniformly refined meshes with THB-splines and initial hole seeding.}
\label{fig_3D_MBB_uniform}
\end{figure*}

\begin{figure*}[ht]\center
\begin{tabular}{p{1.5cm}ccc}
& Linear B-splines & Quadratic B-splines & Cubic B-splines\\
\multirow{ 2}{*}{$\begin{array}{c} l_{\refine}^{0} = 2\\[5pt] 72 \times 24 \times 12\\ \end{array}$} 
&\includegraphics[width=4.29cm]{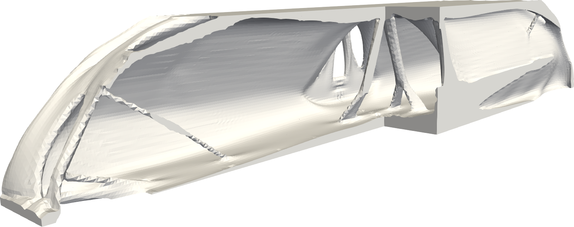}
&\includegraphics[width=4.29cm]{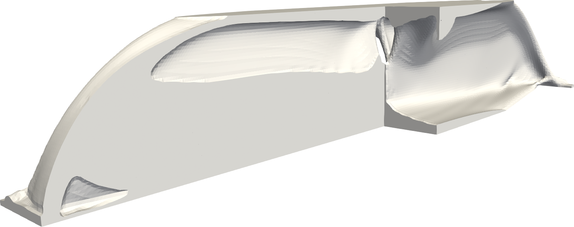}
&\includegraphics[width=4.29cm]{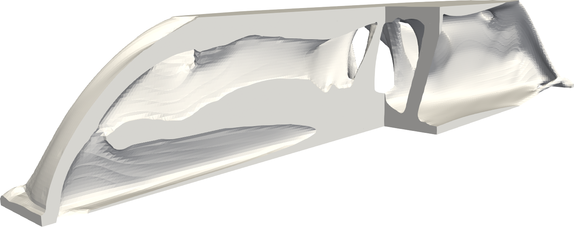}\\
&\includegraphics[width=4.29cm]{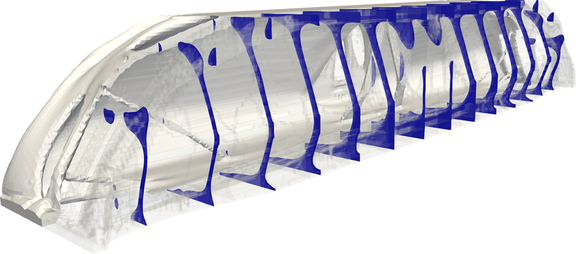}
&\includegraphics[width=4.29cm]{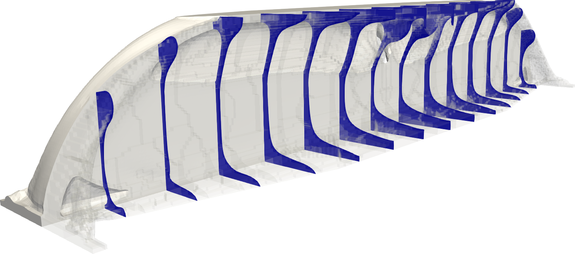}
&\includegraphics[width=4.29cm]{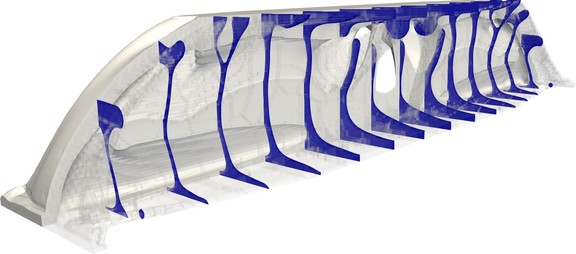}\\
&$\mathcal{S} = 77.82$  &$\mathcal{S} = 79.19$  &$\mathcal{S} = 78.41$ \\
\end{tabular}
\caption{3D beam using adaptively refined meshes with THB-splines and initial hole seeding. Initial refinement levels of $l_{\refine}^{0} = 2$ and maximum refinement level of $l_{\refine,\max} = 3$.}
\label{fig_3D_MBB_adaptive}
\end{figure*}

The optimized designs generated on uniformly and adaptively refined meshes are given in Fig.~\ref{fig_3D_MBB_uniform} and~\ref{fig_3D_MBB_adaptive}, respectively, along with their performance in terms of strain energy. Additionally, we depict vertical cross sections of the designs in blue, to visualize the layouts' internal material arrangement. The efficiency factors are summarized in Table~\ref{table_3D_MBB_compCost}.

In Figure~\ref{fig_3D_MBB_uniform}, uniform meshes with levels of refinement $l_{\refine} = 2$ and $3$ are considered, as coarser meshes cannot resolve the initial hole pattern. The design geometries and their performance in terms of strain energy are quite similar for linear B-splines. Each optimized layout is characterized by two shear webs and a truss-type structure that develops under the loaded surface. For higher B-spline orders, similar designs are observed only when using a finer mesh, i.e., $l_{\refine} = 3$. For coarser meshes with refinement level $l_{\refine} = 2$, a single shear web develops in the center of the design. This behavior suggests that the obtained design depends on both the B-spline order and the refinement level. As the mesh is refined, smaller members are included in the layouts, which exhibit smoother surfaces. This can be observed in particular for the shear webs that become thin and smooth with mesh refinement.

Applying the refinement strategy every 25 iterations, the optimization process is started on a uniform mesh with a refinement level $l_{\refine}^{0} = 2$, see Figure~\ref{fig_3D_MBB_adaptive}. In the course of the optimization process, the mesh is locally refined up to level $l_{\refine,\max} = 3$ or coarsened to a level $l_{\min} = 0$. The adaptive strategy generates designs that display a slightly higher geometric complexity and smoothness than the designs on corresponding initial uniform me\-shes. However, they are visibly dependent on the initial mesh refinement as the obtained design topologies, i.e., two shear webs for linear B-splines versus one shear web for higher order B-splines, are similar to the ones generated on a uniform mesh with $l_{\refine} = 2$. The efficiency ratios are given in Table~\ref{table_3D_MBB_compCost}. They show that a moderate reduction of the computational cost is achieved as the optimization process with adaptivity is started on a rather fine mesh, a reduction by a factor of about 4 over the uniform FEM model and by about 1.5 over the uniform XFEM model.

It is interesting to note that, as for the 2D case, linear B-splines support the development of smaller structural features such as the beams observed for both the uniform and the adaptive designs. At the same time, using these low order B-splines leads to spurious oscillations of the level set field, which result in layouts with rough surfaces. Higher order B-splines promote the generation of smoother designs, as is the case for both quadratic and cubic B-splines.

\begin{table}[ht]\center
\caption{Performance in terms of computational cost for designs in Fig.~\ref{fig_3D_MBB_adaptive}.}
\label{table_3D_MBB_compCost}
\renewcommand*{\arraystretch}{1.5}
\begin{tabular}{l|ccc}\hline
&\multicolumn{3}{c}{$l_{\refine}^{0} = 2$}\\ 
& Lin. & Quad. & Cub. \\ \hline
$N_{\opt}$ 	&501 &520 &520\\ \hline
$E_{\xfem}$ 	&1.26 &1.64&1.43\\
$R_{\xfem}$    	&2.07 &3.43 &2.32\\
$E_{\fem}$     	&3.90 &5.44 &4.64\\
$R_{\fem}$    	&2.72 &4.30 &2.91\\ \hline
\end{tabular}
\end{table}

\subsubsection{Design with simultaneous hole seeding}\label{3DMBB_seeding}
In this section, holes are nucleated in the design domain during the optimization process using the combined level set/density scheme presented in Section~\ref{Hole}. For the combined scheme, the SIMP exponent $\beta_{\SIMP}$ is set to $2.0$. The density shift parameter $\rho_{\shift}$ is initially set to $0.2$ and is updated every 25 optimization steps. Contrary to Section~\ref{3DMBB_initial_hole}, coarser meshes can be used as a starting point for the adaptive optimization, since the initial design geometry is not complex.

\begin{figure*}[ht]\center
\begin{tabular}{p{1.7cm}ccc}
& Linear B-splines & Quadratic B-splines & Cubic B-splines\\
\multirow{ 2}{*}{$\begin{array}{c} l_{\refine}^{0} = 1\\[5pt] 36 \times 12 \times 6\\ \end{array}$} 
&\includegraphics[width=4.19cm]{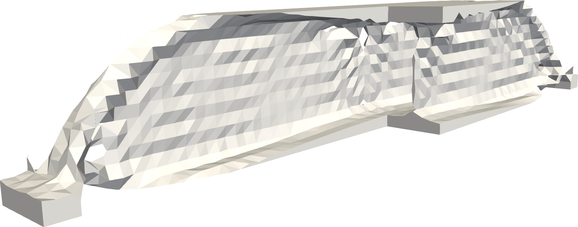}
&\includegraphics[width=4.19cm]{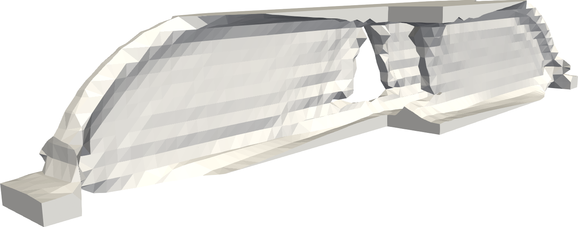}
&\includegraphics[width=4.19cm]{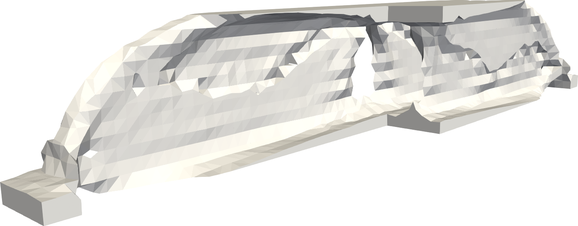}\\
&\includegraphics[width=4.19cm]{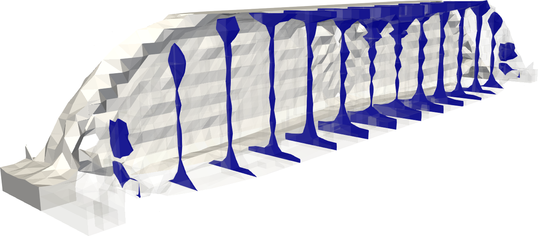}
&\includegraphics[width=4.19cm]{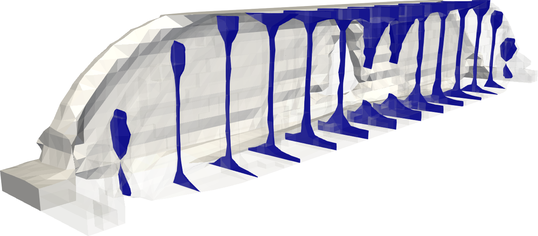}
&\includegraphics[width=4.19cm]{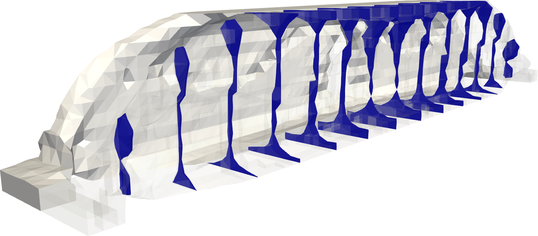}\\
&$\mathcal{S} = 70.84$ &$\mathcal{S} = 71.20$  &$\mathcal{S} = 71.15$ \\[10pt]
 
\multirow{ 2}{*}{$\begin{array}{c} l_{\refine}^{0} = 2\\[5pt] 72 \times 24 \times 12\\ \end{array}$} 
&\includegraphics[width=4.19cm]{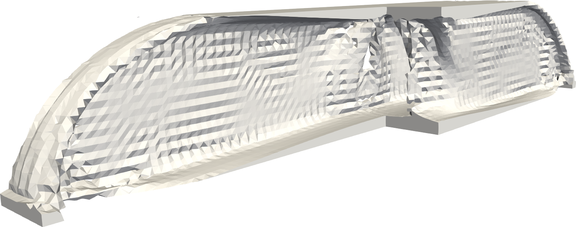}
&\includegraphics[width=4.19cm]{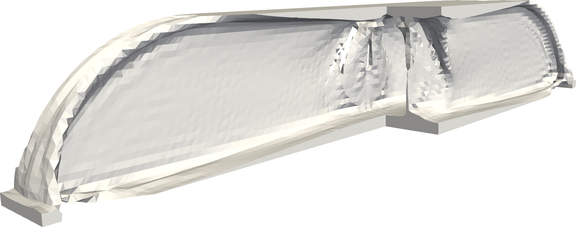}
&\includegraphics[width=4.19cm]{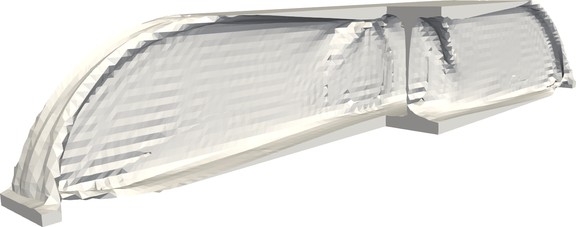}\\
&\includegraphics[width=4.19cm]{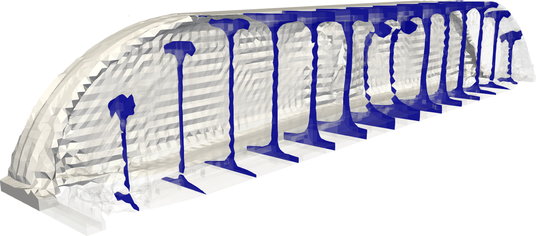}
&\includegraphics[width=4.19cm]{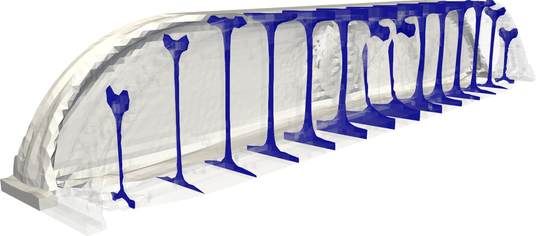}
&\includegraphics[width=4.19cm]{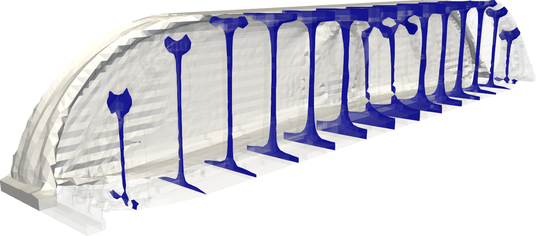}\\
&$\mathcal{S} = 76.50$ &$\mathcal{S} = 76.44$  &$\mathcal{S} = 76.40$ \\[10pt] 

\multirow{ 2}{*}{$\begin{array}{c} l_{\refine}^{0} = 3\\[5pt] 144 \times 48 \times 24\\ \end{array}$} 
&\includegraphics[width=4.19cm]{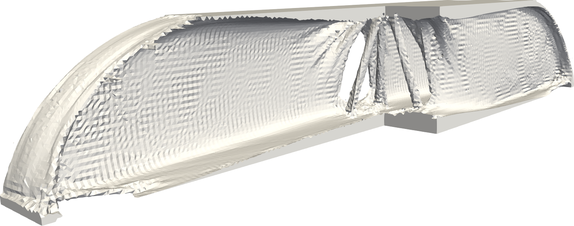}
&\includegraphics[width=4.19cm]{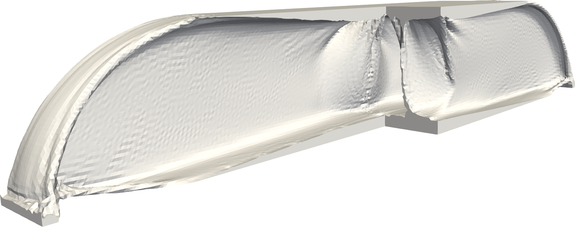}
&\includegraphics[width=4.19cm]{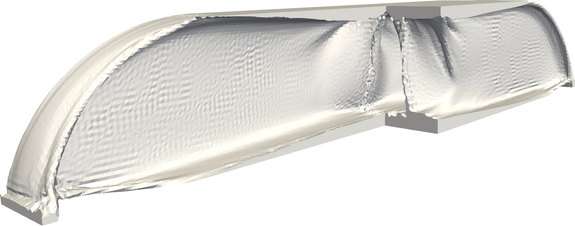}\\
&\includegraphics[width=4.19cm]{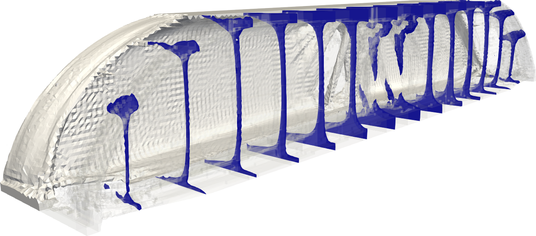}
&\includegraphics[width=4.19cm]{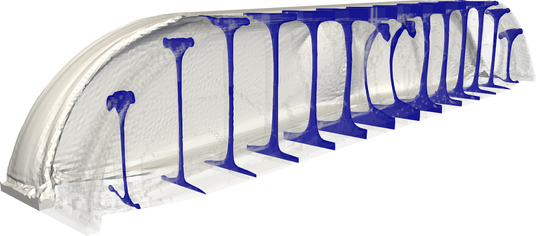}
&\includegraphics[width=4.19cm]{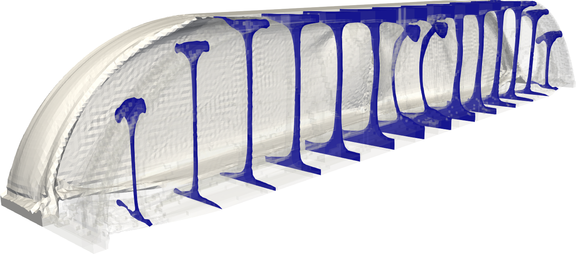}\\
&$\mathcal{S} = 76.98$ &$\mathcal{S} = 76.77$  &$\mathcal{S} = 76.84$ \\
\end{tabular}
\caption{3D beam using uniformly refined meshes with THB-splines and simultaneous hole seeding.}
\label{fig_3D_MBB_uniform_seeding}
\end{figure*}

\begin{figure*}[ht]\center
\begin{tabular}{p{1.7cm}ccc}
& Linear B-splines & Quadratic B-splines & Cubic B-splines\\
\multirow{ 2}{*}{$\begin{array}{c} l_{\refine}^{0} = 1\\[5pt] 36 \times 12 \times 6\\ \end{array}$} 
&\includegraphics[width=4.29cm]{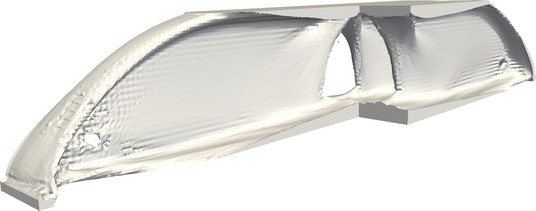}
&\includegraphics[width=4.29cm]{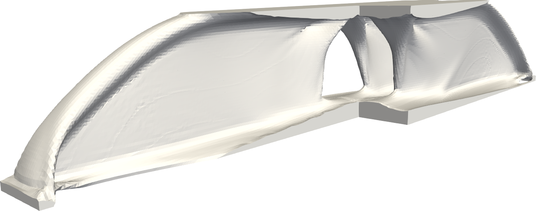}
&\includegraphics[width=4.29cm]{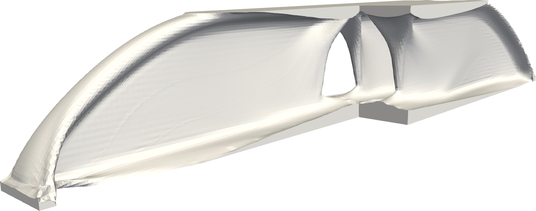}\\
&\includegraphics[width=4.29cm]{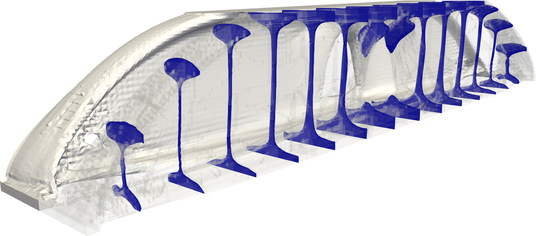}
&\includegraphics[width=4.29cm]{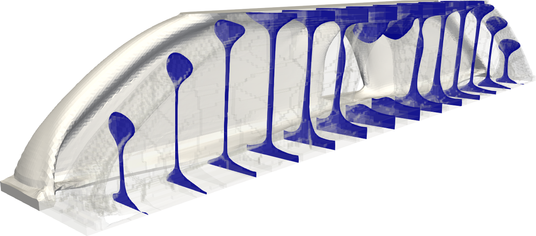}
&\includegraphics[width=4.29cm]{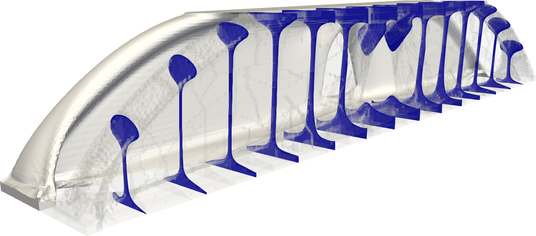}\\
&$\mathcal{S} = 76.19$ &$\mathcal{S} = 74.94$  &$\mathcal{S} = 76.37$ \\
\end{tabular}
\caption{3D beam using adaptively refined meshes with simultaneous hole seeding. Initial refinement levels of $l_{\refine}^{0} = 1$ and maximum refinement level of $l_{\refine,\max} = 3$.}
\label{fig_3D_MBB_adaptive_seeding}
\end{figure*}

The designs obtained on uniformly and adaptively refined meshes are shown in Fig.~\ref{fig_3D_MBB_uniform_seeding} and~\ref{fig_3D_MBB_adaptive_seeding} with their performance in terms of strain energy. Additionally, we depict vertical cross sections of the designs in blue, to visualize the layouts' internal material arrangement. Table~\ref{table_3D_MBB_compCost_seeding} gives the design performance in terms of computational cost.

In Fig.~\ref{fig_3D_MBB_uniform_seeding}, three levels of refinement are considered $l_{\refine} = 1, 2$ and $3$. The uniform designs geometries are quite similar, as each layout is characterized by two shear webs and a truss-type structure that develops under the loaded surface. Also, the performance in terms of strain energy is quite similar for the designs obtained on meshes with refinement levels $l_{\refine} = 2$ and $3$. These results demonstrate that the combined level set/density scheme reduces the dependency on the initial design that was observed with the initial hole seeding approach. Using extremely coarse meshes with $l_{\refine} = 1$ lead to spurious oscillations of the level set field and thus a rough solid/void interface surface. This behavior is more stron\-gly pronounced when using linear B-splines.

Working with an adaptive mesh, the optimization process is started on a uniformly refined mesh with the lowest level of refinement $l_{\refine}^{0} = 1$ and is refined up to $l_{\refine,\max} = 3$ or coarsened to $l_{\min} = 0$, see Fig.~\ref{fig_3D_MBB_adaptive_seeding}. The results demonstrate the ability of the adaptive optimization approach to create designs similar to the ones obtained on uniform meshes with higher level of refinement. The designs also compare well in terms of strain energy performance. The refinement strategy leads to a significant reduction of the computational cost, as the optimization is started on a coarser mesh. Table~\ref{table_3D_MBB_compCost_seeding} shows high efficiency ratios, i.e., a reduction by a factor of about 7 over the uniform FEM model and by about 3 over the uniform XFEM model.

As observed in the previous section, using low order B-splines promotes the development of thinner features, which is not observed with higher polynomial orders. However, linear B-spline designs lack smoothness, even as the mesh is refined. 

\begin{table}[ht]\center
\caption{Performance in terms of computational cost for designs in Fig.~\ref{fig_3D_MBB_adaptive_seeding}.}
\label{table_3D_MBB_compCost_seeding}
\renewcommand*{\arraystretch}{1.5}
\begin{tabular}{l|ccc}\hline
&\multicolumn{3}{c}{$l_{\refine}^{0} = 1$}\\ 
& Lin. & Quad. & Cub. \\ \hline
$N_{\opt}$ 	&292 &209 &275\\ \hline
$E_{\xfem}$ 	&2.87 &4.08 &3.09\\ 
$R_{\xfem}$ 	&3.78 &3.75 &3.92\\
$E_{\fem}$ 	&6.70 &8.78 &7.56\\
$R_{\fem}$ 	&3.78 &3.75 &3.92\\ \hline
\end{tabular}
\end{table}

The effect of hole seeding during the optimization process is not as pronounced as in 2D. In 3D holes can be created more easily by effectively changing the thickness, i.e., changing the shape in the third dimension. The main advantage of the combined scheme lies in the ability to use coarser meshes as starting point for the adaptive topology optimization, which leads to larger efficiency ratios $E_{\xfem}$, i.e., around $3$ for the combined scheme versus around $1.5$ with the initial seeding strategy when compared to uniform XFEM meshes.

Both the 2D and 3D beam examples show that using linear B-splines allows for the development of thinner structural members. However, in the 3D case, the latter details develop at the cost of spurious oscillations in the level set field, leading to a rough representation of the surface. Higher order B-splines promote the development of smoother designs. In 2D and in 3D, the designs generated with quadratic and cubic B-splines match closely. Therefore, all upcoming simulations are carried out with quadratic B-splines. They lead to smooth designs without reducing the complexity of the geometry achieved. Moreover, they present smaller support and a lower level of complexity than cubic B-splines, which leads to simpler implementation and a reduced overall computational cost.

\subsection{Short beam under uniform pressure}\label{TableDesign}
This final example aims at demonstrating the capabilities of the proposed mesh adaptivity to generate geometrically complex designs with limited computational resources. The optimization problem formulation is modified to generate designs with increasingly lower volume fractions and show that mesh adaptation allows for the development of extremely thin structural members.

\begin{figure}[h]\center
\includegraphics[scale=1]{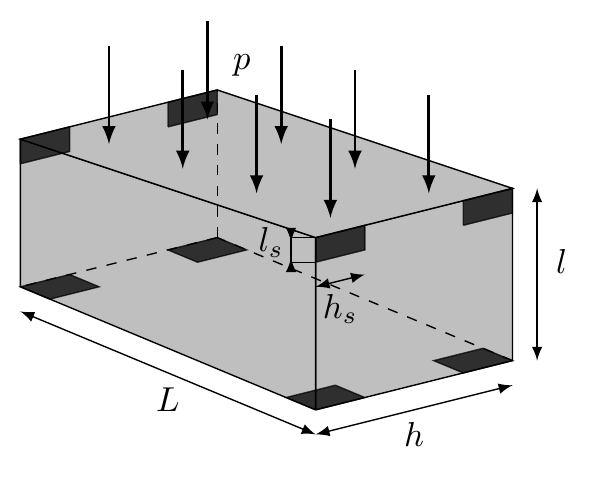}
\caption{Three dimensional short beam under uniform pressure.}
\label{fig_table}
\end{figure}

We consider a short beam subject to a uniform pressure on its top. The design domain with boundary conditions is shown in Fig.~\ref{fig_table} and is defined by a rectangular cuboid with a length $L = 2.0$, a height $l = 0.6$ and a thickness $h = 1.0$. The structural response is described by an isotropic linear elastic model with a Young's modulus $E = 6.9\ 10^{10}$ and a Poisson's ratio $\nu = 0.3$. The presence of material is enforced on the top of the table, where a uniform pressure load $p = 0.67\ 10^{6}$ is applied. The structure is supported on the side at the top corners, and on the bottom at each corner with $l_s = 0.01$ and $h_s = 0.2$.

We optimize the structure for minimum compliance and mass, no constraints are imposed. Applying the adaptive refinement strategy, we focus on the influence of the mass of the designs on the achieved computational cost reduction. Taking advantage of the symmetry, we only study a quarter of the design domain. The problem parameters are given in Table~\ref{table_3DTable_parameter} in self-consistent units. 
\begin{table}[h]\center
\caption{Parameter list for the 3D short beam under uniform pressure design optimization problem.}
\label{table_3DTable_parameter}
\renewcommand{\arraystretch}{1.5}
\begin{tabular}{p{0.45\columnwidth}p{0.45\columnwidth}}
\hline
Parameter & Value\\ \toprule
$w_s^{0}$ & 0.5\\
$w_m$ & 10.0\\
$c_p$ &0.05\\
$c_{\phi}$ &5.0\\
$c_v$ & 1.0\\
$\mathcal{S}^0$ & 61.93\\
$\mathcal{M}^0$ & 0.3\\
$\mathcal{P}^0$ & 3\\
$\phi_{\scale}$ & $3h_{init}$\\
$\gamma_N$ & $100\ E/h$\\
$\gamma_G$ & $0.01$\\
$\phi_{\target}$ &$1.5\ h_{\init}$\\
$\nabla \phi_{\target}$ &0.75\\ 
$\mathcal{I}_{\max}$ & 1\\
$\gamma_{I}$ & 4.61\\
$\alpha_1 = \alpha_2 = \alpha_3$ & 0.5\\ \hline
\end{tabular}
\end{table} 

We investigate the effect of the designs mass on the adaptive strategy. To generate designs with various strain energy/mass ratios, we vary the compliance weighting parameter $w_s$ value. To prevent an overly aggressive removal of mass early on in the optimization process, we apply a continuation on $w_s$, which is initially set to $w_s^{0} = 0.5$ and is gradually decreased during the optimization process to the prescribed values $w_s^{\final} = 0.1, 0.05, 0.005$ and $0.0005$. We should note that the strategy for reducing $w_s$ impact the computational costs. Large $w_s$ values put more emphasis on the strain energy. Thus, stiffer structures, using more mass, are generated, which in turn requires more elements to be resolved at a particular refinement level.

Based on the observations made for the 2D and 3D beam problems in Sections~\ref{2DBeam} and~\ref{3DBeam}, we study this problem with the level set/density combined scheme only, as it allows to initiate the design process on a rather coarse mesh. The SIMP exponent $\beta_{\SIMP}$ is set to $1.5$. The density shift parameter $\rho_{\shift}$ is initially set to $0.2$ and is updated every 25 optimization steps. Furthermore the design variable field is interpolated with quadratic truncated hierarchical B-splines. We consider four levels of refinement $l_{\refine} = 0, 1, 2,$ and $3$, corresponding to $24 \times 14 \times 12$, $48 \times 28 \times 24$, $96 \times 56 \times 48$ and $192 \times 112 \times 96$ element meshes. The problem is solved using adaptively refined meshes, starting from an initial mesh with a refinement level $l_{\refine}^{0} = 0$ and allowing refinement up to $l_{\refine,\max} = 3$ and coarsening to $l_{\min} = 0$. Solving the problem on uniformly refined XFEM meshes is omitted.

\begin{figure*}[ht]\center
\begin{tabular}{m{0.15\textwidth}m{0.375\textwidth}m{0.375\textwidth}} 
$\begin{array}{lll} w_s^{\final}& =&  0.1\\[5pt] l_{\refine}^{0} &=& 0\\ \end{array}$
&\includegraphics[width=6cm]{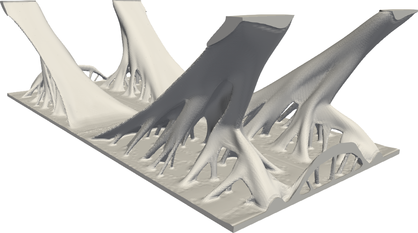}
&\includegraphics[width=6cm]{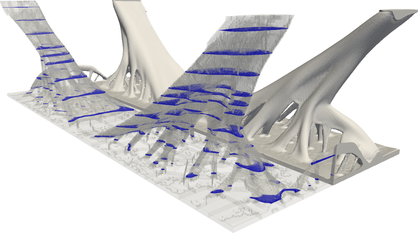}\\

$\begin{array}{lll} w_s^{\final}& =&  0.05\\[5pt] l_{\refine}^{0} &=& 0\\ \end{array}$
&\includegraphics[width=6cm]{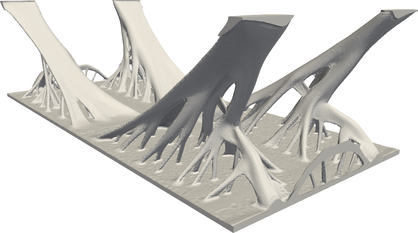}
&\includegraphics[width=6cm]{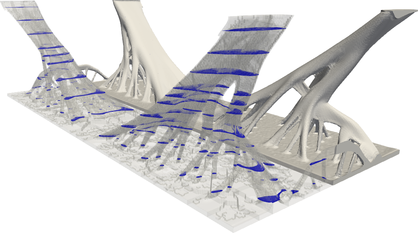}\\

$\begin{array}{lll} w_s^{\final}& =&  0.005\\[5pt] l_{\refine}^{0} &=& 0\\ \end{array}$
&\includegraphics[width=6cm]{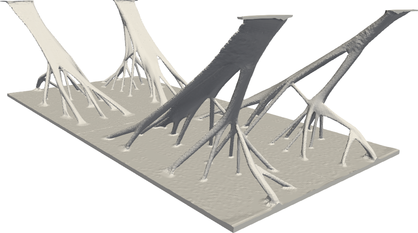}
&\includegraphics[width=6cm]{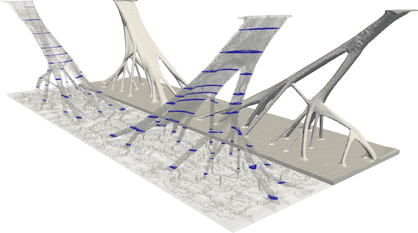}\\

$\begin{array}{lll} w_s^{\final}& =&  0.0005\\[5pt] l_{\refine}^{0} &=& 0\\ \end{array}$
&\includegraphics[width=6cm]{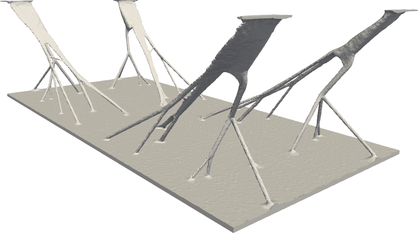}
&\includegraphics[width=6cm]{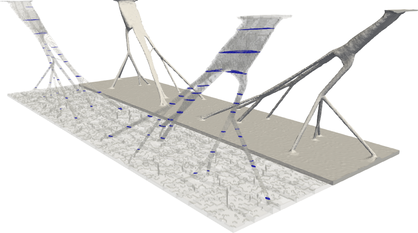}\\
\end{tabular}
\caption{Three dimensional short beam under uniform pressure design using adaptively refined mesh with quadratic THB-splines and simultaneous hole seeding. Initial refinement level of $l_{\refine}^{0} = 0$ and maximum refinement level of $l_{\refine,\max} = 3$.}
\label{fig_table_adaptive_seeding}
\end{figure*}

The designs generated on adaptively refined meshes with decreasing volume fractions are shown in Fig.~\ref{fig_table_adaptive_seeding}. Additionally, we depict horizontal cross sections of the designs in blue, to visualize the layouts' internal material arrangement. The performance in terms of strain energy, mass and the volume ratios are given in Table~\ref{fig_table_seeding_objective}. As the mass of the design decreases, very thin structural members are generated. These designs demonstrate the ability of the adaptive mesh strategy to resolve fine members and details starting from rather coarse meshes.

Table~\ref{table_table_seeding_compCost} gives the iteration $k_{\peak}$ at which the peak number of free DOFs is reached, the peak number of free DOFs and the number of free DOFs at the end of the optimization process. The peak number of DOFs is similar for each problem setup and is recorded around iteration 95, i.e. after the second refinement step. The final number of unconstrained DOFs show that as the mass of the designs becomes smaller, the associated number of free DOFs in the system drops, leading to a reduced computational time. 

Table~\ref{table_table_seeding_compCost} also provides the efficiency factors evaluated with $N_{opt}$ being the number of iterations to convergence. Here, we only consider the efficiency factors with respect to an equivalent FEM model, i.e., $E_{\fem}$ and $R_{\fem}$ considering a uniformly refined mesh with $l_{\refine} = 3$. The efficiency ratios $E_{\fem}$ increase as the mass of the designs decreases, i.e., the lower the volume fraction, the larger the computational saving. Since all cases start with $w_s = 0.5$ and since the peak DOFs count is recorded around iteration 95 for each design, the $R_{\fem}$ values are almost the same for all cases. The strategy for reducing $w_s$ could be improved to increase the computational efficiency due to mesh adaptation.

\begin{table}[ht]\center
\caption{Performance in terms of strain energy, mass and volume fraction for designs in Fig.~\ref{fig_table_adaptive_seeding}.}
\label{fig_table_seeding_objective}
\renewcommand*{\arraystretch}{1.5}
\begin{tabular}{c|c|c|c}\hline
$w_s^{\final}$ & Strain energy & Mass & Volume fraction \\ \hline
0.1       &\ 349.38  &0.0232 &5.57\\ 
0.05     &\ 469.24  &0.0189 &4.14\\ 
0.005   &1442.14  &0.0115 &1.68\\ 
0.0005 &4421.07  &0.0098 &1.19\\ \hline
\end{tabular}
\end{table}

\begin{table*}[ht]\center
\caption{Performance in terms of computational cost for designs in Fig.~\ref{fig_table_adaptive_seeding}.}
\label{table_table_seeding_compCost}
\renewcommand*{\arraystretch}{1.5}
\begin{tabular}{l|c|c|c|c}\hline
& $w_s^{\final} = 0.1$ & $w_s^{\final} = 0.05$ & $w_s^{\final} = 0.005$ & $w_s^{\final} = 0.0005$\\ \hline
$k_{\peak}$      			&95 &94 &94 &95\\ 
Peak \# DOFs  			&1910815 &1852155 &1867655 &1836490\\
End \# DOFs			&1289080 &1128070 &758425 &758425\\ \hline
$N_{\opt}$ 			&218 &220 &252 &273\\ 
$E_{\fem}$ 			&23.16 &25.56 &31.01 &32.73\\
$R_{\fem}$ 			&9.79 &10.21 &10.10 &10.32\\ \hline
\end{tabular}
\end{table*}

\section{Conclusions}\label{Conclusions}
This paper presents an adaptive mesh refinement strategy using hierarchical B-splines to perform level set topology optimization. The geometry of the design is described by a LSF and the analysis is performed using the XFEM, considering a generalized Heaviside enrichment strategy. The problem of simultaneously achieving a precise description of the geometry and an accurate evaluation of the physical responses at an acceptable computational cost is addressed by implementing a hierarchical mesh refinement, leading to an adaptive discretization of the design and state variable fields. While the state variable field interpolation is restricted to first order B-splines, the design variable field is interpolated with higher order B-splines, up to the third order. Truncated B-splines are considered for local refinement. After a user-specified number of optimization steps, local refinement, along with an increase/decrease of the minimum/maximum refinement level, is performed. The local mesh refinement is controlled by a user-defined geometric criterion. In this work, the mesh is refined around the solid/void interface and in the solid domain. To allow for a sufficient freedom in the design, we consider an initial hole seeding approach and a hole nucleation approach based on a combined level set/density scheme. The proposed method is applied to 2D and 3D structural solid-void topology optimization problems considering mass and strain energy. Optimized designs are generated on uniformly and adaptively refined meshes to provide a comparison in terms of complexity of the geometry achieved, performance of the designs and computational cost.

The numerical experiments suggest that using higher order B-splines promotes the development of smooth designs and limits the need for filtering techniques required when using low order interpolation functions. Moreover, introducing truncation allows us to accurately impose bounds on the design variables. In this study, imposing bounds on the optimization variables is required for the combined level set/\-density scheme.

Additionally, the numerical simulations show that the ability to nucleate holes during the optimization is crucial for increasing the computational advantage provided by mesh adaptation. The initial hole seeding approach limits the performance of the adaptive strategy, as rather fine meshes are required to resolve the complex initial hole pattern. Using the combined level set/density scheme, the optimization process can be initiated on coarser meshes, which leads to a more pronounced reduction of the computational cost associated with the design process.

Although the choice of the initial mesh refinement level influences the generated designs, the numerical studies demonstrate that the adaptive strategy allows for the development of optimized complex features, such as thin shear webs in 3D, even when starting from rather coarse meshes. Moreover, the adaptive designs exhibit similar performance when compared to the uniform ones. The evaluation of efficiency factors  in terms of the number of unconstrained DOFs shows that the adaptive strategy provides an improvement in terms of computational cost. In general, the efficiency factors increase both with the order of the B-spline interpolation and with the coarseness of the initial mesh used. For 3D problems, the achieved cost reduction strongly depends on the volume fraction of the optimized structures and increases as the mass of the design drops. 

Future work will focus on: (i) extending the capabilities of the mesh adaptivity, and (ii) addressing problems with more complex physics. With regards to the adaptivity, the method could be extended to handle the design and state variable fields separately. Other refinement criteria should also be considered, such as finite element error indicators. Additionally, the B-spline interpolation should be extended to higher order B-splines for the state variable field, as is the case for the design variables. Finally, as more complex physics are considered, allowing for a more dynamic refinement process could be crucial to resolve, for example, transient problems.

\section*{Acknowledgements}
The first, third, fourth, and fifth authors acknowledge the support for this work from the Defense Advanced Research Projects Agency (DARPA) under the TRADES program (agreement HR0011-17-2-0022). The opinions and conclusions presented in this paper are those of the authors and do not necessarily reflect the views of the sponsoring organizations.



\end{document}